\documentclass[final,3p,times] {elsarticle}

\usepackage{amsmath,amssymb}

\usepackage[]{amsmath}
\usepackage{graphics}
\usepackage{amssymb}
\usepackage{amsthm}
\usepackage{setspace}
\usepackage{epsfig}
\usepackage{subfigure}
\usepackage{booktabs}
\usepackage{epstopdf}
\usepackage{float}
\usepackage{mathrsfs}

\usepackage{mathrsfs}
\usepackage{cases}

\biboptions{numbers,sort&compress}
\usepackage[pdftex,colorlinks,linkcolor=red,anchorcolor=blue,citecolor=green]{hyperref}

\usepackage{amsmath}
\usepackage{times}
\usepackage{anysize}
\marginsize{2.5cm}{2.5cm}{1cm}{2cm}

\selectfont

\newtheorem{remark}{Remark}
\newtheorem{pro}{{\bf Proposition}}[section]
\journal {XXX}
\begin{document}

\begin{frontmatter}

\title{Numerical inverse scattering transform for the derivative nonlinear $\rm Schr\ddot{o}dinger$ equation}
\author[label1]{ Shikun Cui}
\author[label2]{ Zhen Wang \corref{cor1}}\ead{wangzmath@163.com}
\cortext[cor1]{School of Mathematical Sciences, Beihang University, Beijing, 100191, China}
\address[label1]{School of Mathematical Sciences, Dalian University of Technology, Dalian, 116024, China}
\address[label2]{School of Mathematical Sciences, Beihang University, Beijing, 100191, China}

\begin{abstract}

In this paper, we develop the numerical inverse scattering transform (NIST) for solving the derivative nonlinear $\rm Schr\ddot{o}dinger$ (DNLS) equation. The key technique involves formulating a Riemann-Hilbert problem (RHP) that is associated with the initial value problem and solving it numerically. Before solving the RHP, two essential operations need to be carried out. Firstly, high-precision numerical calculations are performed on the scattering data. Secondly, the RHP is deformed using the Deift-Zhou nonlinear steepest descent method.
The DNLS equation has a continuous spectrum consisting of the real and imaginary axes and features three saddle points, which introduces complexity not encountered in previous NIST approaches. In our numerical inverse scattering method, we divide the $(x,t)$-plane into three regions and propose specific deformations for each region. These strategies not only help reduce computational costs but also minimize errors in the calculations.
Unlike traditional numerical methods, the NIST does not rely on time-stepping to compute the solution. Instead, it directly solves the associated Riemann-Hilbert problem. This unique characteristic of the NIST eliminates convergence issues typically encountered in other numerical approaches and proves to be more effective, especially for long-time simulations.
\end{abstract}

\begin{keyword}
Inverse scattering transform, Riemann-Hilbert problem, Derivative nonlinear $\rm Schr\ddot{o}dinger$ equation, Numerical method
\end{keyword}
\end{frontmatter}
\tableofcontents

\section{Introduction}
The nonlinear Schr$\rm \ddot{o}$dinger equation (NLS) is a significant integrable equation that finds applications in describing the behavior of optical solitons, Bose-Einstein condensation, and other physical phenomena\cite{HSH,HJS,Kodama,deOliveira}. The derivative nonlinear Schr$\rm \ddot{o}$dinger (DNLS) equation is an important variant of the NLS equation.
Within the class of DNLS equations, three well-known equations are the Kaup-Newell (KN) equation\cite{KN}
\begin{equation}\label{DNLS1}
iq_t+q_{xx}+i(|q|^2q)_x=0,
\end{equation}
the Chen-Lee-Liu (CLL) equation\cite{CLL}
\begin{equation}\label{DNLS2}
iq_t+q_{xx}+i|q|^2q_x=0,
\end{equation}
and the Gerdzhikov-Ivanov (GI) equation\cite{GI}
\begin{equation}\label{DNLS3}
iq_t+q_{xx}-iq^2\bar{q}_x+\frac{1}{2}|q|^4q=0.
\end{equation}. These equations have their own unique properties and are extensively studied in various fields.


We refer to the KN equation, the CLL equation, and the GI equation as DNLSE $\rm \uppercase\expandafter{\romannumeral1}$, DNLSE $\rm \uppercase\expandafter{\romannumeral2}$, and DNLSE $\rm \uppercase\expandafter{\romannumeral3}$, respectively. These equations are interconnected through gauge transformations\cite{gauge_Kundu,gauge_Wadati,gauge_Clarkson}, allowing for conversion between them. Therefore, our main focus lies in investigating the specific properties and characteristics of DNLSE $\rm \uppercase\expandafter{\romannumeral3}$.


The inverse scattering transform (IST), originally proposed by GGKM, is a significant approach for solving integrable equations with initial value problems\cite{ist_Gardner}.
In the case of the DNLS equation, the inverse scattering transform was introduced by Kaup and Newell\cite{KN}.
A more contemporary method that has gained widespread popularity is the Riemann-Hilbert method, which serves as a modern version of the inverse scattering transform.
By employing the inverse scattering transform, Zhou and Huang were able to obtain $N$-soliton solutions for the DNLS equation\cite{ist_zhouhuang}.
Additionally, Chen and Lam extended the application of the inverse scattering transform to the DNLS equation with non-zero boundary conditions\cite{ist_chenlam}.
Utilizing the Riemann-Hilbert method, Guo and Ling derived the $N$-soliton solution for the coupled derivative Schr$\rm \ddot{o}$dinger equation\cite{ist_gbl}.
Zhang and Yan further developed the inverse scattering transform for the DNLS equation, considering both zero and non-zero boundary conditions\cite{IST_yan}.
Ma and Kuang expanded the scope of the inverse scattering transform by applying it to the nonlocal DNLS equation\cite{ist_mxx}.
Lastly, Liu et al. successfully obtained triple-pole soliton solutions of the DNLS equation using the inverse scattering transform\cite{ist_liuetal}.


The inverse scattering transform is not only applicable in solving the DNLS equation but also in exploring the existence of solutions\cite{ist_2_Liujq,ist_2_Pelinovsky0,ist_2_Pelinovsky,ist_2_Jenkins}.
Furthermore, it plays a crucial role in studying the long-term behavior of the solution.
Deift and Zhou introduced a nonlinear steepest descent method for determining the long-term behavior of integrable equations\cite{deift1}.
Miller and Bilman proposed the robust inverse scattering transform to analyze the asymptotic properties of rogue waves\cite{Miller}.
Xu et al. investigated the long-term asymptotics of the DNLS equation with step-like initial values\cite{ist_3_xujian}.
Arruda et al. studied the long-term asymptotics of the DNLS equation on the half-line\cite{ist_3_Arruda}.
Jenkis et al. presented the soliton resolution for the derivative nonlinear Schr$\rm \ddot{o}$dinger equation\cite{jenkis}.
Tian examined the long-term asymptotic behavior of the generalized derivative nonlinear Schr$\rm \ddot{o}$dinger equation\cite{ist_3_Tian}.
Yang and Fan analyzed the long-term asymptotic behavior of the DNLS equation with initial data of finite density type\cite{ist_3_yang}.


Methods for solving the DNLS equation can be classified into two categories: analytical methods\cite{KN,ist_zhouhuang,ist_chenlam,ist_gbl,IST_yan,ist_mxx,ist_liuetal,DB_fan,DB_gbl,CP1,yangbo1} and numerical methods\cite{Boyd,Yangjk,JieS}. Analytical methods, such as the inverse scattering transform, are primarily useful for calculating pure soliton solutions of the DNLS equation. However, they fail to provide accurate results when the reflection coefficient is non-zero, limiting their applicability to analyzing the long-term asymptotic behavior of the solution\cite{deift1,ist_3_Arruda,ist_3_xujian,jenkis,ist_3_Tian,ist_3_yang}.
In contrast, numerical methods offer a viable framework for solving general initial value problems in the DNLS equation. Traditional numerical methods rely on time-stepping to compute the solution, which accumulates calculation errors over time and often leads to convergence issues. Hence, it becomes crucial to develop a reliable numerical method capable of accurately calculating the long-term evolution of the solution.


Recently, there has been a notable development in the field of solving integrable equations with the introduction of the numerical inverse scattering transform (NIST) by Deconinck, Olver, and Trogdon\cite{nist_Trogdon}. This powerful numerical method has been successfully applied to various equations such as the focusing and defocusing nonlinear Schr$\rm\ddot{o}$dinger equations by Trogdon and Olver\cite{nist_Trogdon1}, the Toda lattice by Bilman and Trogdon\cite{nist_Trogdon2}, the sine-Gordon equation by Deconinck, Trogdon, and Yang\cite{nist_Trogdon3}, and the Kundu-Eckhaus equation by Cui and Wang \cite{nist_cui}.
The NIST can be divided into two main components: numerical direct scattering and numerical inverse scattering. In the numerical direct scattering step, numerical methods are employed to compute the scattering data. Subsequently, in the numerical inverse scattering step, the calculated scattering data is utilized to recover the solutions for the DNLS equation.
The NIST offers distinctive characteristics, which can be summarized as follows:
\begin{itemize}
  \item Unlike traditional numerical methods, the NIST does not require time-stepping, thus making the computational cost independent of time and space.
  \item The error generated by the NIST remains unaffected by increasing time, thereby avoiding numerical convergence issues.
  \item While traditional numerical methods are more effective for short-time simulations, the NIST proves to be more effective for long-time simulations.
  \item It is important to note that currently, the NIST is primarily applicable to integrable equations, and its implementation is complex.
\end{itemize}

The NIST represents a significant advancement in numerical methods for integrable equations, demonstrating its advantages such as independence from time-stepping, error stability, and increased effectiveness in long-time simulations. However, it is worth mentioning that its application is confined to integrable equations, and its implementation can be challenging.

In this paper,  the numerical inverse scattering transform of the DNLS equation is developed.
In numerical direct scattering, we need to solve a quadratic eigenvalue problem rather than a standard eigenvalue problem, which is different from the AKNS system. Introducing $\widehat{\psi}=k\psi$, we standardize the quadratic eigenvalue problem, then the high-precision numerical method is used to solve the standard eigenvalue problem. The continuous spectrum of the DNLS equation is $i\mathbb{R}\cup\mathbb{R}$ rather than $\mathbb{R}$, which has never appeared in the previous NIST.
Considering that the original Riemann-Hilbert problem is oscillatory, we must deform it.
There are three saddle points for the DNLS equation, so deforming the original Riemann-Hilbert problem is complex and difficult. Based on the relationship between $x$ and $t$, we divide the $(x,t)$-plane into three regions and propose three schemes for deforming the original Riemann-Hilbert problem. These deformations are proposed after considering the computational cost and accuracy comprehensively.
At the same time, the rationality of the deformations are demonstrated.


Our paper is structured as follows:

$\bullet$ Section \ref{sec:method intruction}: Construction of the Riemann-Hilbert Problem

We review the key findings related to the construction process of the Riemann-Hilbert problem.

$\bullet$ Section \ref{sec_direct}: Numerical Direct Scattering

We implement numerical direct scattering methods.
Through these techniques, we calculate essential values such as the reflection coefficient, discrete eigenvalues, and normalized constant for the DNLS equation.

$\bullet$ Section \ref{sec_inverse}: Numerical Inverse Scattering

We employ numerical inverse scattering methods.
In section \ref{deform_1_1} and section \ref{deform_1_2}, we focus on deforming the Riemann-Hilbert problem, assuming no soliton evolution from the initial value of the DNLS equation.
In section \ref{deform_2_1} and section \ref{deform_2_2}, we examine the deformation of the Riemann-Hilbert problem, considering scenarios where the initial value of the DNLS equation can evolve into solitons.
In section \ref{sec_cheby}, we utilize the Chebyshev collocation method to solve the deformed Riemann-Hilbert problem.

$\bullet$ Section \ref{sec_result1}: Numerical Results and Comparison

We present the numerical results obtained using the Numerical Inverse Scattering Transform (NIST).
Additionally, we compare the performance of the NIST with the Fourier spectral method (FSM).

$\bullet$ Section \ref{sec_con}: Conclusions and Discussions

We summarize the conclusions derived from our study.
Furthermore, we engage in discussions regarding the performance and effectiveness of the NIST in comparison to the FSM.
We also propose potential avenues for future research, and emphasize the contributions our study makes.

\section{Construction of the Riemann-Hilbert Problem}\label{sec:method intruction}
\subsection{Spectral analysis}
Considering the initial value problem for the DNLSE $\rm \uppercase\expandafter{\romannumeral3}$,
\begin{equation}\label{initial_problem}
\left\{
\begin{array}{lr}
iq_t+q_{xx}-iq^2\bar{q}_x+\frac{1}{2}|q|^4q=0, \\
q(x,0)=q_0(x),
\end{array}
\right.
\end{equation}
where $q_0(x)$ is a given function defined in Schwartz space $\mathcal{S}(\mathbb{R})$, subscripts $x$ and $t$ represent the partial derivatives with respect to space and time respectively, and $\bar{q}$ is the complex conjugate of $q$.

The DNLS $\rm \uppercase\expandafter{\romannumeral3}$  is completely integrable and admits the following Lax pair\cite{ist_2_Pelinovsky},
\begin{equation}\label{lax_1}
\left\{
\begin{array}{lr}
\psi_x + ik^{2}\sigma_{3}\psi=U\psi, \vspace{0.15cm} \\
\psi_t + ik^{4}\sigma_{3}\psi=V\psi,
\end{array}
\right.
\end{equation}
where $k\in\mathbb{C}$ is the spectral parameter, $\psi$ is the eigenfunction, $$U=kQ+\frac{i}{2}|q|^2\sigma_3, \
V=2k^3Q + ik^2|q|^2\sigma_{3}+\frac{i}{4}|u|^4\sigma_3-i kQ_x\sigma_3+\frac{1}{2}(q\bar{q}_x-q_x\bar{q})\sigma_3,\
Q=\left(
\begin{array}{cc}
0 & q \\
-\bar{q} & 0 \\
\end{array}
\right).
$$
Here $\sigma_3$ is one of the Pauli matrices
$$\sigma_1=\left(
\begin{array}{cc}
0 & 1 \\
1 & 0 \\
\end{array}
\right), \sigma_2=\left(
\begin{array}{cc}
0 & i \\
-i & 0 \\
\end{array}
\right)\ {\rm and}\ \sigma_3=\left(
\begin{array}{cc}
1 & 0 \\
0 & -1 \\
\end{array}
\right).$$

The eigenfunction $\psi$ admits the following asymptotics
$$ \psi\rightarrow{\rm e}^{-{i}\theta\sigma_3}\ {\rm as} \ |x|\rightarrow\infty, {\rm where}\ \theta=k^{2 }x+2k^{4}t.$$

Doing transformation $\mu=\psi{\rm e}^{i\theta\sigma_{3}}$,
we get a new Lax pair
\begin{equation}\label{lax_2}
\left\{
\begin{array}{lr}
\mu_x + {i} k^{2}[\sigma_{3}, \mu]=U\mu,
\vspace{0.15cm}
\\ \mu_t + 2{i} k^{4}[\sigma_{\rm3}, \mu]=V\mu,
\end{array}
\right.
\end{equation}
where $[\sigma_3, \mu]=\sigma_3\mu-\mu\sigma_3$.

The eigenfunction $\mu$ admits the following asymptotics $$\mu\rightarrow I\ {\rm as}\ |x|\rightarrow\infty,$$ where $I$ is the identity matrix.
\begin{pro}
The solution of the {\rm DNLSE} $\rm \uppercase\expandafter{\romannumeral3}$ is given by
\begin{equation}\label{coeffs_6}
q=2{i}\mu^{(1)}_{(1,2)},
\end{equation}
where  $\mu^{(1)}_{(1,2)}$ represents the first row and second column element of matrix $\mu^{(1)}$, $\mu^{(1)}
=\lim_{k\rightarrow\infty}(k\mu).$
\end{pro}
\begin{proof}
Expand $\mu$ in $k=\infty$, and we get
\begin{equation}\label{eqqqq}
\mu=\mu^{(0)}+\frac{\mu^{(1)}}{k}+\frac{\mu^{(2)}}{k^2}+\cdots
\end{equation}
Substitute  Eq.(\ref{eqqqq}) into Eq.(\ref{lax_2}), we can prove Eq.(\ref{coeffs_6}) by matting  the $O(k^0)$, $O(k^1)$ and $O(k^2)$ terms.

\end{proof}

The Lax pair (\ref{lax_2}) can be written into a total differential form,
$${\rm d}({\rm e}^{i\theta\hat{\sigma}_3}\mu)={\rm e}^{i\theta\hat{\sigma}_3}(U{\rm d}x+V{\rm d}t)\mu,$$
where ${\rm e}^{i\theta\hat{\sigma}_3}\mu={\rm e}^{i\theta\sigma_3}\mu {\rm e}^{-{i}\theta\sigma_3}.$

The eigenfunctions of the Lax pair (\ref{lax_2}) are obtained, which satisfy the Volterra integral equation,
\begin{equation}\label{lax_3}
\begin{array}{lr}
\mu_{1}(x,t,k)=I+ \int_{-\infty}^x{\rm e}^{{-i} k^2(x-y)\hat{\sigma}_{3}}U(y,t,k)\mu_{1}(y,t,k) dy, \\
\mu_{2}(x,t,k)=I+ \int_{+\infty}^x{\rm e}^{{-i} k^2(x-y)\hat{\sigma}_{3}}U(y,t,k)\mu_{2}(y,t,k) dy.
\end{array}
\end{equation}


The eigenfunctions of Eq. (\ref{lax_3}) exhibit specific symmetries and analyticity properties. These characteristics play a crucial role in constructing the Riemann-Hilbert problem associated with the initial value problem described by equation (\ref{initial_problem}).
These symmetries and analyticity properties allow for the formulation of a well-posed mathematical problem known as the Riemann-Hilbert problem. The Riemann-Hilbert problem involves finding a function that satisfies certain conditions on the boundary of a given domain.
By leveraging the symmetries and analyticity of the eigenfunctions, we can successfully construct the Riemann-Hilbert problem and subsequently analyze the associated initial value problem (\ref{initial_problem}). This enables a deeper understanding of the dynamics and behavior of the system under investigation.

\begin{figure}[H]
  \centering
  \includegraphics[width=2.35in,height=2.05in]{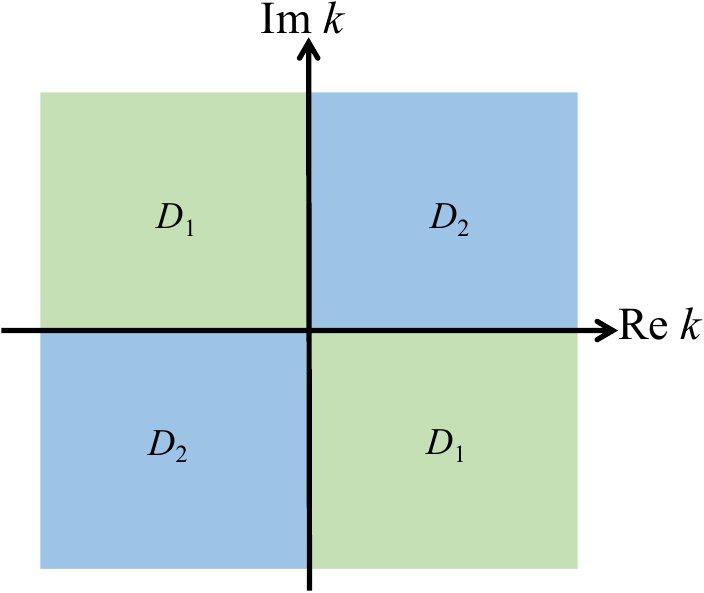}\\
  \caption{The area $D_1$ and area $D_2$ on complex-$k$ plane.}\label{pic:area1}
\end{figure}

\begin{pro}[Analyticity]
The eigenfunction $\mu_j(x,t,k)$ $(j=1,2)$ satisfies the following analyticity,
\begin{itemize}
\vspace{-0.0cm}
\item $[\mu_1]_2$ and $[\mu_2]_1$ are bounded and analytic in area $D_1$=$\{k\in\mathbb{C}|\  {\rm Im}(k^2)<0 \},$
\vspace{-0.0cm}
\item $[\mu_1]_1$ and $[\mu_2]_2$ are bounded and analytic in area $D_2$=$\{k\in\mathbb{C}|\  {\rm Im}(k^2)>0 \},$
\vspace{-0.0cm}
\item ${\rm det}[\mu_{1}(x,t,k)] = {\rm det} [\mu_{2}(x,t,k)] =1,$
\end{itemize}
where $[\mu_j]_k$ is the $k$-th column of $\mu_j(x,t,k)$. The area $D_1$ and the area $D_2$ are shown in Figure \ref{pic:area1}.
\end{pro}

\begin{pro}[Symmetry]
The eigenfunction $\mu_j(x,t,k)$ $(j=1,2)$ satisfies the following symmetries,\begin{itemize}
  \vspace{-0.0cm}
  \item $-\sigma \overline{\mu_j(x,t,\bar{k})} \sigma=\mu_j$,
  \item $\sigma_3 \mu_j(x,t,-k) \sigma_3=\mu_j$,
\end{itemize}
where $\sigma$= $\left(
                   \begin{array}{cc}
                     0 & 1 \\
                     -1 & 0 \\
                   \end{array}
                 \right)$,
                 $\sigma_{\rm3}$= $\left(
                   \begin{array}{cc}
                     1 & 0 \\
                     0 & -1 \\
                   \end{array}
                 \right)$.
\end{pro}

\begin{proof}
Considering the $x$-part of the Lax pair (\ref{lax_2})
\begin{equation}\label{proof_1}
(\ \mu_{j}(k)\ )_x+ik^{2}[\sigma_{3},\mu_j(k) ]=(kQ+q_1)\mu_j(k),
\end{equation}
where $q_1=\begin{cases}
0\sigma_3,& {\rm for\ the\ DNLS\  \uppercase\expandafter{\romannumeral1}},\\
-\frac{i}{4}|q|^{2}\sigma_{3},& {\rm for\ the\ DNLS\  \uppercase\expandafter{\romannumeral2}}, \\
\frac{i}{2}|q|^{2}\sigma_{3},& {\rm for\ the\ DNLS\  \uppercase\expandafter{\romannumeral3}}.
\end{cases}$

When we replace $k$ with $\bar{k}$ and perform a conjugate transformation of Eq. (\ref{proof_1}), the new equation becomes:
\begin{equation}\label{proof_2}
(\ \overline{\mu_{j}(\bar{k})}\ )_x-{i} k^{2}[\sigma_{3},\overline{\mu_j(\bar{k})} ]=(k\overline{Q}-q_1)\overline{\mu_j(\bar{k})}.
\end{equation}

Multiply $\sigma$ on both sides of Eq. (\ref{proof_2}),
and we get
\begin{equation}\label{proof_3}
(\ \sigma\overline{\mu_{j}(\bar{k})}\sigma \ )_x+{i} k^{2}[\sigma_{3},\sigma\overline{\mu_j(\bar{k})}\sigma]=(kQ+q_1)\sigma\overline{\mu_j(\bar{k})}\sigma.
\end{equation}

Both $ -\sigma\overline{\mu_{j}(\bar{k})}\sigma$ and  $\mu_{j}(k)$ are the solution of Eq. (\ref{proof_1}), and they have the same asymptotics,
$${\rm so}\quad -\sigma\overline{\mu_{j}(x,t,\bar{k})}\sigma=\mu_{j}(x,t,k),\ (j=1, 2).$$

Replace $k$ with $-k$ in Eq. (\ref{proof_1}), and we get
\begin{equation}\label{dingli_1}
(\ \mu_{j}(-k)\ )_x+ik^{2}[\sigma_{3},\mu_j(-k) ]=(-kQ+q_1)\mu_j(-k).
\end{equation}

Multiplying $\sigma_3$ on both sides of Eq. (\ref{dingli_1}), we get
$$(\sigma_{3}\ \mu_{j}(-k)\sigma_{3}\ )_x+ik^{2}[\sigma_{3},\sigma_{3}\mu_j(-k)\sigma_{3} ]=(kQ+q_1)\ \sigma_{3}\mu_j(-k)\sigma_{3}.$$

Both $ \sigma_{3}\mu_{j}(-k)\sigma_{3}$ and  $\mu_{j}(k)$ are the solution of Eq. (\ref{proof_1}), and they have the same asymptotics
$$ \sigma_{3}\mu_{j}(-k)\sigma_{3} = \mu_{j}(k),\ (j=1, 2).$$

\end{proof}

\vspace{0.2cm}
Both $\mu_1(x,t,k)$ and $\mu_2(x,t,k)$ are the solution of the Lax pair (\ref{lax_2}), so there exists a linear relation between  $\mu_1(x,t,k)$ and $\mu_2(x,t,k)$,
\begin{equation}\label{matrix_1}
\mu_2(x,t,k)=\mu_{1}(x,t,k) {\rm e}^{-i\theta\hat{\sigma}_{3}} S(k),
\end{equation}
where $S(k)$ is independent of $t$ and $x$.

\begin{pro}
\ $S(k)$ satisfies the following symmetries: $S(k)=-\sigma\overline{S(\bar{k})}\sigma$, $S(k)=\sigma_{\rm3}S(-k)\sigma_{3},$
$$S(k)=\left(
         \begin{array}{cc}
           \tilde{a}(k) & b(k) \\
           \tilde{b}(k) & a(k)  \\
         \end{array}
       \right),
$$
$$\tilde{a}(k)=\overline{a(\bar{k})},$$
$$\tilde{b}(k)=-\overline{b(\bar{k})}.$$
\end{pro}

\begin{proof}
With the help of $\mathbf{Proposition\ 2.3}$, the symmetries of $S(k)$ can be proved easily.

\end{proof}

\subsection{A Riemann-Hilbert problem associated with initial value problem}
Rewriting Eq. (\ref{matrix_1}) into
\begin{equation}\label{rhp_1}
( [\mu_2]_1, [\mu_2]_2 )=( [\mu_1]_1, [\mu_1]_2 ){\rm e}^{-i\theta\hat{\sigma}_3}\left(
         \begin{array}{cc}
           \hat{a}(k) & b(k) \\
           \hat{b}(k) & a(k) \\
         \end{array}
       \right),
\end{equation}
where $[\mu_1]_2$ and $[\mu_2]_1$ are bounded and analytic in $D_1$, $[\mu_1]_1$ and $[\mu_2]_2$ are bounded and analytic in $D_2$.
The area $D_1$ and $D_2$ are shown in Figure \ref{pic:area1}.

We define $M(x,t,k)$ by $$M(x,t,k)=\left\{
\begin{array}{lr}
\Big( \frac{[\mu_2]_1}{\tilde{a}(k)},[\mu_1]_2, \Big) ,\quad  k\in D_1, \\
\Big( [\mu_1]_1 ,  \frac{[\mu_2]_2}{a(k)} \Big) ,\quad k\in D_2,
\end{array}
\right.$$
then we get a Riemann-Hilbert problem (RHP):

\noindent$\textbf{Riemann-Hilbert problem \uppercase\expandafter{\romannumeral0}}$
\begin{itemize}
  \vspace{-0.0cm}
  \item $M(x,t,k)$ is meromorphic in $\mathbb{C}\setminus\Sigma_1$,
  \vspace{-0.0cm}
  \item $M_{+}(x,t,k)=M_{-}(x,t,k)T(x,t,k)$, $k\in\Sigma_1$,
  \vspace{-0.0cm}
  \item $M(x,t,k)\rightarrow I$, $k\rightarrow\infty$,
\end{itemize}
where $\Sigma_1=i\mathbb{R}\cup \mathbb{R}$ , the jump matrix $$T(x,t,k)=\left(
\begin{array}{cc}
 1& {r}(k){\rm e}^{-2i\theta(k)} \\
  -\tilde{r}(k){\rm e}^{2i\theta(k)} & 1-r(k)\tilde{r}(k)  \\
  \end{array}
  \right),$$
where $r(k)=\frac{b(k)}{a(k)}$ and $\tilde{r}(k)=\frac{\tilde{b}(k)}{\tilde{a}(k)}$.

\begin{remark}
If $a(k)$ exists the one-order zero point $\kappa_j$ in the first quadrant of $\mathbb{C}$,
\textbf{Riemann-Hilbert problem} will satisfy the pole conditions,
\begin{itemize}
  \item${\rm Res}\{M(x,t,k), k=\kappa_j\} =\lim_{k\rightarrow\kappa_j}M(x,t,k)\left(
                                                                                  \begin{array}{cc}
                                                                                    0 & C_j{\rm e}^{-2i\theta(\kappa_j)}  \\
                                                                                    0& 0 \\
                                                                                  \end{array}
                                                                                \right),$
  \item {\rm Res}\{M(x,t,k), $k=\bar{\kappa}_j$\} = $\lim_{k\rightarrow\bar{\kappa}_j}M(x,t,k)\left(
                                                                                  \begin{array}{cc}
                                                                                    0 & 0\\
                                                                                    \bar{C}_j{\rm e}^{2{ i}\theta(\bar{\kappa}_j)} & 0 \\
                                                                                  \end{array}
                                                                                \right),$
  \item {\rm Res}\{M(x,t,k), $k=-\bar{\kappa}_j$\} = $\lim_{k\rightarrow-\bar{\kappa}_j}M(x,t,k)\left(
                                                                                  \begin{array}{cc}
                                                                                    0 & 0 \\
                                                                                    \bar{C}_j{\rm e}^{2 i\theta(-\bar{\kappa}_j)} & 0 \\
                                                                                  \end{array}
                                                                                \right),$
  \item {\rm Res}\{M(x,t,k), $k=-\kappa_j$\} = $\lim_{k\rightarrow-\kappa_j}M(x,t,k)\left(
                                                                                  \begin{array}{cc}
                                                                                    0 & C_j{\rm e}^{-2 i\theta(-\kappa_j)} \\
                                                                                   0 & 0 \\
                                                                                  \end{array}
                                                                                \right).$
  \end{itemize}

\end{remark}
\begin{pro}
The solution of the $\rm DNLSE\ \uppercase\expandafter{\romannumeral1}$ $\rm (\ref{initial_problem})$ is given by the $\textbf{Riemann-Hilbert problem \uppercase\expandafter{\romannumeral0}}$
$$q(x,t)=2im_{(12)},$$
where $m_{(12)}$ represents the first row and second column element of matrix $m$, $m=\lim_{k\rightarrow\infty} \Big( kM(x,t,k) \Big)$.
\end{pro}

\section{Numerical direct scattering}\label{sec_direct}
In this section, the numerical methods are used to calculate the scattering data.
The scattering data consists of the reflection coefficient $r(k)$, the discrete eigenvalue $\kappa_j(j=1,2,...,n)$ and the normalized constant $C_j$.
It is important to note that the condition $a(\kappa_j)=0$ is equivalent to $\kappa_j$ being a discrete eigenvalue of the DNLS equation.
If $k=\kappa_j$ is the one-order pole of $a(k)$, we define the normalized constant by
$$C_j=\frac{b_j}{a'(\kappa_j)},$$ where $a'(k)$ is the derivative of $a(k)$ and $b_j$ is a constant.
Set $k=\kappa_j$, Eq. (\ref{rhp_1}) gives $$[\mu_2]_2(x,t,\kappa_j)=b_j{\rm e}^{-2i\theta(\kappa_j)} [\mu_1]_1(x,t,\kappa_j).$$

\subsection{Calculating the discrete eigenvalues}\label{sec_eigenvalue}

There are two main methods available for calculating the discrete eigenvalues of the DNLS equation\cite{Yousefi}. The first method involves an iterative approach to find the zero points of $a(k)$, while the second method revolves around solving a matrix eigenvalue problem.
When $\kappa_j$ is identified as a discrete eigenvalue for the DNLS equation, $\mathbf{Proposition\ 2.3}$ asserts that $\bar{\kappa}_j$, $-\bar{\kappa}_j$, and $-\kappa_j$ also qualify as discrete eigenvalues. This property simplifies the search for zero points of a(k) to the first quadrant of the complex plane.
However, even limiting the search to the first quadrant of the complex-$k$ plane can still be complex and computationally expensive when using the iterative method. Therefore, the more efficient approach of solving the matrix eigenvalue problem is typically employed to compute the eigenvalues for the DNLS equation.

The eigenvalues for the DNLS equation are calculated from the following eigenvalue problems
\begin{numcases}{}
\psi_x + ik^{2}\sigma_{3}\psi=kQ\psi,\label{dnls_lax_1} \quad &{\rm for\ the\ DNLSE\  \uppercase\expandafter{\romannumeral1}}, \\
\psi_x + ik^{2}\sigma_{3}\psi=kQ\psi+\frac{i}{4}Q^2\sigma_3\psi, \quad &{\rm for\ the\ DNLSE\  \uppercase\expandafter{\romannumeral2}},\label{dnls_lax_2}\\
\psi_x + ik^{2}\sigma_{3}\psi=kQ\psi-\frac{i}{2}Q^2\sigma_3\psi, \quad &{\rm for\ the\ DNLSE\  \uppercase\expandafter{\romannumeral3}} \label{dnls_lax_3}.
\end{numcases}
Eqs. (\ref{dnls_lax_1}), (\ref{dnls_lax_2}) and (\ref{dnls_lax_3}) can be classified as quadratic eigenvalue problems, which we refer to as DNLSE eigenvalue problems. In order to solve these problems more effectively, we must convert the DNLSE eigenvalue problems into standard eigenvalue problems.

\begin{pro}\label{test}
The DNLSE eigenvalue problems {\rm (\ref{dnls_lax_1}), (\ref{dnls_lax_2}) and (\ref{dnls_lax_3})} can be written as the standard eigenvalue problems by introducing transformation $\widehat{\psi}=k\psi$.
\end{pro}
\begin{proof}
We only take the eigenvalue problem (\ref{dnls_lax_1}) as an example, the eigenvalue problems (\ref{dnls_lax_2}) and (\ref{dnls_lax_3}) are similar to (\ref{dnls_lax_1}).
Set $\widehat{\psi}=k\psi$, we get
\begin{equation}
\begin{aligned}\label{linear_1}
&\psi_x + {i} k\sigma_{3}\widehat{\psi}=Q\widehat{\psi}, \\
&\widehat{\psi}=k\psi,
\end{aligned}
\end{equation}
where $\psi$ and $\widehat{\psi}$ are two vectors of size $2\times1$.

Eq. (\ref{linear_1}) is rewritten into
\begin{equation}\label{linear_21}
\left(
  \begin{array}{cccc}
    0        & -iq  & i\partial_x & 0 \\
   -i\bar{q} & 0      & 0 & -i\partial_x \\
    1 & 0 & 0 & 0 \\
    0 & 1 & 0 & 0 \\
  \end{array}
\right)\left(
               \begin{array}{c}
                 \widehat{\psi} \\
                 \psi \\
               \end{array}
             \right)=k\left(
               \begin{array}{c}
                 \widehat{\psi} \\
                 \psi \\
               \end{array}
             \right),
\end{equation}
where $q$ is the potential function, $\partial_x$ represents the differential operation, $\widehat{\psi}$ and $\psi$ are the eigenfunctions.
It is obvious that Eq. (\ref{linear_21}) is a standard eigenvalue problem.

\end{proof}

\textbf{Proposition \ref{test}} gives a method to standardize the DNLSE eigenvalue problems.
The $\rm DNLSE\ {\uppercase\expandafter{\romannumeral1}}$ eigenvalue problem is transformed into Eq. (\ref{linear_21}).
The $\rm DNLSE\ {\uppercase\expandafter{\romannumeral2}}$ eigenvalue problem is transformed into
\begin{equation}\label{linear_22}
\left(
  \begin{array}{cccc}
    0        & -{i}q  & i\partial_x-\frac{1}{4}|q|^2 & 0 \\
   -i\bar{q} & 0      & 0 & -i\partial_x-\frac{1}{4}|q|^2 \\
    1 & 0 & 0 & 0 \\
    0 & 1 & 0 & 0 \\
  \end{array}
\right)\left(
               \begin{array}{c}
                 \widehat{\psi} \\
                 \psi \\
               \end{array}
             \right)=k\left(
               \begin{array}{c}
                 \widehat{\psi} \\
                 \psi \\
               \end{array}
             \right).
\end{equation}
The $\rm DNLSE\ {\uppercase\expandafter{\romannumeral3}}$ eigenvalue problem is transformed into
\begin{equation}\label{linear_23}
\left(
  \begin{array}{cccc}
    0        & -iq  & i\partial_x+\frac{1}{2}|q|^2 & 0 \\
   -i\bar{q} & 0    & 0 & -i\partial_x+\frac{1}{2}|q|^2 \\
    1 & 0 & 0 & 0 \\
    0 & 1 & 0 & 0 \\
  \end{array}
\right)\left(
               \begin{array}{c}
                 \widehat{\psi} \\
                 \psi \\
               \end{array}
             \right)=k\left(
               \begin{array}{c}
                 \widehat{\psi} \\
                 \psi \\
               \end{array}
             \right).
\end{equation}

By transforming the quadratic eigenvalue problems into standard eigenvalue problems, we simplify the solution process and can leverage existing numerical methods for solving standard eigenvalue problems to find the eigenvalues and eigenvectors of the DNLSE equations.
There are some methods can be used to solve these eigenvalue problems\cite{Boyd,Yangjk,Trogdon111,Cui}, we use the Chebyshev collocation method\cite{Cui} to solve the DNLSE eigenvalue problems.
Before introducing the numerical method, some marks are introduced in Table \ref{table_111}.
\begin{table}[h]
\centering
\caption{Annotations of some marks used in this paper.}
\begin{tabular}{c|c|c}
  \hline\hline\noalign{\smallskip}
  Mark& Formula of the mark& Note \\
  \hline
    $\vec{\chi}$ & $\vec{\chi}=\bigg(-1, {\rm cos}(\frac{n_1-2}{n_1-1}\pi),\ \cdots \, \ {\rm cos}(\frac{1}{n_1-1}\pi), 1 \bigg)^{\top}$& Chebyshev nodes \\
    $n_1$&& The number of Chebyshev nodes\\
  $T_k(x)$      &$T_k={\rm cos}(k\cdot {\rm arccos}x)$, $k=0,\ldots, n_1-1$& Chebyshev polynomials\\
  $T(x)$&$T(x)=[T_0(x),T_1(x),\ldots,T_{n_1-1}(x)]$&\\
  $H(x)$       & $H(x)={\rm tanh}(\hat{a}x)$, $0<\hat{a}<1$   & Hyperbolic tangent map\\
  $H^{-1}(x)$ & &The inverse mapping of $H(x)$ \\
  $\mathcal{D}$& &Chebyshev derivative matrix \\
  $f(\vec{\chi})$ &$f(\vec{\chi})=\left(
                     \begin{array}{c}
                       f(x)\mid_{x=-1} \\
                       f(x)\mid_{x={\rm cos}(\frac{n_{1}-2}{n_{1}-1}\pi)} \\
                       \vdots \\
                       f(x)\mid_{x=1} \\
                     \end{array}
                   \right)$ & \\
   $\mathcal{F}$ &The discrete cosine transform of $f(\vec{\chi})$& Chebyshev coefficient matrix \\
  $\mathrm{D}\Big[f(\vec{\chi})\Big]$&${\rm D}\Big[f(\vec{\chi})\Big]=\left(
                                                                        \begin{array}{ccc}
                                                                          f(x)\mid_{x=-1}  &  &  \\
                                                                             & \ddots &  \\
                                                                           &  &f(x)\mid_{x=1}  \\
                                                                        \end{array}
                                                                      \right)$& Diagonal matrix\\
  \hline
\end{tabular}\label{table_111}
\end{table}

We need to transform the DNLSE eigenvalue problems into matrix eigenvalue problems.
For a given function $g(x)$ defined in real field $\mathbb{R}$,  we can approximate $g(x)$ and $\frac{\partial g}{\partial x}$ by
\begin{numcases}{}
  g(x)=T^{\mathbb{R}}(x)\mathcal{F}g(\hat{\chi}), x\in\mathbb{R}, \label{chebyshev_1} \\
 \frac{\partial g(x)}{\partial x}=\frac{\partial H(x)}{\partial x}T^{\mathbb{R}}(x)\mathcal{D}\mathcal{F}g(\hat{\chi}), {\rm for} \ x\in\mathbb{R},
 \label{chebyshev_2}
\end{numcases}
where $T^{\mathbb{R}}(x)=T( H(x) )=[T_0( H(x) ), \cdots, T_{n_1-1}( H(x) )]$ and $\hat{\chi}=H^{-1}(\vec{\chi})$.

Using Eq. (\ref{chebyshev_1}) and Eq. (\ref{chebyshev_2}), $\psi=\left(\begin{array}{c}
                                                   \psi_1 \\
                                                   \psi_2 \\
                                                 \end{array}\right)$,
$\widehat{\psi}=\left(\begin{array}{c}
                                                   \widehat{\psi_1} \\
                                                   \widehat{\psi_2} \\
                                                 \end{array}\right)$, $\psi_x$ and the potential function  $q(x)$ are approximated by the Chebyshev polynomials with $n_1$ nodes,
\begin{equation}\label{cheb_1_1}
\begin{aligned}
&\psi_j=T^{\mathbb{R}}(x)\mathcal{F}\psi_j(\hat{\chi}), \\ &\widehat{\psi_j}=T^{\mathbb{R}}(x)\mathcal{F}\widehat{\psi_j}(\hat{\chi}), \\
&\frac{\partial \psi_j}{\partial x}=\frac{\partial H(x)}{\partial x}T^{\mathbb{R}}(x)\mathcal{D}\mathcal{F}\psi_j(\hat{\chi}), \\
&q(x)=T^{\mathbb{R}}(x)\mathcal{F}q(\hat{\chi}).
\end{aligned}
\end{equation}

Substitute Eq. (\ref{cheb_1_1}) into Eqs. (\ref{linear_21}), (\ref{linear_22}) and (\ref{linear_23}), and set $x=\hat{\chi}=H^{-1}(\vec{\chi})$,
then the DNLSE eigenvalue problem is transformed into a matrix eigenvalue problem
\begin{equation}\label{eigen_1}
\left[
\setlength{\arraycolsep}{0.05pt}
  \begin{array}{cccc}
 0 & i{\rm D}\Big[q(\hat{\chi})\Big]  &\quad  {\rm  D}\Big[q_1(\hat{\chi})\Big] +i{\rm  D}\Big[\frac{\partial H(\vec{\chi})}{\partial x} \Big]\mathcal{F}^{-1}\mathcal{D}\mathcal{F}
  & 0 \\
 i{\rm D}[\bar{q}(\hat{\chi})] & 0&  0 &{\rm D}\Big[q_1(\hat{\chi})\Big]-i{\rm  D}\Big[\frac{\partial H(\vec{\chi})}{\partial x} \Big]\mathcal{F}^{-1}\mathcal{D}\mathcal{F} \\
 I&0&0&0 \\
 0&I&0&0
  \end{array}
\right]\left[
         \begin{array}{c}
          \widehat{\psi_1}(\hat{\chi}) \\
          \widehat{\psi_2}(\hat{\chi}) \\
          \psi_1(\hat{\chi}) \\
          \psi_2(\hat{\chi}) \\
         \end{array}
       \right]=k\left[
         \begin{array}{c}
          \widehat{\psi_1}(\hat{\chi}) \\
          \widehat{\psi_2}(\hat{\chi}) \\
          \psi_1(\hat{\chi}) \\
          \psi_2(\hat{\chi}) \\
         \end{array}
       \right],
\end{equation}
where
$q_1=\begin{cases}{}
0 & {\rm for\ the\ DNLSE\ {\uppercase\expandafter{\romannumeral1}}}, \\
-\frac{1}{4}|q|^2&{\rm for\ the\ DNLSE\  {\uppercase\expandafter{\romannumeral2}}}, \\
\frac{1}{2}|q|^2& {\rm for\ the\ DNLSE\  {\uppercase\expandafter{\romannumeral3}}}.
\end{cases}$

The shape of the matrix eigenvalue problems (\ref{eigen_1}) is $4n_1\times4n_1$, and we can solve it by the `eig' function in MATLAB.

\begin{figure}[h]
  \centering
  \includegraphics[width=2.1in,height=1.9in]{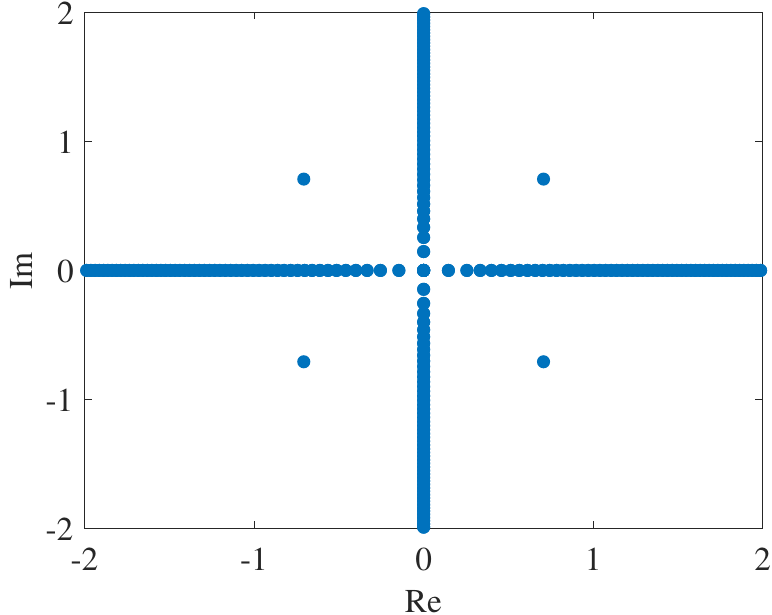}\\
  \caption{The calculated eigenvalues of the DNLSE $\rm \uppercase\expandafter{\romannumeral1}$  with (\ref{potential_0}) potential.}\label{direct_0}
\end{figure}

For the DNLSE $\uppercase\expandafter{\romannumeral1}$, Kaup and Newell got the soliton solutions by the inverse scattering transform\cite{KN}.
The function
\begin{equation}
q_{KN}(x,0)=4\frac{{\exp}(2x)\Big({\exp}(4x)+{\exp}(-\frac{\pi i}{2})\Big)}{\Big({\exp}(4x)+{\exp}(\frac{\pi i}{2})\Big)^2}
\label{potential_0}
\end{equation}
is the initial value for a pure soliton solution.
There are four discrete eigenvalues $k_1=\frac{\sqrt{2}}{2}+\frac{\sqrt{2}}{2}{i}$, $\bar{k}_1=\frac{\sqrt{2}}{2}-\frac{\sqrt{2}}{2}{i}$, $-\bar{k}_1=-\frac{\sqrt{2}}{2}+\frac{\sqrt{2}}{2}{i}$ and $-k_1=-\frac{\sqrt{2}}{2}-\frac{\sqrt{2}}{2}{i}$ for (\ref{potential_0}).
We calculate the eigenvalues of (\ref{potential_0}) for the DNLSE $\rm \uppercase\expandafter{\romannumeral1}$ under $\hat{a}=0.1$ and $n_1=600$, the calculated results are shown in Figure \ref{direct_0}.
The absolute error between the calculated result and the exact result is $4.80\times10^{-15}$ at $k=k_1$, the absolute error is $1.94\times10^{-15}$ at $k=\bar{k}_1$, the absolute error is $9.89\times10^{-15}$ at $k=-\bar{k}_1$, and the absolute error is $6.34\times10^{-15}$ at $k=-k_1$.

For the DNLSE $\rm \uppercase\expandafter{\romannumeral3}$, Fan got the soliton solutions by the Darboux transformation\cite{DB_fan}.
The DNLSE $\rm \uppercase\expandafter{\romannumeral3}$ exists the follow one-solution solution
\begin{equation}\label{dnls3_solution}
q_{GI}(x,t)=-8\frac{\xi_1\eta_1}{(\xi_1+i\eta_1){\rm e}^{(X_1+iY_1)}+(\xi_1-i\eta_1){\rm e}^{(-X_1+iY_1)}},
\end{equation}
where $\xi_1$ is the real part of the discrete eigenvalue $k_2$ and $\eta_1$ is the imaginary part of the discrete eigenvalue $k_2$,
$$X_1=-4\xi_1\eta_1x-16\xi_1\mu_1(\xi_1^2-\eta_1^2)t,$$
$$Y_1=2(\xi_1^2-\eta_1^2)x+4(\xi_1^4+\eta_1^4-6\xi_1^2\eta_1^2)t.$$
There are four discrete eigenvalues $k_2=1+0.5i$, $\bar{k}_2=1-0.5i$, $-\bar{k}_2=-1+0.5i$ and $-k_2=-1-0.5i$ for
\begin{equation}
q_{GI}(x,0)=\frac{-4}{(1+0.5i){\rm e}^{-2x+1.5ix}+(1-0.5i){\rm e}^{2x+1.5ix}}.
\label{potential_1}
\end{equation}

\begin{figure}[h]
  \centering
  \includegraphics[width=2.1in,height=1.9in]{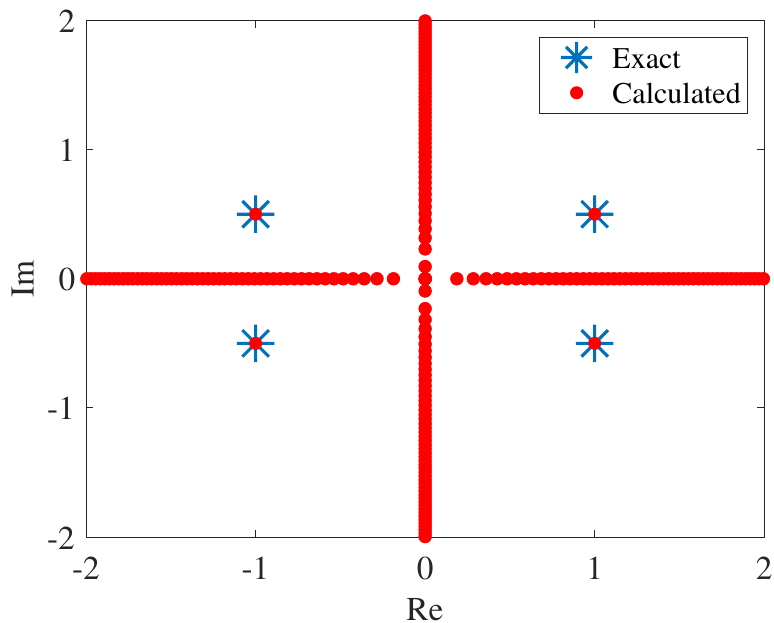}\\
  \caption{The eigenvalues for the DNLSE $\rm \uppercase\expandafter{\romannumeral3}$  with (\ref{potential_1}). Blue part: the exact result; red part: the calculated result.}\label{direct_1}
\end{figure}

We calculate the eigenvalues of $q_{GI}(x,0)$ potential with $\hat{a}=0.1$ and $n_1=400$, and the calculated results are shown in Figure \ref{direct_1}.
The absolute error between the calculated result and the exact result is $1.09\times10^{-15}$ at $k=k_2$, the absolute error is $4.66\times10^{-15}$ at $k=\bar{k}_2$, the absolute error is $1.01\times10^{-14}$ at $k=-\bar{k}_2$, and the absolute error is $2.50\times10^{-15}$ at $k=-k_2$.

\begin{figure}[h]
  \centering
  \includegraphics[width=2.2in,height=1.6in]{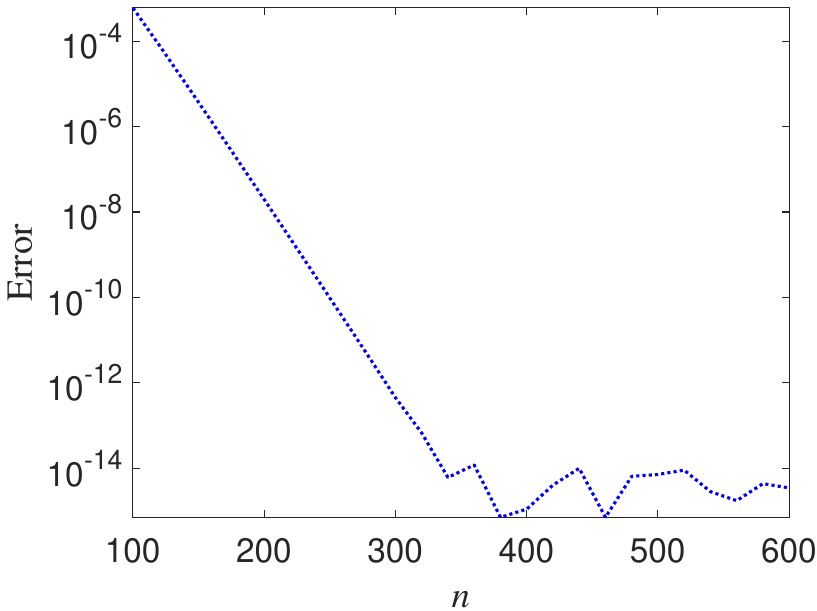}\\
  \caption{Demonstration of spectral convergence for the discrete eigenvalue $k_2$.}\label{direct_2}
\end{figure}

Although the numerical method performs well for the discrete eigenvalues and performs poorly for the continuous spectrum, it is still feasible because the continuous spectrum of the DNLSE is known.
For the DNLSE $\rm \uppercase\expandafter{\romannumeral3}$, we explore the spectral convergence for the discrete eigenvalue $k_2$, the calculated result is shown in Figure \ref{direct_2}. We learn that this method is effective and has machine accuracy.

\begin{remark}
In this numerical method, $\hat{a}$ should be a relatively small number and is related to the potential function $q_0(x)$.
We can initially set $\hat{a}$ to 0.1 and use the convergency to determine whether it is necessary to reduce $\hat{a}$.
\end{remark}
\subsection{Calculating the reflection coefficient}
$S(k)$ is determined by $\mu_1$ and $\mu_2$
 $$\mu_2(x,t,k)=\mu_{1}(x,t,k) {\rm e}^{-{i}\theta\hat{\sigma}_{3}} S(k).$$
Because $S(k)$ is independent of $x$ and $t$, we can set $x=0$ and $t=0$.
We mark $\mu_1(0,0,k)$ and $\mu_2(0,0,k)$ as $\mu_1(k)$ and $\mu_2(k)$ respectively.
The eigenfunctions $\mu_1(k)$ and $\mu_2(k)$ are calculated from the $x$-part of the Lax pair (\ref{lax_2}),
\begin{equation}\label{matrix_44}
(\ \mu_{j}(k)\ )_x+{i}k^{2}[\sigma_{3},\mu_j(k) ]=(kQ+Q_1)\mu_j(k), \ j= 1, 2,
\end{equation}
where $
Q_1=
\begin{cases}
	0\sigma_{3}  \quad &{\rm for\ the\ DNLS\  \uppercase\expandafter{\romannumeral1}}, \\
	-\frac{i}{4}|q(x, 0)|^{2}\sigma_{3}  \quad &{\rm for\ the\ DNLS\  \uppercase\expandafter{\romannumeral2}},\\
\frac{i}{2}|q(x, 0)|^{2}\sigma_{3}  \quad &{\rm for\ the\ DNLS\  \uppercase\expandafter{\romannumeral3}}.
\end{cases}
$
The eigenfunction $\mu_1(k)\rightarrow I$ when $x\rightarrow-\infty$, and the eigenfunction $\mu_2(k)\rightarrow I$ when $x\rightarrow+\infty$.

Set $\Phi_j(k)$=$\mu_j(k)-I$, ($j$=1, 2), then $\Phi_1(k)\rightarrow0$ when $x\rightarrow-\infty$, $\Phi_2(k)\rightarrow0$ when $x\rightarrow+\infty$.
$\Phi_j(k)$ satisfies
\begin{equation}\label{matrix_44}
(\ \Phi_{j}(k)\ )_x+{i} k^{2}[\sigma_{3},\Phi_j(k) ] - ( kQ+Q_1 )\Phi_j(k)=( kQ+Q_1), \ j=1,2.
\end{equation}

Set $$\Phi_j=\left(
             \begin{array}{cc}
               a & b \\
               c & d \\
             \end{array}
           \right),
$$ then $$ \left(\begin{array}{cc}
               a_x & b_x \\
               c_x & d_x \\
             \end{array}
             \right)
             +i k^2
             \left(
             \begin{array}{cc}
               0 & 2b \\
               -2c & 0 \\
             \end{array}
             \right)
             -
             \left(
             \begin{array}{cc}
               q_2 & kq \\
               -k\bar{q} & -q_2 \\
             \end{array}
             \right)\left(
             \begin{array}{cc}
               a & b \\
               c & d \\
             \end{array}
           \right)
           =
           \left(
             \begin{array}{cc}
               q_2 & kq \\
               -k\bar{q} & -q_2 \\
             \end{array}
             \right),
$$
where $q_2=\frac{i}{2}|q(x, 0)|^{2}$.

The following two equations can be obtained
\begin{equation}\label{cheb_1}
\left(\begin{array}{cc}
               a_x  \\
               c_x  \\
\end{array}\right)
             -
             \left(
             \begin{array}{cc}
               q_2 & kq \\
               -k\bar{q} & -q_2+2{i}k^2 \\
             \end{array}
             \right)\left(
             \begin{array}{cc}
               a \\
c \\
             \end{array}
           \right)
           =
           \left(
             \begin{array}{cc}
               q_2\\
               -k\bar{q}\\
             \end{array}
             \right),
\end{equation}

\begin{equation}\label{cheb_2}
\left(\begin{array}{cc}
               b_x \\
               d_x \\
             \end{array}
             \right)
             -
             \left(
             \begin{array}{cc}
               q_2-2{i}k^2& kq \\
               -k\bar{q} & -q_2 \\
             \end{array}
             \right)\left(
             \begin{array}{cc}
               b \\
               d \\
             \end{array}
           \right)
           =
           \left(
             \begin{array}{cc}
              kq \\
             -q_2 \\
             \end{array}
             \right).
\end{equation}
When $j$=1,
$\left(\begin{array}{cc}
               a & b\\
               c & d\\
             \end{array}
             \right)\rightarrow\left(\begin{array}{cc}
               0 & 0\\
               0 & 0\\
             \end{array}
             \right)$ as $x\rightarrow-\infty$,
when $j$=2, $\left(\begin{array}{cc}
               a & b\\
               c & d\\
             \end{array}
             \right)\rightarrow\left(\begin{array}{cc}
               0 & 0\\
               0 &0 \\
             \end{array}
             \right)$ as $x\rightarrow+\infty$.

The Chebyshev collocation method can be used to solve Eq. (\ref{cheb_1}) and Eq. (\ref{cheb_2}).
Using the knowledge in section \ref{sec_eigenvalue}, Eq. (\ref{cheb_1}) becomes
\begin{equation}\label{lax81}
\begin{split}
\left[
      \begin{array}{cc}
       {\rm D}[\frac{{\rm d} H(\vec{\chi})}{{\rm d} x}]\mathcal{F}^{-1}\mathcal{D}\mathcal{F}-{\rm D}[q_2(\hat{\chi})]&-k{\rm D}[q(\hat{\chi})] \\
       k{\rm D}[\bar{q}(\hat{\chi})]&{\rm D}[\frac{{\rm d} H(\vec{\chi})}{{\rm d} x}]\mathcal{F}^{-1}\mathcal{D}\mathcal{F}+{\rm D}[q_2(\hat{\chi})]-2ik^2I\\
      \end{array}
    \right]\left[
             \begin{array}{c}
               a(\hat{\chi}) \\
               c(\hat{\chi}) \\
             \end{array}
           \right]=\left[
             \begin{array}{c}
              q_2(\hat{\chi}) \\
             -k\bar{q}(\hat{\chi}) \\
             \end{array}
           \right].
\end{split}
\end{equation}
Eq. (\ref{cheb_2}) becomes \begin{equation}\label{lax82}
\begin{split}
\left[
      \begin{array}{cc}
       2ik^2I+{\rm D}[\frac{{\rm d} H(\vec{\chi})}{{\rm d} x}]\mathcal{F}^{-1}\mathcal{D}\mathcal{F}-{\rm D}[q_2(\hat{\chi})]&-k{\rm D}[q(\hat{\chi})] \\
       k{\rm D}[\bar{q}(\hat{\chi})]& {\rm D}[\frac{{\rm d} H(\vec{\chi})}{{\rm d} x}]\mathcal{F}^{-1}\mathcal{D}\mathcal{F}+{\rm D}[q_2(\hat{\chi})]\\
      \end{array}
    \right]\left[
             \begin{array}{c}
               b(\hat{\chi}) \\
               d(\hat{\chi}) \\
             \end{array}
           \right]=\left[
             \begin{array}{c}
               kq(\hat{\chi})  \\
              -q_2(\hat{\chi})\\
             \end{array}
           \right].
\end{split}
\end{equation}

\begin{remark}
The mapping $H(x)$ used in this subsection is different from subsection \ref{sec_eigenvalue}. When $j=1$, $H(x)=2{\rm tanh}(\hat{a} x)+1$, when $j=2$, $H(x)=2{\rm tanh}(\hat{a} x)-1$.
\end{remark}

\begin{figure}[h]
  \centering
  \subfigure[The density of $r(k)$ in $k\in\mathbb{R}$]{\includegraphics[width=1.5in,height=1.3in]{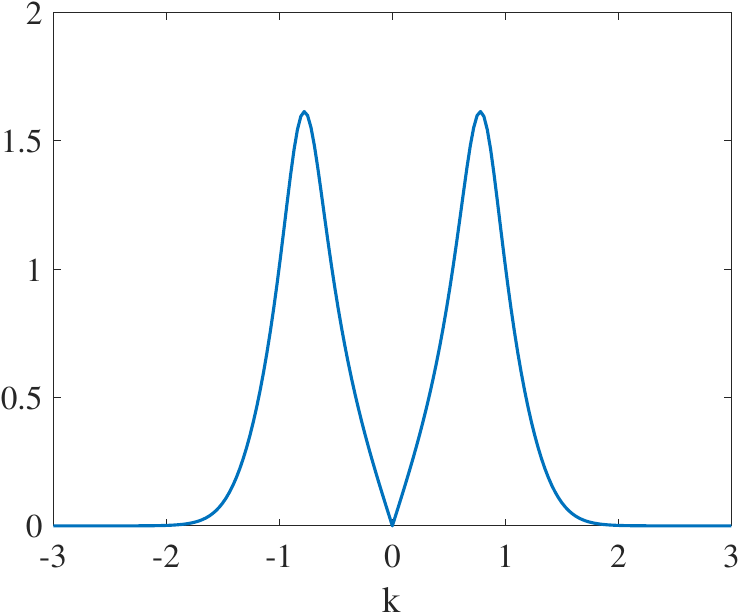}}
  \subfigure[]{\includegraphics[width=1.5in,height=1.3in]{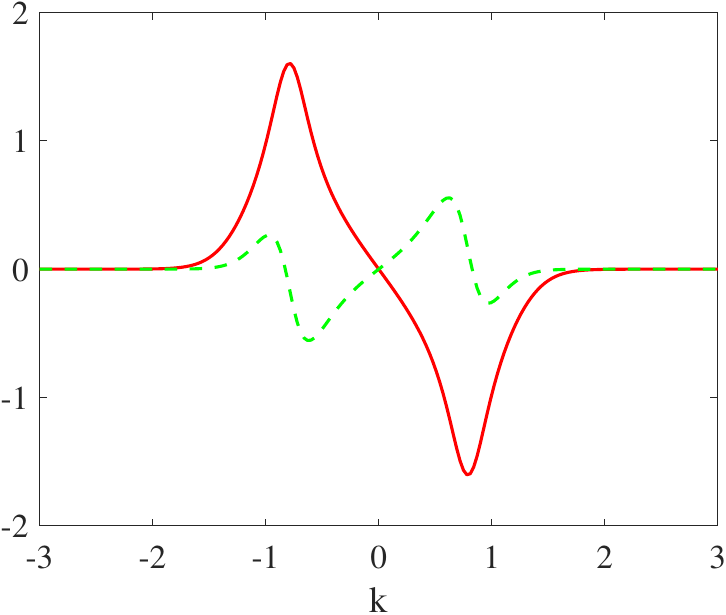}}
  \subfigure[The density of $r(k)$ in $k\in i\mathbb{R}$]{\includegraphics[width=1.5in,height=1.3in]{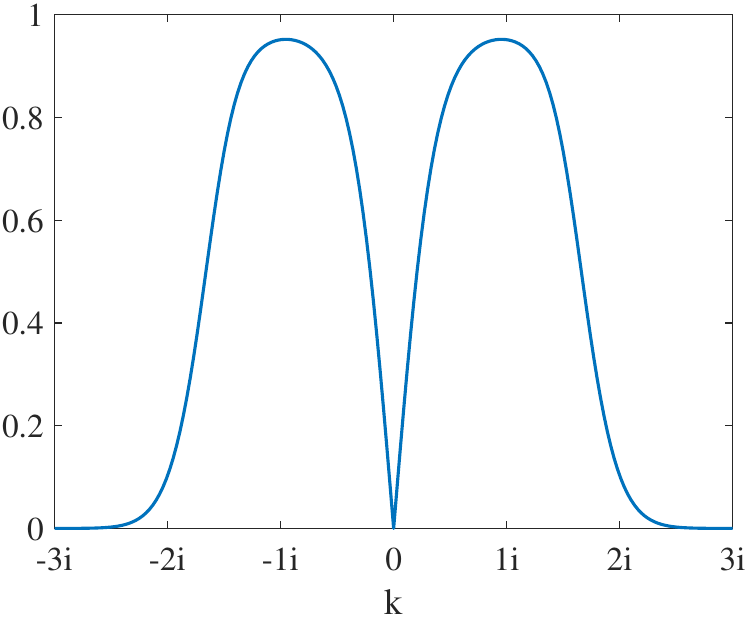}}
  \subfigure[]{\includegraphics[width=1.5in,height=1.3in]{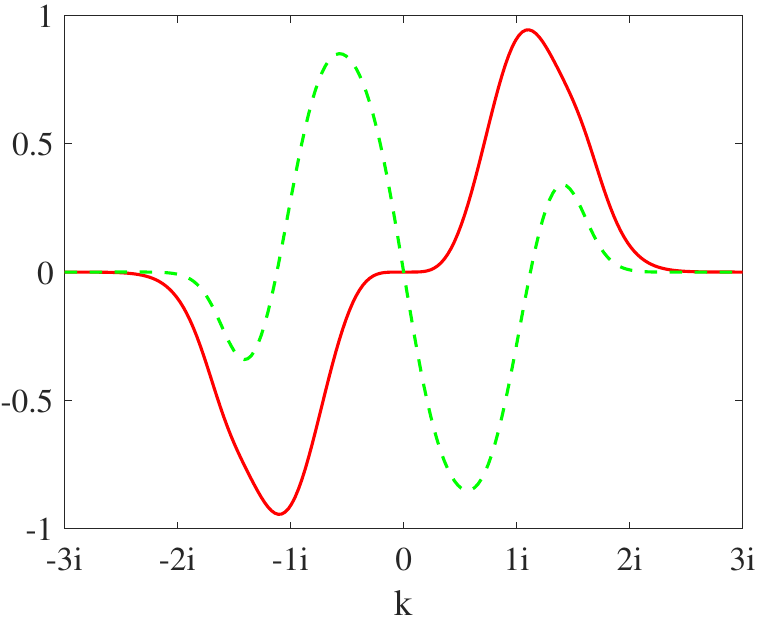}}\\
  \caption{The calculated reflection coefficient $r(k)$ for $q(x,0)=1.5{\rm exp}(-x^2)$. Blue line: the density of $r(k)$; red line: the real part of $r(k)$; green line: the imaginary part of $r(k)$}\label{fig_gauss}
\end{figure}

\begin{figure}[h]
  \centering
  \subfigure[$k=0.5$]{\includegraphics[width=1.5in,height=1.25in]{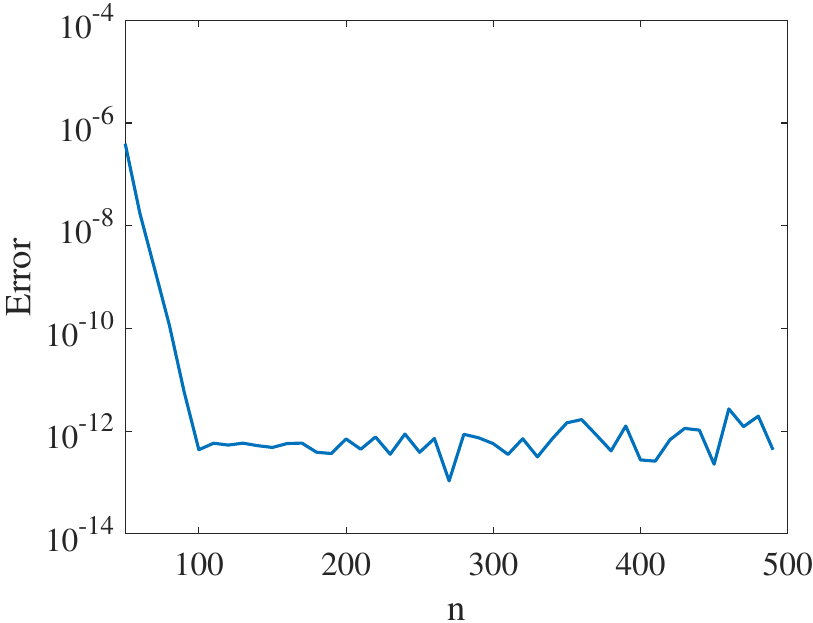}}
  \subfigure[$k=-1$]{\includegraphics[width=1.5in,height=1.25in]{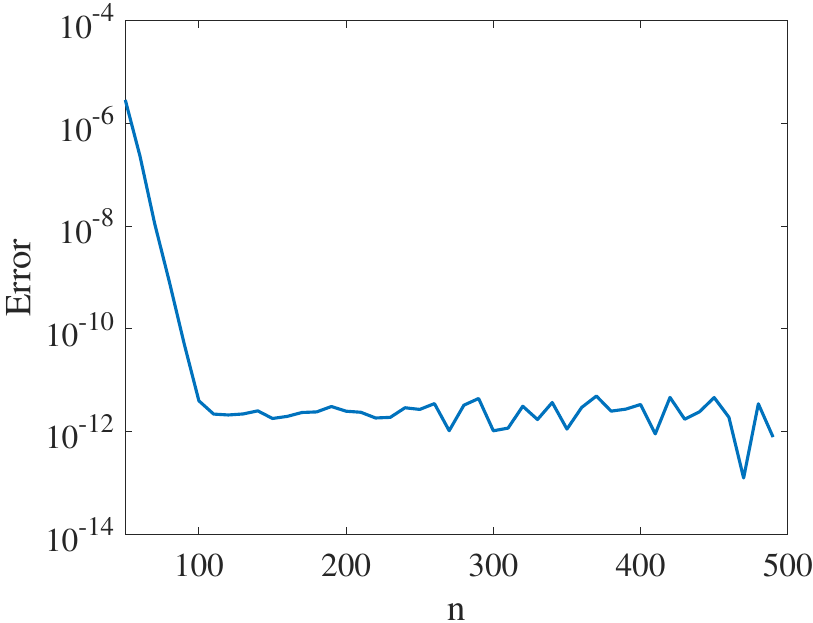}}
  \subfigure[$k=0.5i$]{\includegraphics[width=1.5in,height=1.25in]{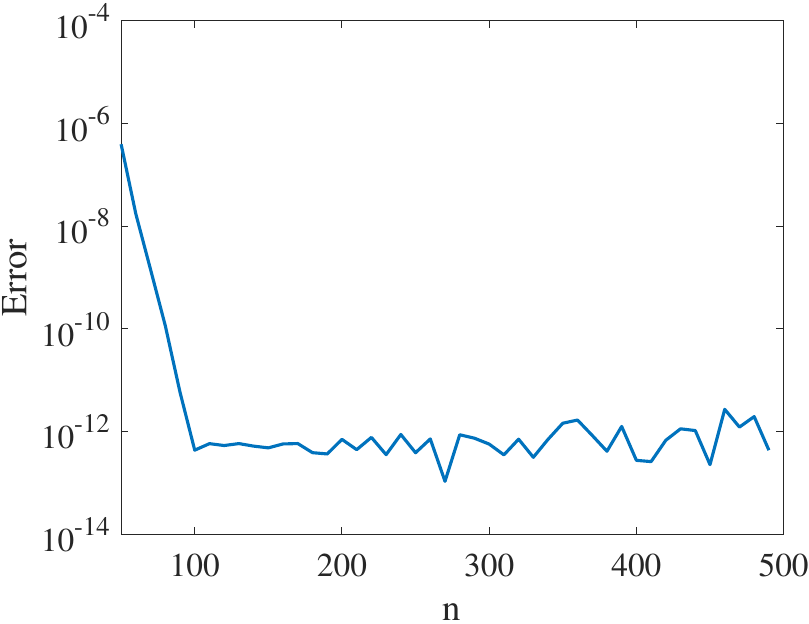}}
  \subfigure[$k=-i$]{\includegraphics[width=1.5in,height=1.25in]{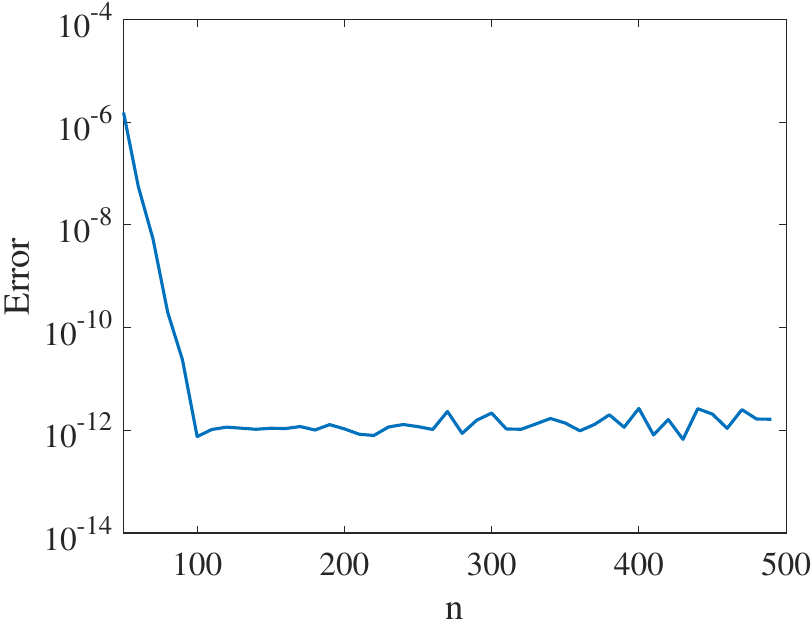}}\\
  \caption{Demonstration of spectral convergence for the reflection coefficient $r(k)$.}\label{fig_gauss_2}
\end{figure}

For the DNLSE $\rm \uppercase\expandafter{\romannumeral3}$, the numerical method ($n_1=201$, $\hat{a}=0.05$) is used to calculate the reflection coefficient $r(k)$ of $q(x,0)=1.5{\rm exp}(-x^2)$ potential. The calculated results are shown in Figure \ref{fig_gauss}.
It is worth noting that the continuous spectrum of the DNLSE is $\mathbb{R}\cup i\mathbb{R}$ rather than $\mathbb{R}$. In Figure \ref{fig_gauss}, the blue line is the density of $r(k)$, the red line is the real part of $r(k)$, the green line is the imaginary part of $r(k)$.
We also explore the convergency of $r(k)$ at $k=0.5$, $k=-1$, $k=0.5i$ and $k=-i$. The calculated results are shown in Figure \ref{fig_gauss_2}. The `Error' of Figure \ref{fig_gauss_2} is defined by $${\rm Error}=\Big| \ r(n_1,k)-r(500,k)\ \Big|,$$
where $r(n_1,k)$ is the reflection coefficient $r(k)$ calculated by $n_1$ Chebyshev nodes.
From Figure \ref{fig_gauss_2}, we learn that the calculated results converge quickly, so the numerical method is effective.

\begin{remark}
Considering that we calculate the reflection coefficient in Schwartz space $\mathcal{S}(\mathbb{R})$, it is also a feasible framework to truncate the calculation interval into a finite interval\cite{nist_Trogdon1,nist_cui}.
\end{remark}
\subsection{Calculating the normalized constant}
When $k=\kappa_j$, we need to calculate the normalized constant $C_j$. If $k=\kappa_j$ is the one-order pole of $a(k)$, the normalized constant is defined by
$$C_j=\frac{b_j}{a'(\kappa_j)},$$
where $\kappa_j$ is the discrete eigenvalue located in the first quadrant of the complex field $\mathbb{C}$.
Set $k=\kappa_j$, the $b_j$  can be obtained by Eq. (\ref{matrix_1}).
Considering that $a(k)$ is analytic in $D_2$, we give three frames for calculating $a'(\kappa_j)$:
\begin{itemize}
  \item (Central Difference Formula): $$a'(\kappa_j)=\lim_{|\Delta|\rightarrow0}\frac{a(\kappa_j+\vec{\Delta})-a(\kappa_j-\vec{\Delta})}{\vec{\Delta}}.$$
  \item (Fourier Transform):$$a'(z)=\mathbb{F}^{-1}[{i}k\mathbb{F}[a(z)]],$$
  where $\mathbb{F}$ represents the Fourier transform and $\mathbb{F}^{-1}$ represents the inverse Fourier transform.
  In the practical calculation, we use the fast Fourier transform (FFT) to calculate $a'(\kappa_j)$.
  \item (Chebyshev Derivative Matrix) \\ $${\rm If}\ a(z)=T^{\mathbb{R}}(z)\mathcal{F}a(\hat{\chi}),\ {\rm then}\ a'(z)=\frac{\partial H(z)}{\partial z}T^{\mathbb{R}}(z)\mathcal{D}\mathcal{F}a(\hat{\chi}).$$
\end{itemize}

In this paper, the third method is used to calculate $a'(\kappa_j)$.
\section{Numerical inverse scattering}\label{sec_inverse}
Numerical inverse scattering consists of two steps, one is deforming the RHP in the spirit of Deift-Zhou nonlinear steepest descent\cite{deift1}, another one is solving the deformed Riemann-Hilbert problem numerically by the Chebyshev collocation method.
Considering the original Riemann-Hilbert problem is oscillatory, the deformation is necessary when time is large.

This section is organized as follows.
In section \ref{deform_1_1}, we deform the RHP under not large time and the non-soliton case.
In section \ref{deform_1_2}, we deform the RHP under large time and the non-soliton case.
In section \ref{deform_2_1}, we deform the RHP under not large time and the soliton case.
In section \ref{deform_2_2}, some extra deformations are given.
In section \ref{sec_cheby}, the numerical method for solving the RHP is introduced.

\begin{remark}
The green line is used to represent the jump contours of the Riemann-Hilbert problem in the complex-$k$ plane, and the red line is used to represent $f(k)=0$ in the complex-$k$ plane, which applies to all figures presented in this section. The formula for $f(k)$ will be given in section \ref{deform_1_1}.
\end{remark}
\subsection{Deformations for the RHP (the non-soliton case)}\label{deform_1_1}
Assuming that the solution of the DNLS equation will not evolve into the soliton.
In this condition, the DNLS equation does not exist discrete eigenvalues, and we do not need to calculate the normalized constant.

In $\mathbb{R}\cup i\mathbb{R}$, ${\rm e}^{2i \theta(k)}$ and ${\rm e}^{2i \theta(k)}$ are oscillatory, especially when $t$ is large.
The jump matrix $T(k)$  has the oscillation term ${\rm e}^{2i\theta(k)}$ and ${\rm e}^{-2i\theta(k)}$ at the same time and the jump contours of \textbf{Riemann-Hilbert problem} is $\mathbb{R}\cup{i}\mathbb{R}$.
Due to the existence of oscillation terms, it is difficulty to numerically solve \textbf{Riemann-Hilbert problem}.
Therefore, we need to deform \textbf{Riemann-Hilbert problem} to keep the jump contours away from $\mathbb{R}\cup{ i}\mathbb{R}$.
The method for deforming the oscillatory Riemann-Hilbert problem is derived from Deift-Zhou nonlinear steepest descent method\cite{deift1}.
\begin{remark}
When time is small, we can use the NIST to calculate the solution of the DNLS equation without any deformation.
When time is big, the deformations of the Riemann-Hilbert problem are necessary.
\end{remark}

We need to find the saddle points of $\varphi(k)={i}\theta={i}(k^2x+2k^4t).$
Define $\lambda=k^2$, then $$\varphi(\lambda)={i}(\lambda x+2\lambda^2t).$$
Set $\varphi'(\lambda)=0$, and we get
\begin{equation}\label{rhp_2}
\lambda_0=-\frac{x}{4t}
\end{equation}
which satisfies $\varphi'(\lambda_0)=0$ and $\varphi''(\lambda_0)\neq0$ .
We get two saddle points $$k_1=\sqrt{ -\frac{x}{4t} }\ {\rm and}\ k_2=-\sqrt{ -\frac{x}{4t} }.$$
It is obvious that $$k_3=0$$ is also the saddle point of $\varphi(k)={i}\theta={i}(k^2x+2k^4t).$

Setting $\lambda={\rm Re}\lambda+i{\rm Im} \lambda $, we get
\begin{equation}\label{andian_1}
{\rm Re}({i}\theta)=-4{\rm Im}\lambda\cdot({\rm Re} \lambda-\lambda_0).
\end{equation}
Substitute $\lambda=k^2$ into Eq. (\ref{andian_1}), we get
\begin{equation}\label{andian_11}
{\rm Re}({i}\theta)=-4{\rm Im}\ k^2\cdot({\rm Re}{\ k^2} -\lambda_0).
\end{equation}
Eq. (\ref{andian_11}) shows that the sign of $\rm Re({i}\theta)$ is determined by ${\rm Im}\ {k^2}$ and $({\rm Re}\ {k^2}-\lambda_0)$.

\begin{figure}[h]
  \centering
  \includegraphics[width=2.2in,height=2.in]{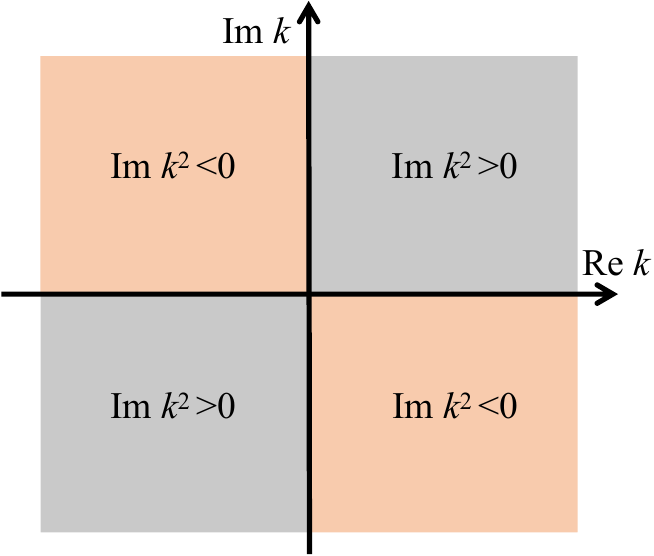}
  \caption{The sign of ${\rm Im}\ k^2$ on complex-$k$ plane. Gray area: $\rm Im \ {\it k^2}>0$. Orange area: $\rm Im \ {\it k^2}<0$.}\label{rekk}
\end{figure}

Given that $$ {\rm Im}\ (k^2)= {\rm Im}\big[ ({\rm Re}\ k+ i{\rm Im}\ k)^2\big]=2{\rm Re}\ k\cdot{\rm Im}\  k,$$
 the sign of ${\rm Im}\  k^2$ is shown in Figure \ref{rekk}.

After determining the sign of ${\rm Im}\ k^2$, we need to determine the sign of ${\rm Re}\ k^2 -\lambda_0$.
We define $f(k)$ by
\begin{equation}\label{fk}
f(k)={\rm Re}\ k^2- \lambda_0={\rm Re}\Big[ ({\rm Re}\ k+i{\rm Im}\ k)^2\Big]-\lambda_0=({\rm Re}\ k)^2-({\rm Im}\ k)^2-\lambda_0.
\end{equation}

\begin{figure}[h]
  \centering
  \subfigure[$x<0$]{\includegraphics[width=1.9in,height=1.8in]{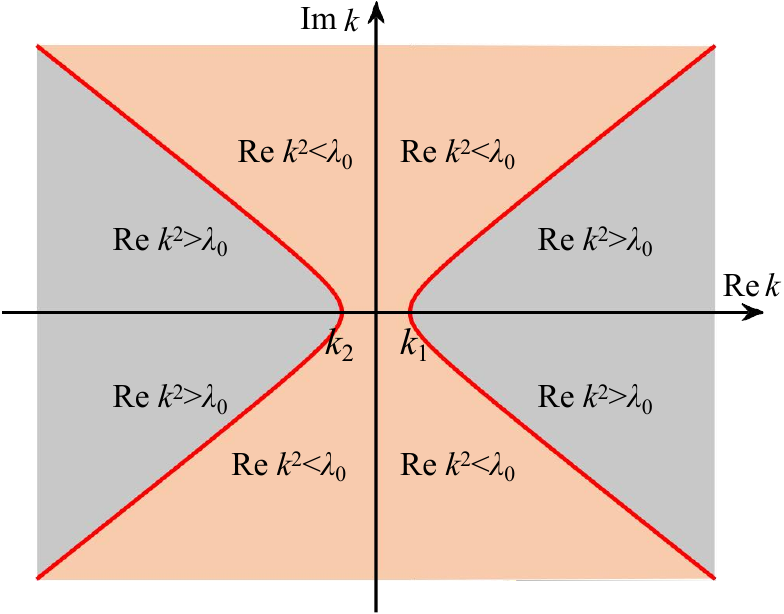}}
  \subfigure[$x=0$]{\includegraphics[width=1.9in,height=1.8in]{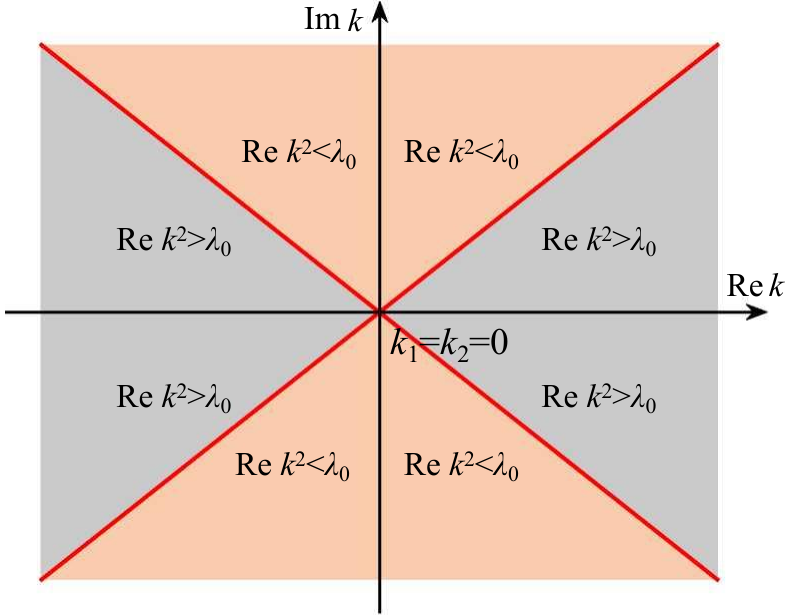}}
  \subfigure[$x>0$]{\includegraphics[width=1.9in,height=1.8in]{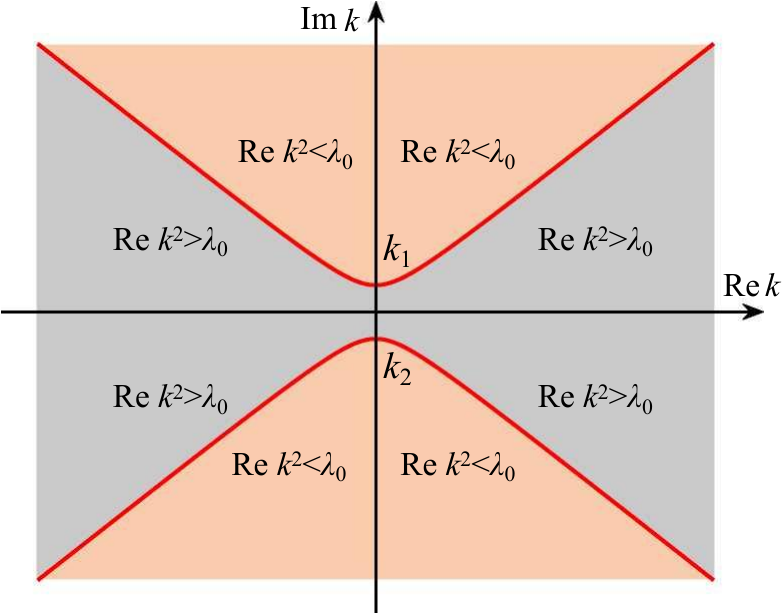}}
  \caption{The sign of ${\rm Re}\ {k^2} -\lambda_0$ in the complex-$k$ plane. Gray area: ${\rm Re}\ {k^2} >\lambda_0$. Orange area: ${\rm Re} {k^2}<\lambda_0$. }\label{rek}
\end{figure}

The shape of $f(k)=0$  and the sign of ${\rm Re}\ k^2 -\lambda_0$  are related to $x$.
 \begin{itemize}
   \item If $x<0$, $f(k)=0$ is a hyperbola whose focus is located in $\mathbb{R}$, whose asymptotes are Re $k$= $\pm$ Im $k$.
   \item If $x=0$, $f(k)=0$ are two straight lines Re $k$= $\pm$ Im $k$.

   \item If $x>0$, $f(k)=0$ is a hyperbola whose focus is located in $i\mathbb{R}$, whose asymptotes are Re $k$= $\pm$ Im $k$.
 \end{itemize}
The sign of ${\rm Re}\ {k^2}-\lambda_0$ is shown in Figure \ref{rek} where the red line is $f(k)=0$.

\begin{figure}[h]
  \centering
  \subfigure[$x<0$]{\includegraphics[width=1.9in,height=1.8in]{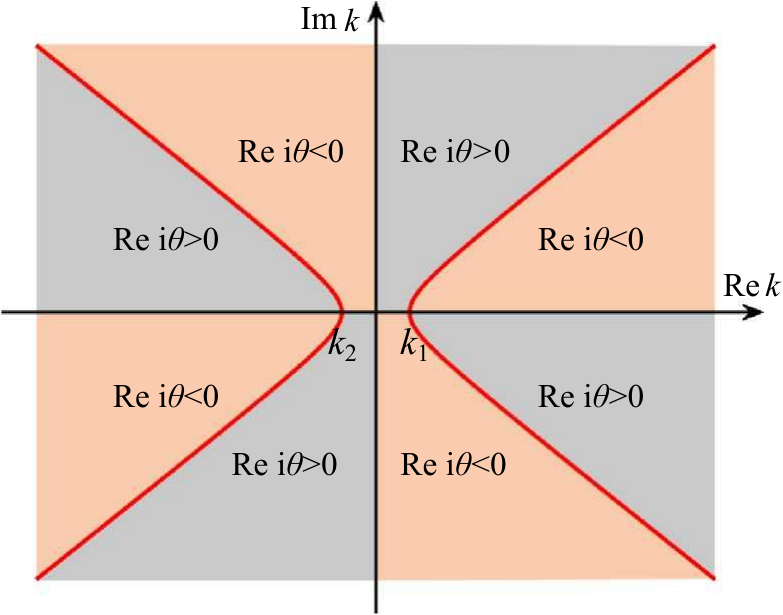}}
  \subfigure[$x=0$]{\includegraphics[width=1.9in,height=1.8in]{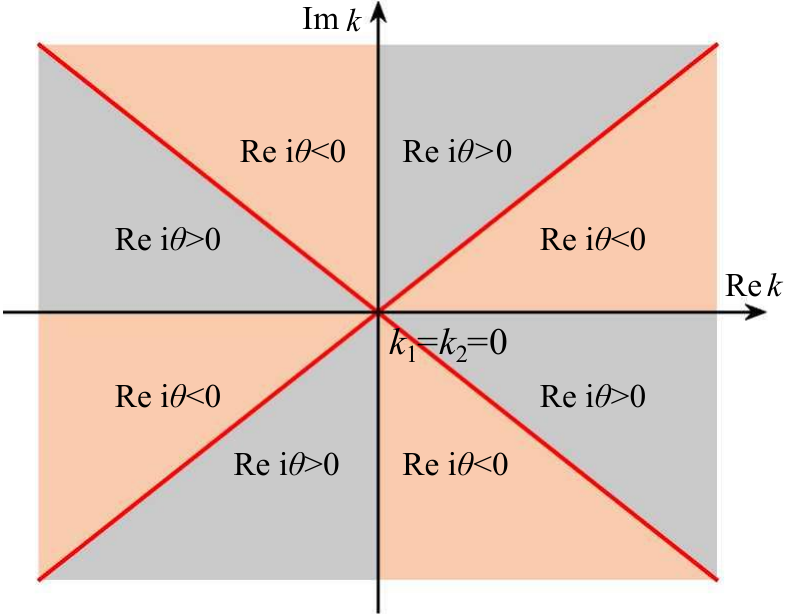}}
  \subfigure[$x>0$]{\includegraphics[width=1.9in,height=1.8in]{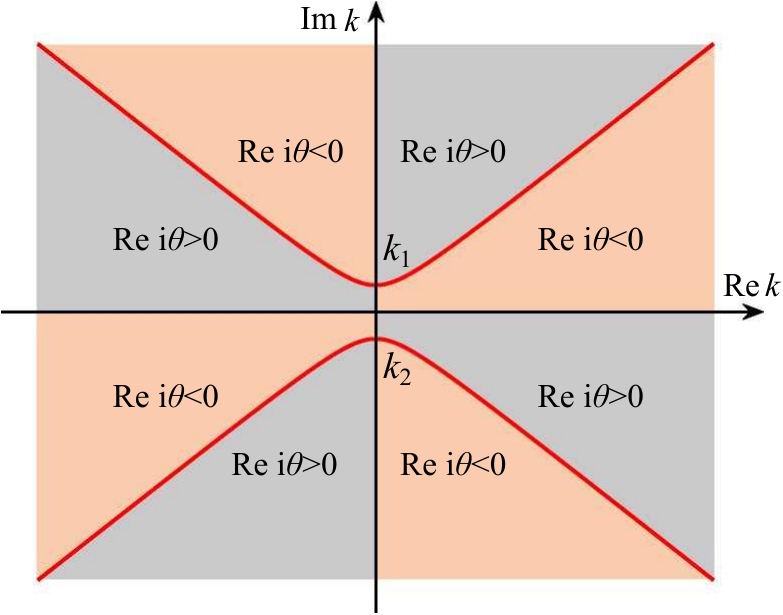}}
  \caption{The sign of ${\rm Re}({i}\theta)$ in the complex-$k$ plane. Gray area $\mathrm{F}_1$: ${\rm Re}(i\theta)>0$. Orange area $\mathrm{F}_2$: ${\rm Re}(i\theta)<0$}\label{rek2}
\end{figure}

The sign of ${\rm Re}(i\theta)$ in the complex-$k$ plane is shown in Figure \ref{rek2}. The gray area of the Figure \ref{rek2} is where ${\rm Re}(i\theta)>0$, and we record it as $\mathrm{F}_1$. The orange area of the Figure \ref{rek2} is where ${\rm Re}(i\theta)<0$, and we record it as $\mathrm{F}_2$.

After determining the sign of ${\rm Re}(i\theta)$, we start to deform \textbf{Riemann-Hilbert problem}.
Two decompositions of the jump matrix $T(x,t,k)$ are considered, where
$$T=\left(
\begin{array}{cc}
 1 & r(k){\rm e}^{-2i \theta(k)} \\
 -\tilde{r}(k){\rm e}^{2i \theta(k)} & 1-\tilde{r}(k)r(k) \\
  \end{array}
  \right).
$$
The first decomposition is $T=LU$, where
$$ L=\left(
\begin{array}{cc}
 1& 0 \\
  -\tilde{r}(k){\rm e}^{2i\theta(k)} & 1 \\
  \end{array}
  \right)
,\ U=  \left(
\begin{array}{cc}
 1 & r(k){\rm e}^{-2i\theta(k)}  \\
  0 & 1 \\
  \end{array}
  \right).
  $$
The second decomposition is $T=ABC$, where
$$\ A=  \left(
\begin{array}{cc}
 1 & \frac{r(k){\rm e}^{-2i \theta(k)}}{ 1-\tilde{r}(k)r(k)} \\
  0 & 1 \\
  \end{array}
  \right),
    \ B=  \left(
\begin{array}{cc}
 \frac{1}{1-\tilde{r}(k)r(k)}& 0 \\
  0 & 1-\tilde{r}(k)r(k) \\
  \end{array}
  \right),
  C=\left(
\begin{array}{cc}
 1& 0 \\
 \frac{-\tilde{r}(k){\rm e}^{2i\theta(k)}}{ 1-\tilde{r}(k)r(k)} & 1 \\
  \end{array}
  \right).$$
In $\mathrm{F}_2$ , ${\rm e}^{2i\theta(k)}$ decays exponentially; in $\mathrm{F}_1$, ${\rm e}^{-2i\theta(k)}$ is decays exponentially.
Matrices $A$ and $U$ can be extended to the orange region $\mathrm{F}_1$, matrices $C$ and $L$ can be extended to the grey region $\mathrm{F}_2$.

There are three saddle points, and the position of the saddle points are determined by $x$ and $t$.
It is very complex to deform the original Riemann-Hilbert problem, and we will give specific expressions in different situations.
For the DNLS equation, we divide the $(x,t)$-plane into three regions($c_1>0$):
\begin{itemize}
  \item Region 1: $-c_1\sqrt{t}\leq x\leq c_1\sqrt{t}$.
  \item Region 2: $x<-c_1\sqrt{t}$.
  \item Region 3: $x>c_1\sqrt{t}$.
\end{itemize}
In different regions for the $(x,t)$-plane, the deformations are different.
After the proper deformations of the jump contours, the oscillatory Riemann-Hilbert problem is transformed into a standard Riemann-Hilbert problem or a weak-oscillatory Riemann-Hilbert problem, the jump matrices are exponentially decaying away from the saddle points.
\subsubsection{Deformation in Region 1(-$c_1\sqrt{t}\leq x\leq c_1\sqrt{t}$)}
In this condition, we direct deform the jump contours at $k=0$, not at the saddle points $k_1$ and $k_2$.

\begin{figure}[h]
  \centering
  \includegraphics[width=2.6in,height=2.2in]{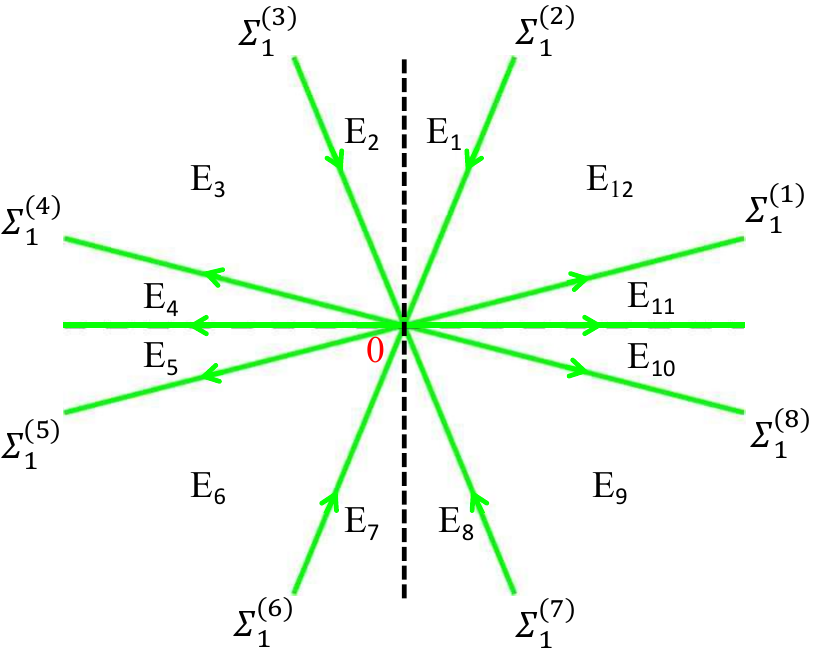}
  \caption{The jump contours for \textbf{Riemann-Hilbert problem \uppercase\expandafter{\romannumeral1}}.}\label{area_2}
\end{figure}

We define $M^{(1)}(x,t,k)$ by
\begin{equation}
M^{(1)}(x,t,k)=M(x,t,k)
\begin{cases}
I, & {\rm if}\ k\in \mathrm{E}_3\cup \mathrm{E}_6\cup  \mathrm{E}_9\cup \mathrm{E}_{12}, \\
C^{-1}, &  {\rm if}\ k\in \mathrm{E}_5\cup \mathrm{E}_{11},\\
A, & {\rm if}\ k\in \mathrm{E}_4\cup \mathrm{E}_{10},\\
L, & {\rm if}\ k\in \mathrm{E}_2\cup \mathrm{E}_{8},\\
U^{-1}, & {\rm if}\ k\in \mathrm{E}_1\cup \mathrm{E}_{7},
\end{cases}
\end{equation}
where $\mathrm{E}_{j}\ (j=1,2, ... , 12)$ is shown in Figure \ref{area_2}.

Then, we get \textbf{Riemann-Hilbert problem \uppercase\expandafter{\romannumeral1}}:
\begin{itemize}
  \vspace{-0.0cm}
  \item $M^{(1)}(x,t,k)$ is meromorphic in $\mathbb{C}\setminus\Sigma_1$,  $\Sigma_1$=$\Sigma^{(1)}_1\cup\cdots \cup\Sigma^{(8)}_1\cup\mathbb{R} $,
  \vspace{-0.0cm}
  \item $M^{(1)}_{+}(x,t,k)=M^{(1)}_{-}(x,t,k)T^{(1)}(x,t,k)$, if $k\in\Sigma_1$,
  \vspace{-0.0cm}
  \item $M^{(1)}(x,t,k)\rightarrow I$, $k\rightarrow\infty$,
\end{itemize}
where $\Sigma_1$ is shown in Figure \ref{area_2},
\begin{equation}
T^{(1)}(x,t,k)=
\begin{cases}
L, &  {\rm if}\ k\in\Sigma^{(3)}_1\cup\Sigma^{(7)}_1,\\
A,&{\rm if}\ k\in\Sigma^{(4)}_1\cup\Sigma^{(8)}_1, \\
C, & {\rm if}\ k\in\Sigma^{(1)}_1\cup\Sigma^{(5)}_1, \\
U, &{\rm if}\ k\in\Sigma^{(2)}_1\cup\Sigma^{(6)}_1, \\
B,&  {\rm if}\ k\in \mathbb{R}.
\end{cases}
\end{equation}

To ensure that this decomposition is reasonable, we give the following proposition:
\begin{pro}[Uniform boundedness]\label{pro1}
If $-c_1\sqrt{t}\leq x\leq c_1\sqrt{t}$, there is a constant $C_0$ s.t.
 \begin{itemize}
   \item $|{\rm e}^{2i \theta(k)}|\leq C_0$, {\rm for\ any}  $k\in \Sigma^{(1)}_1\cup\Sigma^{(3)}_1\cup\Sigma^{(5)}_1\cup\Sigma^{(7)}_1$,
   \item $|{\rm e}^{-2i\theta(k)}|\leq C_0$, {\rm for\ any}  $k\in \Sigma^{(2)}_1\cup\Sigma^{(4)}_1\cup\Sigma^{(6)}_1\cup\Sigma^{(8)}_1$.
 \end{itemize}
\end{pro}
\begin{proof}
Assume that the argument of $\Sigma^{(1)}_1$ is ${\rm arctan}(\frac{1}{10})$ and $c_1=4$.

First, we consider $x=0$.

In this condition, all of $\Sigma^{(j)}_1 (j=2, 4, 6, 8)$ are located in the orange area $\mathrm{F}_2$. In the orange area $\mathrm{F}_2$, ${{\rm Re}(i\theta)}\leq0$, so $|{\rm e}^{2i\theta(k)}|\leq 1$.
Similarly, all of $\Sigma^{(j)}_1 (j=1, 3, 5, 7)$ are located in the grey area $\mathrm{F}_1$. In the grey area $\mathrm{F}_1$, ${{\rm Re}(-i\theta)}\leq0$, then we get $|{\rm e}^{-2i\theta(k)}|\leq 1$.

Then, we consider $-4\sqrt{t}\leq x<0$.

In this condition, $|{\rm e}^{2i\theta(k)}|\leq 1$ for $\forall k\in\Sigma^{(j)}_1 (j=3, 7)$ and
$|{\rm e}^{-2i\theta(k)}|\leq 1$ for $\forall k\in\Sigma^{(j)}_1 (j=2, 6)$.
We need to prove that $|{\rm e}^{2i\theta(k)}|\leq C_0$ for $\forall k\in\Sigma^{(j)}_1 (j=1, 5)$ and $|{\rm e}^{-2i\theta(k)}|\leq C_0$ for $\forall k\in\Sigma^{(j)}_1 (j=4, 8)$, where $C_0$ is a constant.
Considering $\Sigma^{(1)}_1$ and $\Sigma^{(5)}_1$ are symmetric about $k=0$, $\Sigma^{(8)}_1$ and $\Sigma^{(4)}_1$ are symmetric about $k=0$, we only need to prove that $|{\rm e}^{2i\theta(k)}|\leq C_0$ for $\forall k\in\Sigma^{(1)}_1$ and $|{\rm e}^{-2i\theta(k)}|\leq C_0$ for $\forall k\in\Sigma^{(8)}_1$.

Both $\Sigma^{(1)}_1$ and  $\Sigma^{(7)}_1$ consist of two parts, one is located in the gray area $\mathrm{F}_1$ and another one is located in the orange area $\mathrm{F}_2$.
So $\Sigma^{(1)}_1$ and  $\Sigma^{(8)}_1$ can be written in the following form
$$ \Sigma^{(1)}_1=\Big(\Sigma^{(1)}_1\cap\mathrm{F}_2\Big)\cup\Big(\Sigma^{(1)}_1\cap\mathrm{F}_1\Big), \ \Sigma^{(8)}_1=\Big(\Sigma^{(8)}_1\cap\mathrm{F}_2\Big)\cup\Big(\Sigma^{(8)}_1\cap\mathrm{F}_1\Big).$$
Then
\begin{equation}
\begin{cases}
\Sigma^{(1)}_1=\{\ k\ |\ k=10b+bi\},& b\in(0,+\infty), \\
\Sigma^{(1)}_1\cap\mathrm{F}_1=\{\ k\ |\ k=10b+b{i}\},& b\in(0,b_0), \\
\Sigma^{(1)}_1\cap\mathrm{F}_2=\{\ k\ |\ k=10b+b{i}\},& b\in(b_0,+\infty),
\end{cases}
\end{equation}
where $b_0={\rm Im}(k_0)$ and $k_0=a_0+b_0{i}$ is the intersection of $\Sigma^{(1)}_1$ and $f(k)$.
By calculating
\begin{equation}
\begin{cases}
a_0=10b_0, \\
a_0^2-b_0^2=\lambda_0, \\
b_0>0,
\end{cases}
\end{equation}
we get $b_0=\frac{1}{\sqrt{99}}\sqrt{-\frac{x}{4t}}$
and $k_0=\frac{10}{\sqrt{99}}\sqrt{-\frac{x}{4t}}+\frac{1}{\sqrt{99}}\sqrt{-\frac{x}{4t}}{i}.$
So  $k\in \Sigma^{(1)}_1\cap\mathrm{F}_2$ is equivalent to  $b\in(0,\frac{1}{\sqrt{99}}\sqrt{-\frac{x}{4t}})$.

In the orange area $\mathrm{F}_2$, the density of ${\rm e}^{{i}\theta(k)}$ is smaller than one due to ${\rm Re}(\ i\theta(k)\ )<0$.
So there is $|{\rm e}^{2i\theta(k)}|\leq 1$ for $\forall k\in \mathrm{\Sigma^{(1)}_1}\cap\mathrm{F}_2$.
We only need to prove that $|{\rm e}^{2i\theta(k)}|\leq C$ for $\forall k\in\Sigma^{(1)}_1\cap\mathrm{F}_1$.
Substitute $k=10b+bi$ into $i\theta(k)$, we get
\begin{equation}\label{pro_c1}
{\rm Re}{(i\theta)}=-7920tb^4-20xb^2.
\end{equation}
To prove that ${\rm Re}{(i\theta)}$ is bounded in $\Sigma^{(1)}_1\cap\mathrm{F}_1$, we only need to prove that $f_1(b)$=${\rm Re}{(i\theta)}=-7920tb^4-20xb^2$ is bounded for $\forall b\in[0,b_0]$.
For $b\in[0,b_0]$, $max\{f_1(b)\}=f_1(\frac{1}{\sqrt{198}}\sqrt{-\frac{x}{4t}})=\frac{5x^2}{396t}.$
If $-4\sqrt{t}\leq x\leq4\sqrt{t}$, then $f_1(b)\leq\frac{20}{99}$, there is $|{\rm e}^{2i \theta(k)}|\leq{\rm e}^{\frac{40}{99}}$ for $\forall k\in \Sigma^{(1)}_1\cap\mathrm{F}_2$.

Similarly, we can prove that $|{\rm e}^{\rm -2i\it \theta(k)}|\leq1$ for $\forall k\in \Sigma^{(8)}_1\cap\mathrm{F}_1$ and $|{\rm e}^{-2i\theta(k)}|\leq{\rm e}^{\frac{40}{99}}$ for $\forall k\in \Sigma^{(8)}_1\cap\mathrm{F}_2$.

Last, we consider $0<x\leq4\sqrt{t}$.

In this condition, $|{\rm e}^{2i\theta(k)}|\leq 1$ for $\forall k\in \Sigma^{(j)}_1 (j=1, 5)$ and $|{\rm e}^{- 2i\theta(k)}|\leq 1$ for $\forall k\in \Sigma^{(j)}_1 (j=4, 8)$.
We need to prove that $|{\rm e}^{2i\theta(k)}|\leq C_0$ for $\forall k\in \Sigma^{(j)}_1 (j=3, 7)$ and
$|{\rm e}^{-2i\theta(k)}|\leq C_0$ for $\forall k\in \Sigma^{(j)}_1 (j=2, 6)$.
Considering $\Sigma^{(3)}_1$ and $\Sigma^{(7)}_1$ are symmetric about $k=0$, $\Sigma^{(2)}_1$ and $\Sigma^{(6)}_1$ are symmetric about $k=0$, we only need to prove that $|{\rm e}^{2i \theta(k)}|\leq C_0$ for $\forall k\in\Sigma^{(3)}_1$ and $|{\rm e}^{-2i\theta(k)}|\leq C_0$ for $\forall k\in\Sigma^{(2)}_1$.

Both $\Sigma^{(2)}_1$ and  $\Sigma^{(3)}_1$ consist of two parts, one in the gray area $\mathrm{F}_1$ and one in the orange area $\mathrm{F}_2$.
So $\Sigma^{(2)}_1$ and  $\Sigma^{(3)}_1$ can be written in the following form
$$ \Sigma^{(2)}_1=\Big(\Sigma^{(2)}_1\cap\mathrm{F}_1\Big)\cup\Big(\Sigma^{(2)}_1\cap\mathrm{F}_2\Big), \ \Sigma^{(3)}_1=\Big(\Sigma^{(3)}_1\cap\mathrm{F}_1\Big)\cup\Big(\Sigma^{(3)}_1\cap\mathrm{F}_2\Big).$$
Assume that the argument of $\Sigma^{(3)}_1$ is $\frac{\pi}{2}+{\rm arctan}(\frac{1}{10})$.
Then
\begin{equation}
\begin{cases}
\mathrm{\Sigma^{(3)}_1}=\{\ k\ |\ k=-0.1b+bi\},& b\in(0,+\infty), \\
\mathrm{\Sigma^{(3)}_1}\cap\mathrm{F}_1=\{\ k\ |\ k=-0.1b+bi\},& b\in(0,b_0), \\
\mathrm{\Sigma^{(3)}_1}\cap\mathrm{F}_2=\{\ k\ |\ k=-0.1b+bi\},& b\in(b_0,+\infty),
\end{cases}
\end{equation}
where $k_0=a_0+b_0i$ is the intersection of $\Sigma^{(3)}_1$ and $f(k)$, and $b_0={\rm Im}(k_0)$.
By calculating
\begin{equation}
\begin{cases}
a_0=-0.1b_0, \\
a_0^2-b_0^2=\lambda_0, \\
b_0>0,
\end{cases}
\end{equation}
we get $b_0=\frac{10}{\sqrt{99}}\sqrt{\frac{x}{4t}}$
and $k_0=-\frac{1}{\sqrt{99}}\sqrt{\frac{x}{4t}}+\frac{10}{\sqrt{99}}\sqrt{\frac{x}{4t}}i.$

In the orange area $\mathrm{F}_2$, the density of ${\rm exp}(\ i\theta(k)\ )$ is smaller than one due to ${\rm Re}(\  i\theta(k)\ )<0$.
So there is $|{\rm e}^{\rm 2i\it \theta(k)}|\leq 1$ for $\forall k\in\Sigma^{(3)}_1\cap\mathrm{F}_2$.
We only need to prove that $|{\rm e}^{2i \theta(k)}|\leq C_0$ for $\forall k\in\Sigma^{(3)}_1\cap\mathrm{F}_1$.
Substitute $k=-0.1b+bi$ into $i\theta(k)$, and we get
\begin{equation}\label{pro_c1}
{\rm Re}{(i\theta)}=-0.792tb^4+0.2xb^2.
\end{equation}
To prove that ${\rm Re}{(i\theta)}$ is bounded in $\Sigma^{(3)}_1\cap\mathrm{F}_1$, we only need to prove that $f_2(b)$=${\rm Re}{(i\theta)}=-0.792tb^4+0.2xb^2$ is bounded for $\forall b\in(0,b_0)$.
For $b\in(0,b_0)$, $max\{f_2(b)\}=f_2(\frac{10}{\sqrt{198}}\sqrt{-\frac{x}{4t}})=\frac{5x^2}{396t}.$
If $-4\sqrt{t}\leq x\leq4\sqrt{t}$, then $f_2(b)\leq\frac{20}{99}$ and $|{\rm e}^{2i\theta(k)}|\leq{\rm e}^{\frac{40}{99}}$ for $\forall k\in \Sigma^{(3)}_1\cap\mathrm{F}_2$.
Similarly, we can prove that $|{\rm e}^{-2i\theta(k)}|\leq1$ for $\forall k\in \Sigma^{(2)}_1\cap\mathrm{F}_1$ and $|{\rm e}^{-2i\theta(k)}|\leq{\rm e}^{\frac{40}{99}}$ for $\forall k\in \Sigma^{(2)}_1\cap\mathrm{F}_2$.

Our proof is complete, the constant $C_0$ is ${\rm e}^{\frac{40}{99}}$.

\end{proof}

Our decomposition transforms the original Riemann-Hilbert problem into a week-oscillatory Riemann-Hilbert problem.
\textbf{Proposition \ref{pro1}} gives the uniform boundedness of the jump matrices.
This decomposition can both ensure numerical accuracy and minimize computational costs.

\subsubsection{Deformation in Region 2 ($x<-c_1\sqrt{t}$)}
If we directly deform the jump contours at $k=0$, the jump matrices are not uniformly bounded, which will lead to numerical accuracy issues.
When $x<-c_1\sqrt{t}$, the saddle points $k=k_1$ and $k=k_2$ are located in $\mathbb{R}$, we need to decompose the jump contours at the saddle points $k=0$, $k=k_1$ and $k=k_2$.

\begin{figure}[h]
  \centering
  \includegraphics[width=3.3in,height=2.in]{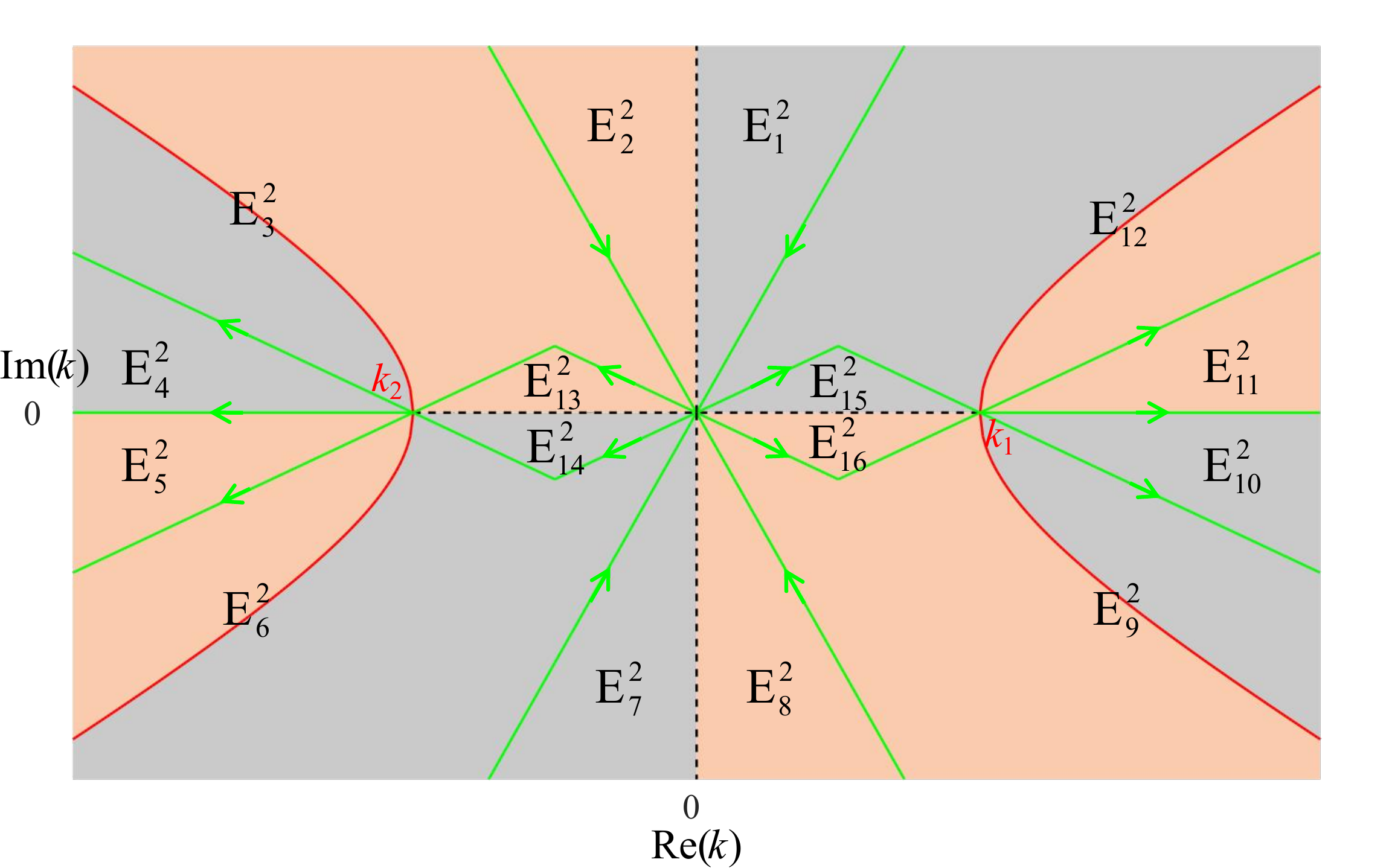}
  \caption{The area $\mathrm{E}^{(2)}_{j}\ (j=1,2, ... , 16)$ in the complex-$k$ plane.}\label{area_t1}
\end{figure}

We define $M^{(2)}(x,t,k)$ by
\begin{equation}
M^{(2)}(x,t,k)=M(x,t,k)
\begin{cases}
I, & {\rm if}\ k\in \mathrm{E}^2_3\cup \mathrm{E}^2_6\cup  \mathrm{E}^2_9\cup \mathrm{E}^2_{12}, \\
U^{-1}, &  {\rm if}\ k\in \mathrm{E}^2_1\cup\mathrm{E}^2_7\cup\mathrm{E}^2_{14}\cup\mathrm{E}^2_{16},\\
L, & {\rm if}\ k\in \mathrm{E}^2_2\cup \mathrm{E}^2_8\cup\mathrm{E}^2_{13}\cup\mathrm{E}^2_{15},\\
A, & {\rm if}\ k\in \mathrm{E}^2_4\cup \mathrm{E}^2_{10},\\
C^{-1}, & {\rm if}\ k\in \mathrm{E}^2_5\cup \mathrm{E}^2_{11},
\end{cases}
\end{equation}
where the area $\mathrm{E}^{2}_{j}\ (j=1,2, ... , 16)$ is shown in Figure \ref{area_t1}.

\begin{figure}[h]
  \centering
  \includegraphics[width=3.3in,height=2.in]{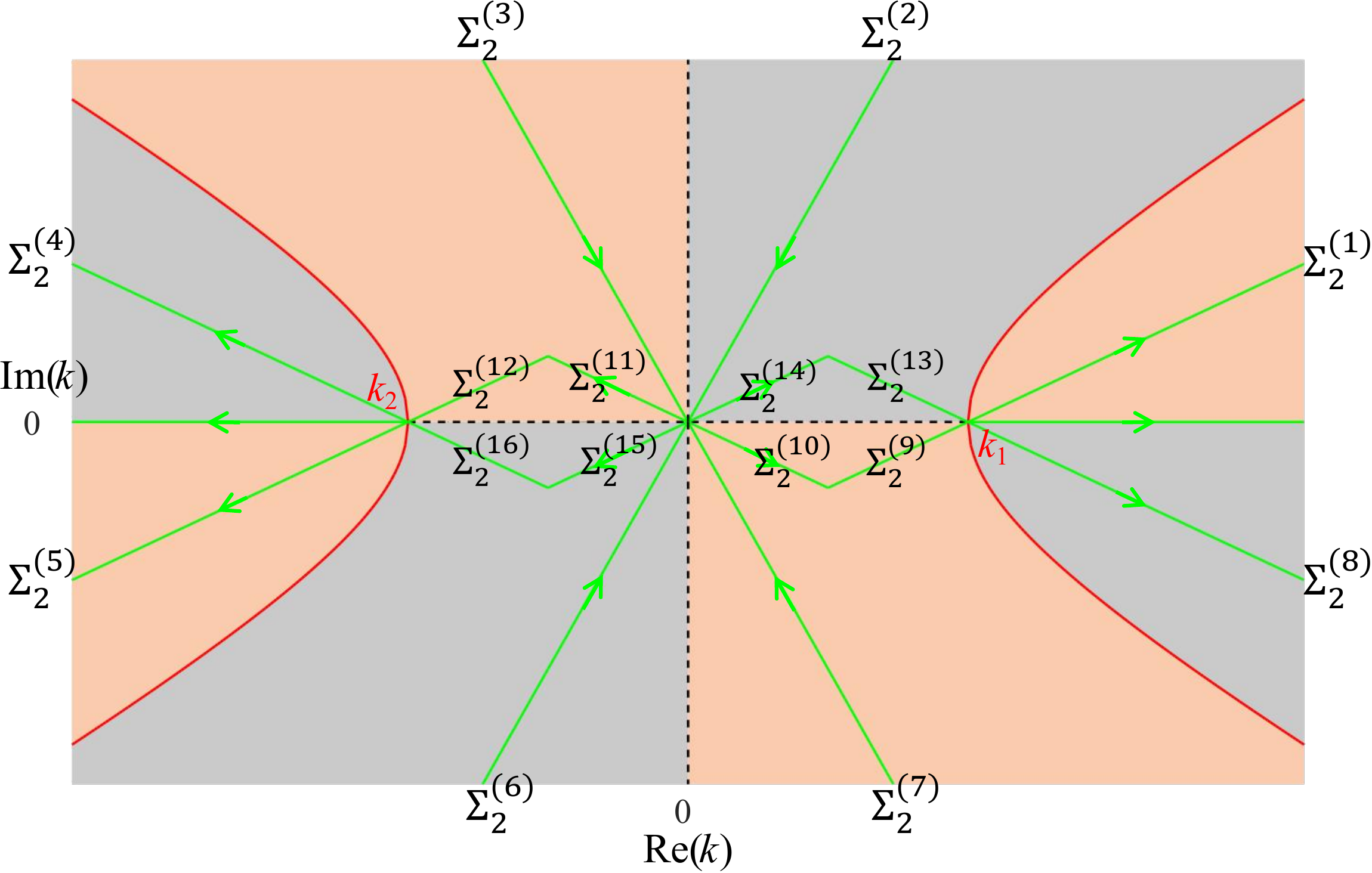}
  \caption{The jump contours for \textbf{Riemann-Hilbert problem \uppercase\expandafter{\romannumeral2}} in the complex-$k$ plane.}\label{jump_t1}
\end{figure}

We get \textbf{Riemann-Hilbert problem \uppercase\expandafter{\romannumeral2}}:\begin{itemize}
  \vspace{-0.0cm}
  \item $M^{(2)}(x,t,k)$ is meromorphic in $\mathbb{C}\setminus\Sigma_2 $,
  $\mathrm{\Sigma_2}$=$\Sigma^{(1)}_2\cup \cdots \cup\Sigma^{(16)}_2\cup(-\infty,k_2)\cup(k_1,+\infty)$,
  \vspace{-0.0cm}
  \item $M^{(2)}_{+}(x,t,k)=M^{(2)}_{-}(x,t,k)T^{(2)}(x,t,k)$, if $k\in\mathrm{\Sigma_2}$,
  \vspace{-0.0cm}
  \item $M^{(2)}(x,t,k)\rightarrow I$, $k\rightarrow\infty$,
\end{itemize}
where the jump contour $\Sigma^{(j)}_2 (j=1,2, ... ,16)$ is shown in Figure \ref{jump_t1},
\begin{equation}
T^{(2)}(x,t,k)=
\begin{cases}
A,&{\rm if}\ k\in\Sigma^{(4)}_2\cup\Sigma^{(8)}_2,\\
C,&{\rm if}\ k\in\Sigma^{(1)}_2\cup\Sigma^{(5)}_2,\\
U,&{\rm if}\ k\in\Sigma^{(2)}_2\cup\Sigma^{(6)}_2\cup\Sigma^{(13)}_2
\cup\Sigma^{(14)}_2\cup\Sigma^{(15)}_2\cup\Sigma^{(16)}_2,\\
L,&{\rm if}\ k\in\Sigma^{(3)}_2\cup\Sigma^{(7)}_2\cup\Sigma^{(9)}_2
\cup\Sigma^{(10)}_2\cup\Sigma^{(11)}_2\cup\Sigma^{(12)}_2,\\
B,&{\rm if}\ k\in(-\infty,k_2)\cup(k_1,+\infty).
\end{cases}
\end{equation}

There are 18 jump contours for \textbf{Riemann-Hilbert problem \uppercase\expandafter{\romannumeral2}}.
Compare to \textbf{Riemann-Hilbert problem \uppercase\expandafter{\romannumeral1}}, \textbf{Riemann-Hilbert problem  \uppercase\expandafter{\romannumeral2}} has 8 extra jump contours. The extra jump contours are $\Sigma^{(j)}_2 (j=9,10, ... ,16)$.

\subsubsection{Deformation in Region 3($x>c_1\sqrt{t}$)}
When $x>c_1\sqrt{t}$, the saddle points are located in ${i}\mathbb{R}$, we need to deform the jump contours at the saddle points $k=0$, $k=k_1$ and $k=k_2$.

\begin{figure}[h]
  \centering
  \includegraphics[width=3.3in,height=2.in]{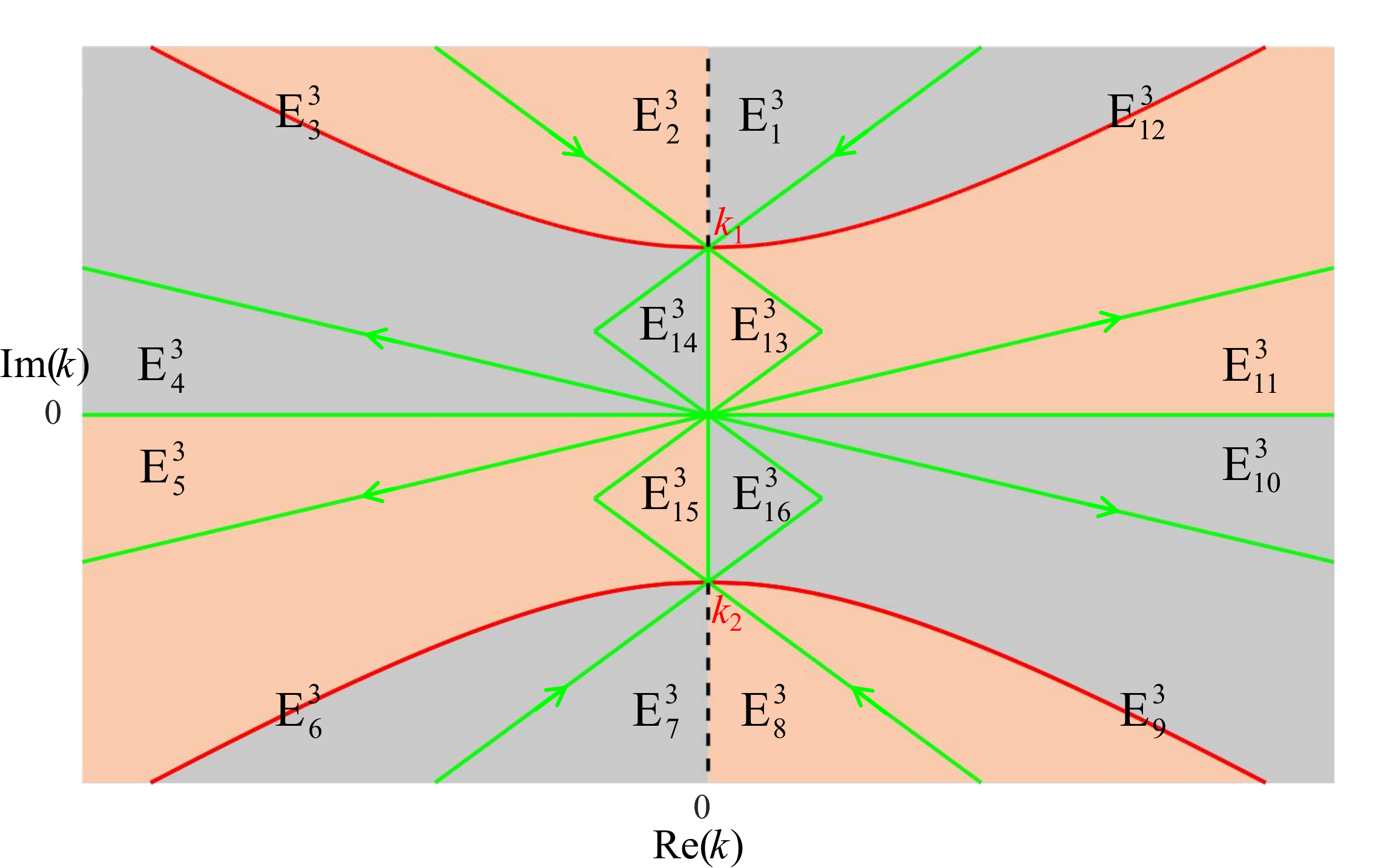}
  \caption{The area $\mathrm{E}^{3}_{j}\ (j=1,2, ... , 16)$ in the complex-$k$ plane.}\label{area_t2}
\end{figure}

We define $M^{(3)}(x,t,k)$ by
\begin{equation}
M^{(3)}(x,t,k)=M(x,t,k)
\begin{cases}
I,      &  {\rm if}\ k\in \mathrm{E}^3_3\cup \mathrm{E}^3_6\cup  \mathrm{E}^3_9\cup \mathrm{E}^3_{12}, \\
U^{-1}, &  {\rm if}\ k\in \mathrm{E}^3_1\cup\mathrm{E}^3_7\cup\mathrm{E}^3_{14}\cup\mathrm{E}^3_{16},\\
L,      &  {\rm if}\ k\in \mathrm{E}^3_2\cup \mathrm{E}^3_8\cup\mathrm{E}^3_{13}\cup\mathrm{E}^3_{15},\\
A,      &  {\rm if}\ k\in \mathrm{E}^3_4\cup \mathrm{E}^3_{10},\\
C^{-1}, &  {\rm if}\ k\in \mathrm{E}^3_5\cup \mathrm{E}^3_{11},
\end{cases}
\end{equation}
where the area $\mathrm{E}^{3}_{j}\ (j=1,2, ... , 16)$ is shown in Figure \ref{area_t2}.

\begin{figure}[h]
  \centering
  \includegraphics[width=3.3in,height=2.in]{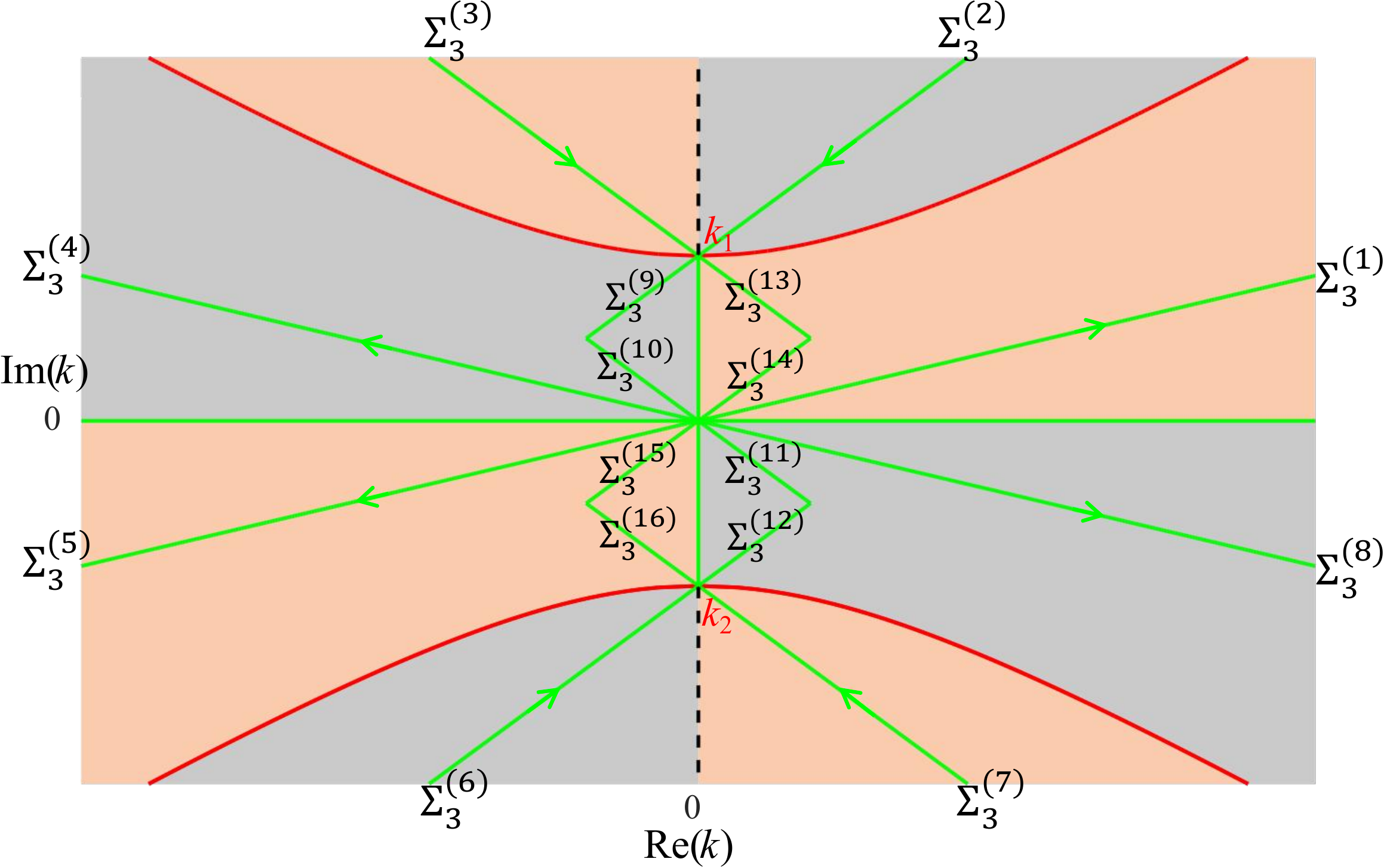}
  \caption{The jump contours for \textbf{Riemann-Hilbert problem \uppercase\expandafter{\romannumeral3}} in the complex-$k$ plane.}\label{jump_t2}
\end{figure}

We get \textbf{Riemann-Hilbert problem \uppercase\expandafter{\romannumeral3}}:
\begin{itemize}
  \vspace{-0.0cm}
  \item $M^{(3)}(x,t,k)$ is meromorphic in $\mathbb{C}\setminus\Sigma_3$, $\Sigma_3=\Sigma^{(1)}_3\cup \cdots \cup\Sigma^{(16)}_3\cup(0, k_2)\cup(0, k_1)\cup\mathbb{R}$,
  \vspace{-0.0cm}
  \item $M^{(3)}_{+}(x,t,k)=M^{(3)}_{-}(x,t,k)T^{(3)}(x,t,k)$, if $k\in\Sigma_3$,
  \vspace{-0.0cm}
  \item $M^{(3)}(x,t,k)\rightarrow I$, $k\rightarrow\infty$,
\end{itemize}
where the jump contour $\Sigma^{(j)}_3 (j=1,2, ... ,16)$ is shown in Figure \ref{jump_t2},
\begin{equation}
T^{(3)}(x,t,k)=
\begin{cases}
A, &  {\rm if}\ k\in\Sigma^{(4)}_3\cup\Sigma^{(8)}_3\cup\Sigma^{(9)}_3
\cup\Sigma^{(10)}_3\cup\Sigma^{(11)}_3\cup\Sigma^{(12)}_3,\\
C,&{\rm if}\ k\in\Sigma^{(1)}_3\cup\Sigma^{(5)}_3\cup\Sigma^{(13)}_3
\cup\Sigma^{(14)}_3\cup\Sigma^{(15)}_3\cup\Sigma^{(16)}_3,\\
L, & {\rm if}\ k\in\Sigma^{(3)}_3\cup\Sigma^{(7)}_3,\\
U, &{\rm if}\ k\in\Sigma^{(2)}_3\cup\Sigma^{(6)}_3,\\
B,&  {\rm if}\ k\in(0, k_2)\cup(0, k_1)\cup\mathbb{R}.
\end{cases}
\end{equation}

\subsection{Removing the jump matrix $B$ (the non-soliton case)}\label{deform_1_2}
Both \textbf{Riemann-Hilbert problem \uppercase\expandafter{\romannumeral1}}, \textbf{Riemann-Hilbert problem \uppercase\expandafter{\romannumeral2}} and \textbf{Riemann-Hilbert problem \uppercase\expandafter{\romannumeral3}} are not oscillatory. When $t$ is not too large, we can solve them by the Chebyshev collocation method.
When $t$ is large, we need to remove the jump matrix $B$.

\subsubsection{Removing the jump matrix $B$ in Region 1}\label{remove_1}
For \textbf{Riemann-Hilbert problem \uppercase\expandafter{\romannumeral1}}, we need to remove the jump on $(0, -\infty)$ and $(0, +\infty)$.
We consider a scalar Riemann-Hilbert problem:
\begin{itemize}
  \item $\delta_{j}(k) (j=1,2)$ is meromorphic in $\mathbb{C}\setminus{i}\mathbb{R}$,
  \item $\begin{cases}
  \delta_{j,+}(k)=\delta_{j,-}(k)(1+r(k)\overline{r(\bar{k})}),
  \begin{cases}
  k\in(0, -\infty)\ {\rm for}\ j=1, \\
  k\in(0, +\infty)\ {\rm for}\ j=2,
  \end{cases} \\
  \delta_{j,+}(k)=\delta_{j,-}(k),   \begin{cases}
  k\in(0, +\infty)\ {\rm for}\ j=1, \\
  k\in(0, -\infty)\ {\rm for}\ j=2,
  \end{cases}
  \end{cases}$
  \item $\delta_j(k)\rightarrow1$, $k\rightarrow\infty$.
\end{itemize}
When $j=1$, the Riemann-Hilbert problem exists the solution
\begin{equation}\label{delta_1}
\delta_1(k)={\exp}\Big[{\frac{1}{2\pi i}}\int_{0}^{-\infty}\frac{{\rm ln}(1+r(z)\overline{r(\bar{z})})}{z-k}{\rm d}z\Big],
\end{equation}
when $j=2$, the Riemann-Hilbert problem exists the solution
\begin{equation}\label{delta_1}
\delta_2(k)={\exp}\Big[{\frac{1}{2\pi i}}\int_{0}^{+\infty}\frac{{\rm ln}(1+r(z)\overline{r(\bar{z})})}{z-k}{\rm d}z\Big].
\end{equation}

Set $\Delta_{(j)}=\left(
              \begin{array}{cc}
               \frac{1}{\delta_j(k)} & 0 \\
                0 & \delta_j(k) \\
              \end{array}
            \right) (j=1,2)$, $\Delta_{(1)}$ can be used to remove the jump on $(0, -\infty)$,
             $\Delta_{(2)}$ can be used to remove the jump on $(0, +\infty)$.

For the RHP
\begin{equation}\label{scalar_rhp1}
M^{(1)}_{+}(x,t,k)=M^{(1)}_{-}(x,t,k)B, k\in(-\infty,0),
\end{equation}
we multiply the solution of the RHP (\ref{scalar_rhp1}) by $\Delta_{(1)}^{-1}$ to remove the jump of $(-\infty, 0)$.
In $(-\infty,0)$, $$M^{(1)}_{+}\Delta_{(1),+}^{-1}=M^{(1)}_{-}B\Delta_{(1),+}^{-1}
=M^{(1)}_{-}\Delta_{(1),-}^{-1}\Delta_{(1),-}B\Delta_{(1),+}^{-1}=M^{(1)}_{-}\Delta_{(1),-}^{-1}.$$
Indeed, we see that there is no jump for $M^{(1)}\Delta_{(1)}^{-1}$.
We can prove that $\delta_1(k)=\delta_2(-k)$ by the symmetries for $r(k)$.

Doing the transformation
\begin{equation}\label{deform_1}
M^{(4)}(x,t,k)=M^{(1)}(x,t,k)\Delta_1^{-1},
\end{equation}
where $\Delta_1^{-1}=\Delta_{(1)}^{-1}\Delta_{(2)}^{-1}=\left(
                                                     \begin{array}{cc}
                                                      \frac{1}{\delta_1(k)\delta_1(-k)}  &0 \\
                                                       0 & \delta_1(k)\delta_1(-k) \\
                                                     \end{array}
                                                   \right)$,
we get \textbf{Riemann-Hilbert\ problem \ \uppercase\expandafter{\romannumeral4}}:
\begin{itemize}
  \item $M^{(4)}(x,t,k)$ is meromorphic in $\mathbb{C}\setminus\Sigma_4$,
  \qquad $\Sigma_4=\Sigma^{(1)}_1\cup \cdots \cup\Sigma^{(8)}_1$,
  \vspace{-0.0cm}
  \item $M^{(4)}_{+}(x,t,k)=M^{(4)}_{-}(x,t,k)\begin{cases}
\Delta_1 L\Delta_1^{-1}, &{\rm if}\ k\in\Sigma^{(3)}_1\cup\Sigma^{(7)}_1,\\
\Delta_1A\Delta_1^{-1}, &{\rm if}\ k\in\Sigma^{(4)}_1\cup\Sigma^{(8)}_1,\\
\Delta_1C\Delta_1^{-1}, & {\rm if}\ k\in\Sigma^{(1)}_1\cup\Sigma^{(5)}_1,\\
\Delta_1U\Delta_1^{-1}, &{\rm if}\ k\in\Sigma^{(2)}_1\cup\Sigma^{(6)}_1.
\end{cases}$
  \vspace{-0.0cm}
  \item $M^{(4)}(x,t,k)\rightarrow I$, $k\rightarrow\infty$,
\end{itemize}
where $\Sigma_1^{(j)}(j=1,2,...,8)$ is shown in Figure \ref{area_2}.

\begin{remark}
For the DNLSE \uppercase\expandafter{\romannumeral1}, $k=0$ is not the singularity of $\delta_1(k)$, so the deformation (\ref{deform_1}) is sufficient.
\end{remark}
\begin{proof}
We divide $\delta_1(0)$ into two parts,
$$\delta_1(0)={\rm exp}\Big[{\frac{1}{2\pi{i}}}\int^{\varepsilon}_{-\infty}\frac{{\rm log}(1+r(z)\overline{r(\bar{z})})}{z}{\rm d}z\Big]+{\rm exp}\Big[{\frac{1}{2\pi{i}}}\int^{0}_{\varepsilon}\frac{{\rm log}(1+r(z)\overline{r(\bar{z})})}{z}{\rm d}z\Big],$$
where $\varepsilon$ is a small number.

Using the asymptotic of the eigenfunctions, we learn that $r(z)=O(\frac{1}{z})\ {\rm as}\ z\rightarrow\infty,$
so  $${\rm exp}\Big[{\frac{1}{2\pi{i}}}\int^{\varepsilon}_{-\infty}\frac{{\rm log}(1+r(z)\overline{r(\bar{z})})}{z}{\rm d}z\Big]={\rm exp}\Big[{\frac{1}{2\pi{i}}}\int^{\varepsilon}_{-\infty}\frac{{\rm log}(1+|r(z)|^2)}{z}{\rm d}z\Big]$$ converges.

Considering that $r(0)=0$, so
$$\frac{{\rm d}[{\rm log}(1+|r(z)|^2)]}{{\rm d}z}\mid_{z=0}=\lim_{z\rightarrow0}\frac{{\rm log}(1+|r(z)|^2)-0}{z-0},$$
so $\int^{0}_{\varepsilon}\frac{{\rm log}(1+|r(z)|^2)}{z}{\rm d}z$ convergences when $\frac{{\rm d}[{\rm log}(1+|r(z)|^2)]}{{\rm d}z}|_{z=0}$ is bounded.

\end{proof}
\subsubsection{Removing the jump matrix $B$ in Region 2}
In this condition, we need to remove the jump on $(k_2,-\infty)$ and $(k_1,+\infty)$.
We consider a scalar Riemann-Hilbert problem:
\begin{itemize}
  \item $\delta_{j}(k) (j=3,4)$ is meromorphic in $\mathbb{C}\setminus\mathbb{R}$,
  \item $\begin{cases}
  \delta_{j,+}(k)=\delta_{j,-}(k)(1+r(k)\overline{r(\bar{k})}),
  \begin{cases}
  k\in(k_2,-\infty)\ {\rm for}\ j=3, \\
  k\in(k_1,+\infty)\ {\rm for}\ j=4,
  \end{cases} \\
  \delta_{j,+}(k)=\delta_{j,-}(k),   \begin{cases}
  k\in(k_2,+\infty)\ {\rm for}\ j=3, \\
  k\in(k_1,-\infty)\ {\rm for}\ j=4,
  \end{cases}
  \end{cases}$
  \item $\delta_j(k)\rightarrow1$, $k\rightarrow\infty$.
\end{itemize}
The Riemann-Hilbert problem exists the solution
\begin{equation}\label{delta_3}
\delta_j(k)=
\begin{cases}
{\rm exp}\Big[{\frac{1}{2\pi{i}}}\int_{k_2}^{-\infty}\frac{{\rm ln}(1+r(z)\overline{r(\bar{z})})}{z-k}{\rm d}z\Big], & {\rm for} j=3, \\
{\rm exp}\Big[{\frac{1}{2\pi{i}}}\int_{k_1}^{+\infty}\frac{{\rm ln}(1+r(z)\overline{r(\bar{z})})}{z-k}{\rm d}z\Big], & {\rm for} j=4,
\end{cases}
\end{equation}
and $\delta_3(k)=\delta_4(-k)$.
We can use $\left(
              \begin{array}{cc}
                \frac{1}{\delta_3(k)} & 0 \\
                0 & \delta_3(k) \\
              \end{array}
            \right)$ to remove the jump on $(k_2,-\infty)$,
and we can use $\left(
              \begin{array}{cc}
                \frac{1}{\delta_4(k)} & 0 \\
                0 & \delta_4(k) \\
              \end{array}
            \right)$ to remove the $\Delta_{(2)}$ can be used to remove the jump on $(k_1,+\infty)$.
It is worth noting that $\delta(k)=\delta_3(k)\delta_3(-k)$ has singular points $k=k_1$ and $k=k_2$, so the jump contour should be away from singularities.

\begin{figure}[h]
  \centering
  \includegraphics[width=2.6in,height=1.5in]{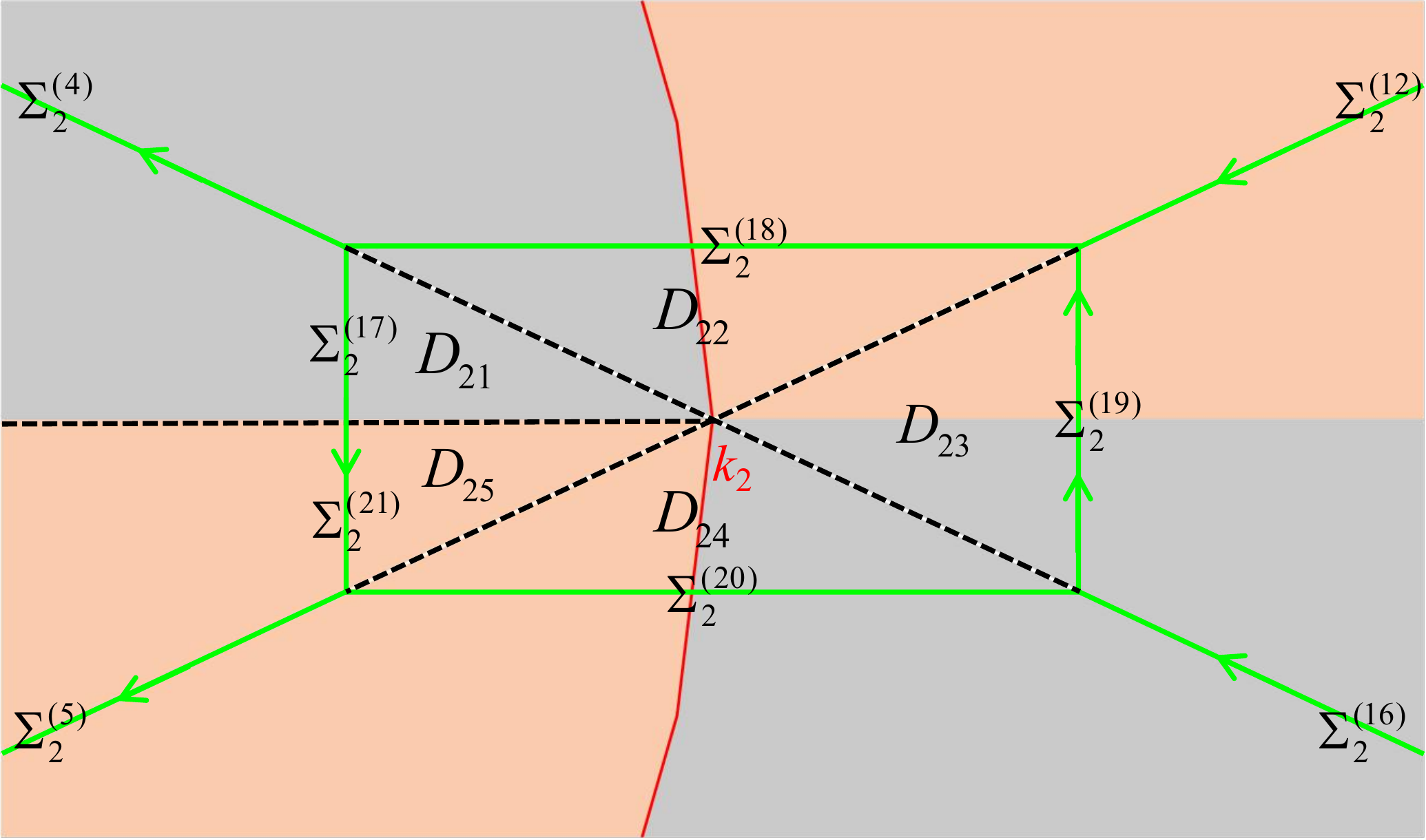}
  \includegraphics[width=2.6in,height=1.5in]{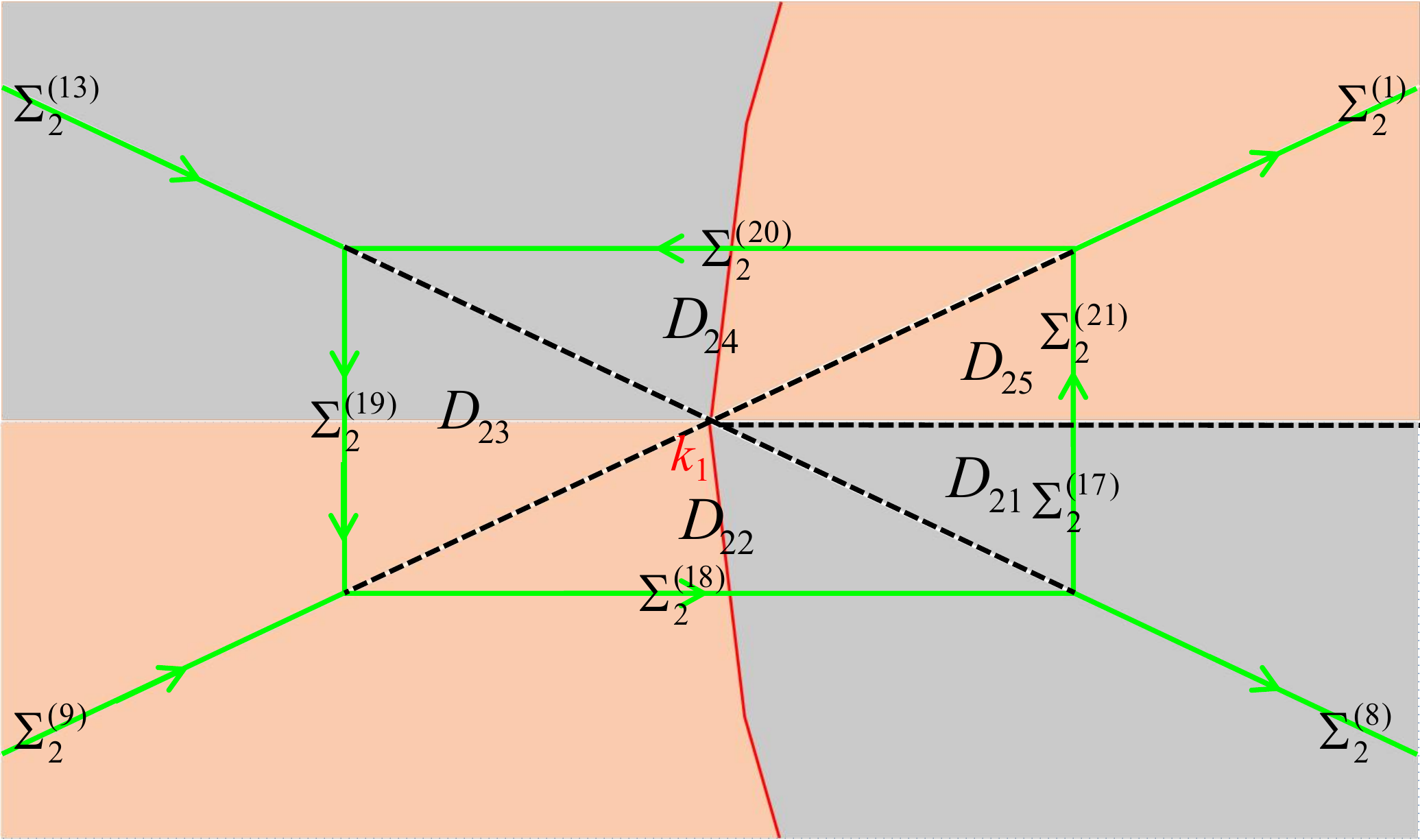}
  \caption{The jump contours for the \textbf{Riemann-Hilbert problem \uppercase\expandafter{\romannumeral6}} in the complex-$k$ plane.}\label{delta3_new}
\end{figure}

We define $M^{(5)}(x,t,k)$ by
\begin{equation}
	M^{(5)}(x,t,k)=M^{(2)}(x,t,k)
	\begin{cases}
		I
		, &\quad  {\rm if}\ k\in D_{21},
		\\
		A, &\quad {\rm if}\ k\in D_{22},
		\\
		L^{-1}A, &\quad {\rm if}\ k\in D_{23},
		\\
		C^{-1}B^{-1}, &\quad {\rm if}\ k\in D_{24},
		\\
		B^{-1}, &\quad {\rm if}\ k\in D_{25}.
		\\
		\Delta_2^{-1},  &\quad {\rm if}\ k\in {\rm others},
	\end{cases}
\end{equation}
where $\Delta_2=\left(
 \begin{array}{cc}
\frac{1}{\delta_3(k)\delta_3(-k)}  &0 \\
 0 & \delta_3(k)\delta_3(-k)\\
 \end{array}
 \right)$, the area $D_{2i} (i=1,2, ... , 5)$ is shown in Figure \ref{delta3_new}.

We get \textbf{Riemann-Hilbert problem \uppercase\expandafter{\romannumeral5}}:
\begin{itemize}
  \vspace{-0.0cm}
\item $M^{(5)}(x,t,k)$ is meromorphic in $\mathbb{C}\setminus\Sigma_5$,  $\Sigma_5$=$\Sigma_{2}^{(1)}\cup \cdots \cup\Sigma_{2}^{(21)}$,
\vspace{-0.0cm}
\item $M^{(5)}_{+}(x,t,k)=M^{(5)}_{-}(x,t,k)T^{(5)}(x,t,k)$, if $k\in\Sigma_5$,
\vspace{-0.0cm}
\item $M^{(5)}(x,t,k)\rightarrow I$, $k\rightarrow\infty$,
\end{itemize}
where $\Sigma_{2}^{(1)},\cdots, \Sigma_{2}^{(21)}$ are shown in Figure \ref{jump_t1} and Figure \ref{delta3_new},
\begin{equation}
	T^{(5)}(x,t,k)=
	\begin{cases}
\Delta_2C\Delta_2^{-1},  &{\rm if}\ k\in\Sigma^{(1)}_2\cup\Sigma^{(5)}_2,\\
\Delta_2U\Delta_2^{-1}, &{\rm if}\ k\in\Sigma^{(2)}_2\cup\Sigma^{(6)}_2\cup\Sigma^{(13)}_2
\cup\Sigma^{(14)}_2\cup\Sigma^{(15)}_2\cup\Sigma^{(16)}_2,\\
\Delta_2L\Delta_2^{-1}, &{\rm if}\ k\in \Sigma^{(3)}_2\cup\Sigma^{(7)}_2\cup\Sigma^{(9)}_2
\cup\Sigma^{(10)}_2\cup\Sigma^{(11)}_2\cup\Sigma^{(12)}_2,\\
\Delta_2A\Delta_2^{-1}, &{\rm if}\ k\in\Sigma^{(4)}_2\cup\Sigma^{(8)}_2,\\
\Delta_2, \quad &{\rm if}\ k\in \Sigma_2^{(17)}, \\
\Delta_2A, \quad &{\rm if}\ k\in \Sigma_2^{(18)}, \\
\Delta_2L^{-1}A, \quad &{\rm if}\ k\in \Sigma_2^{(19)}, \\
\Delta_2C^{-1}B^{-1}, \quad &{\rm if}\ k\in \Sigma_2^{(20)}, \\
\Delta_2B^{-1}, \quad &{\rm if}\ k\in \Sigma_2^{(21)}.
	\end{cases}
\end{equation}

\subsubsection{Removing the jump matrix $B$ in Region 3}
In this condition, we need to remove the jump on $(-\infty, 0)$, $(+\infty, 0)$, $(0, k_1)$ and $(0,k_2)$.
First, we remove the jump on $(0, k_1)$ and $(0, k_2)$.
Then, we remove the jump on $(-\infty, 0)$ and $(+\infty, 0)$, the method for removing the jump on $(-\infty, 0)$ and $(+\infty, 0)$ is introduced in section \ref{remove_1}.

In order to remove the jump on $(0, k_1)$ and $(0, k_2)$, we need two functions
$$\delta_5(k)={\rm exp}\Big[{\frac{1}{2\pi{i}}}\int^{k_2}_{0}\frac{{\rm log}(1+r(z)\overline{r(\bar{z})})}{z-k}{\rm d}z\Big],
\delta_6(k)={\rm exp}\Big[{\frac{1}{2\pi{i}}}\int^{k_1}_{0}\frac{{\rm log}(1+r(z)\overline{r(\bar{z})})}{z-k}{\rm d}z\Big].$$

\begin{figure}[h]
  \centering
  \includegraphics[width=2.6in,height=1.5in]{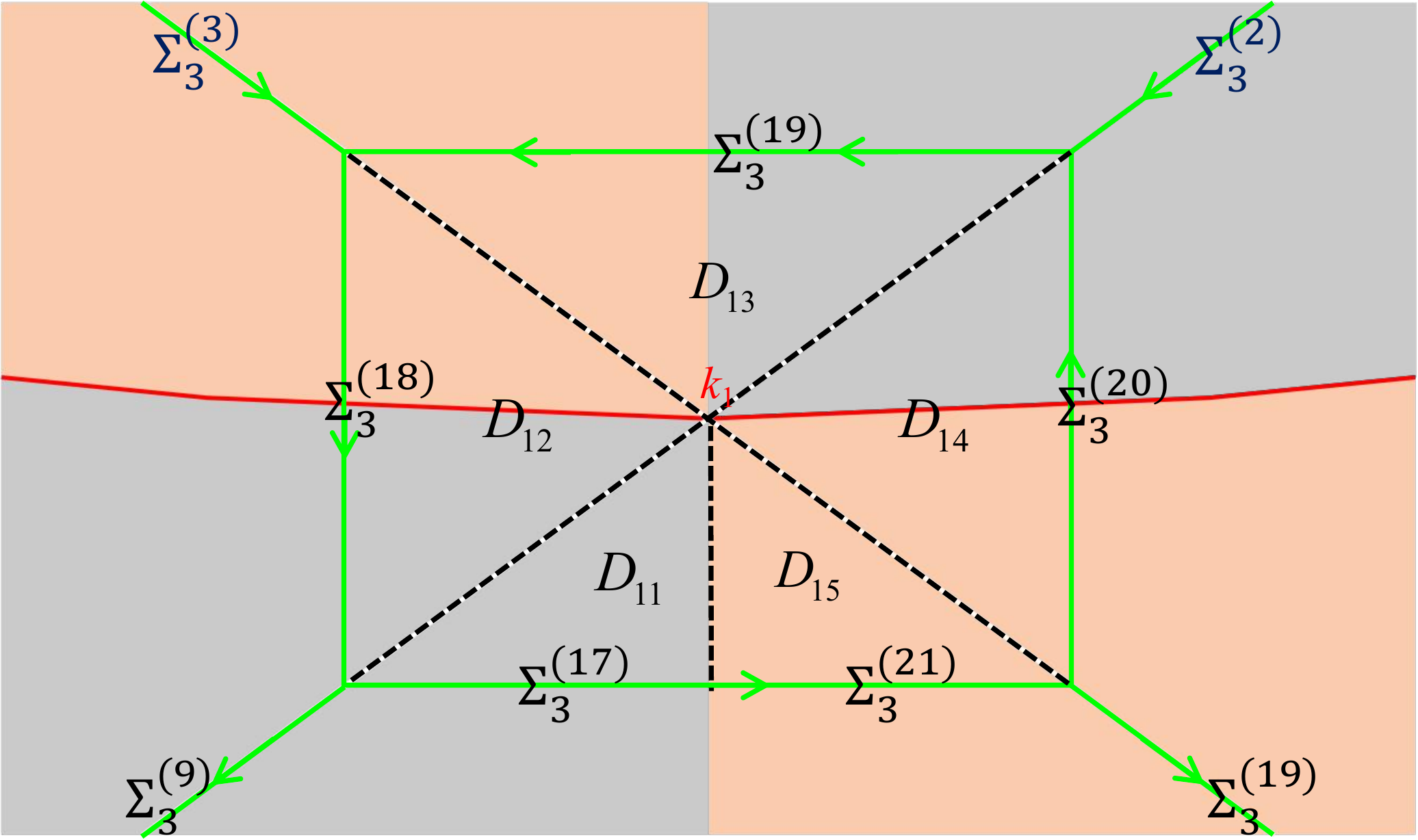}
  \includegraphics[width=2.6in,height=1.5in]{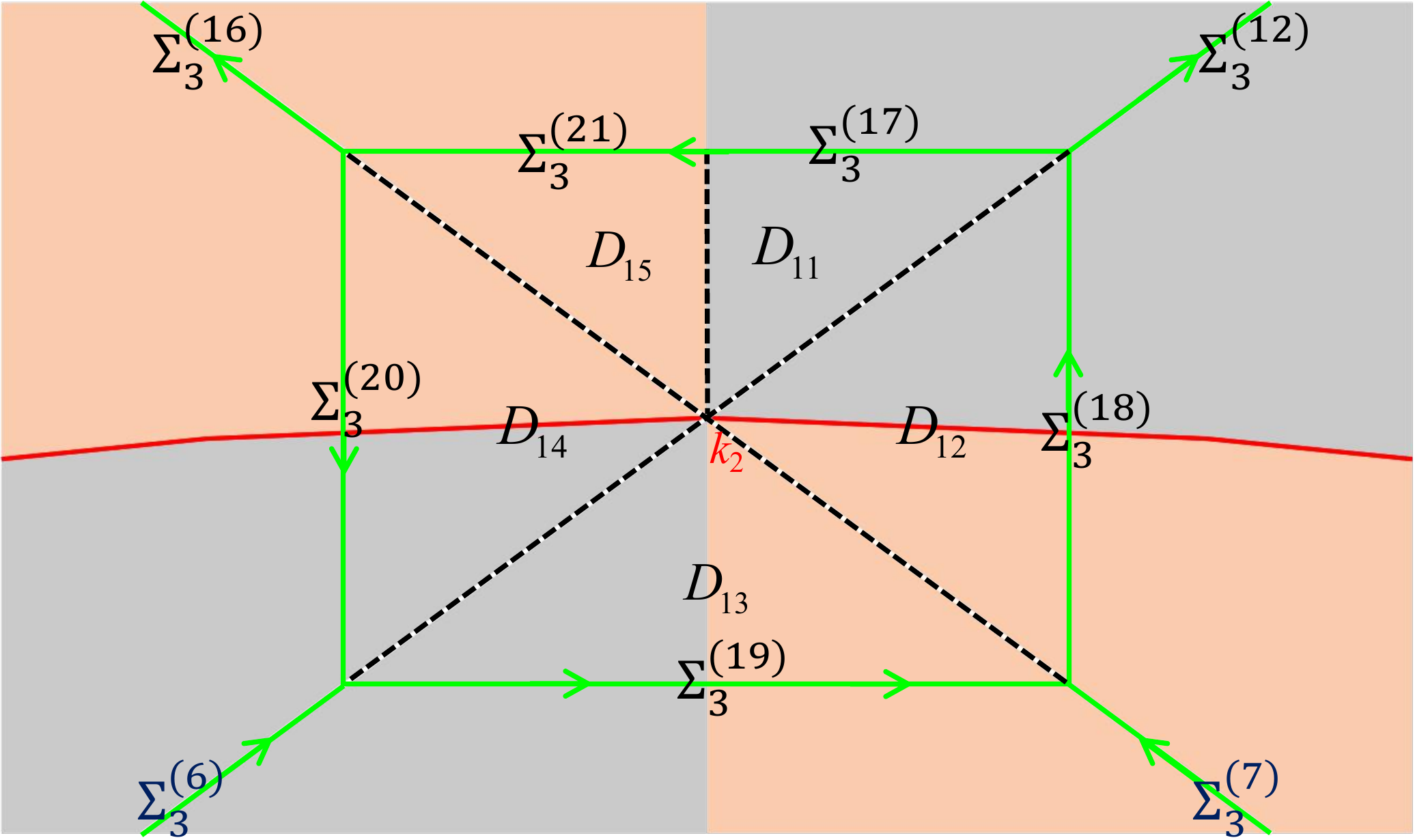}
  \caption{The jump contours for the \textbf{Riemann-Hilbert problem \uppercase\expandafter{\romannumeral6}} in the complex-$k$ plane.}\label{delta3}
\end{figure}

First, we define $M^{(6)}(x,t,k)$ by
\begin{equation}
	M^{(6)}(x,t,k)=M^{(3)}(x,t,k)
	\begin{cases}
		I
		, &\quad  {\rm if}\ k\in D_{11},
		\\
		A, &\quad {\rm if}\ k\in D_{12},
		\\
		L^{-1}A, &\quad {\rm if}\ k\in D_{13},
		\\
		C^{-1}B^{-1}, &\quad {\rm if}\ k\in D_{14},
		\\
		B^{-1}, &\quad {\rm if}\ k\in D_{15}.
		\\
		\Delta_3^{-1},  &\quad {\rm if}\ k\in {\rm others},
	\end{cases}
\end{equation}
where $\Delta_3=\left(
 \begin{array}{cc}
 \frac{1}{\delta_5(k)\delta_5(-k)}  &0 \\
 0 & \delta_5(k)\delta_5(-k) \\
 \end{array}
 \right)$, and the area $D_{1i} (i=1,2, ... , 5)$ is shown in Figure \ref{delta3}.

Then, we do a transformation $$M^{(6)}(x,t,k)=M^{(6)}(x,t,k)\Delta_1^{-1}.$$
We get \textbf{Riemann-Hilbert problem \uppercase\expandafter{\romannumeral6}}:
\begin{itemize}
  \vspace{-0.0cm}
\item $M^{(6)}(x,t,k)$ is meromorphic in $\mathbb{C}\setminus\Sigma_6$,  $\Sigma_6$=$\Sigma_{3}^{(1)}\cup \cdots \cup\Sigma_{3}^{(21)}$,
\vspace{-0.0cm}
\item $M^{(6)}_{+}(x,t,k)=M^{(6)}_{-}(x,t,k)T^{(6)}(x,t,k)$, if $k\in\Sigma_6$,
\vspace{-0.0cm}
\item $M^{(6)}(x,t,k)\rightarrow I$, $k\rightarrow\infty$,
\end{itemize}
where $\Sigma_{3}^{(1)},\cdots, \Sigma_{3}^{(21)}$ are shown in Figure \ref{jump_t2} and Figure \ref{delta3},
\begin{equation}
	T^{(6)}(x,t,k)=
	\begin{cases}
\Delta_1\Delta_3U\Delta_3^{-1}\Delta_1^{-1},  &{\rm if}\ k\in\Sigma^{(2)}_3\cup\Sigma^{(6)}_3,\\
\Delta_1\Delta_3C\Delta_3^{-1}\Delta_1^{-1}, &{\rm if}\ k\in\Sigma^{(1)}_3\cup\Sigma^{(5)}_3\cup\Sigma^{(13)}_3
\cup\Sigma^{(14)}_3\cup\Sigma^{(15)}_3\cup\Sigma^{(16)}_3,\\
\Delta_1\Delta_3A\Delta_3^{-1}\Delta_1^{-1}, &{\rm if}\ k\in \Sigma^{(4)}_3\cup\Sigma^{(8)}_3\cup\Sigma^{(9)}_3
\cup\Sigma^{(10)}_3\cup\Sigma^{(11)}_3\cup\Sigma^{(12)}_3,\\
\Delta_1\Delta_3L\Delta_3^{-1}\Delta_1^{-1}, &{\rm if}\ k\in\Sigma^{(3)}_3\cup\Sigma^{(7)}_3,\\
\Delta_1\Delta_3\Delta_1^{-1}, \quad &{\rm if}\ k\in \Sigma_3^{(17)}, \\
\Delta_1\Delta_3A\Delta_1^{-1}, \quad &{\rm if}\ k\in \Sigma_3^{(18)}, \\
\Delta_1\Delta_3L^{-1}A\Delta_1^{-1}, \quad &{\rm if}\ k\in \Sigma_3^{(19)}, \\
\Delta_1\Delta_3C^{-1}B^{-1}\Delta_1^{-1}, \quad &{\rm if}\ k\in \Sigma_3^{(20)}, \\
\Delta_1\Delta_3B^{-1}\Delta_1^{-1}, \quad &{\rm if}\ k\in \Sigma_3^{(21)}.
	\end{cases}
\end{equation}

\subsection{Deformations for the RHP (the soliton case)}\label{deform_2_1}
Lastly, we deal with the addition of solitons for the case of the DNLS equation.
We need to transform the jump conditions into the jump conditions.
We assume that the set $\{\kappa_j\}(j=1,2, ...)$ is bounded away from ${i}\mathbb{R}\cup\mathbb{R}$.

Fix $\varepsilon>0$, so that the circles $A^{(1)}_j=\{k\in\mathbb{C}: |k-\kappa_j|=\varepsilon\}$, $A^{(2)}_j=\{k\in\mathbb{C}: |k-\bar{\kappa}_j|=\varepsilon\}$, $A^{(3)}_j=\{k\in\mathbb{C}: |k+\bar{\kappa}_j|=\varepsilon\}$ and $A^{(4)}_j=\{k\in\mathbb{C}: |k+\kappa_j|=\varepsilon\}$ do not intersect each other and $i\mathbb{R}\cup\mathbb{R}$.

We define $N(x,t,k)$ by \begin{equation}
N(x,t,k)=M(x,t,k)
\begin{cases}
\left(
                   \begin{array}{cc}
                     1 & \frac{C_j{\rm e}^{-2i\theta{(\kappa_j)}}}{k-\kappa_j} \\
                     0 & 1 \\
                   \end{array}
                 \right), & {\rm if}\ |k-\kappa_j|<\varepsilon, \\
\left(
\begin{array}{cc}
1 & 0 \\
\frac{\bar{C}_j{\rm e}^{{2i}\theta{(\bar{\kappa}_j)}}}{k-\bar{\kappa}_j} & 1 \\
\end{array}
\right),  &{\rm if}\ |k-\bar{\kappa}_j|<\varepsilon, \\
\left(
\begin{array}{cc}
1 & 0\\
\frac{\bar{C}_j{\rm e}^{{2i}\theta{(-\bar{\kappa}_j)}}}{k+\bar{\kappa}_j}  & 1 \\
\end{array}
\right),  &{\rm if}\ |k+\bar{\kappa}_j|<\varepsilon, \\
\left(
                   \begin{array}{cc}
                     1 & \frac{C_j{\rm e}^{-{2i}\theta{(-\kappa_j)}}}{k+\kappa_j}  \\
                     0 & 1 \\
                   \end{array}
                 \right), & {\rm if}\ |k+\kappa_j|<\varepsilon, \\
I, &{\rm if}\ k\in {\rm others},
\end{cases}
\end{equation}
and the pole conditions are replaced by the jump conditions
\begin{equation}
N_{+}(x,t,k)=N_{-}(x,t,k)
\begin{cases}
\left(
                   \begin{array}{cc}
                     1 &  \frac{C_j{\rm e}^{-{2i}\theta{(\kappa_j)}}}{k-\kappa_j} \\
                     0 & 1 \\
                   \end{array}
                 \right), &k\in A^{(1)}_j, \\
\left(
\begin{array}{cc}
1 & 0 \\
\frac{\bar{C}_j{\rm e}^{{2i}\theta{(\bar{\kappa}_j)}}}{k-\bar{\kappa}_j} & 1 \\
\end{array}
\right), &k\in A^{(2)}_j, \\
\left(
\begin{array}{cc}
1 & 0 \\
\frac{\bar{C}_j{\rm e}^{{2i}\theta{(-\bar{\kappa}_j)}}}{k+\bar{\kappa}_j} & 1 \\
\end{array}
\right),  &k\in A^{(3)}_j, \\
\left(
                   \begin{array}{cc}
                     1 & \frac{C_j{\rm e}^{-{2i}\theta{(-\kappa_j)}}}{k+\kappa_j} \\
                     0 & 1 \\
                   \end{array}
                 \right), &k\in A^{(4)}_j,
\end{cases}
\end{equation}
where $A^{(1)}_j$ and $ A^{(4)}_j$ have counter-clockwise orientations, $A^{(2)}_j$ and $ A^{(3)}_j$ have clockwise orientations.

We get \textbf{Riemann-Hilbert problem N0}
\begin{itemize}
  \item $N(x,t,k)$ is meromorphic in $\mathbb{C}\setminus( i\mathbb{R}\cup \mathbb{R}\cup A^{(1)}_j\cup A^{(2)}_j \cup A^{(3)}_j \cup A^{(4)}_j )$,
  \item $N_{+}(x,t,k)=N_{-}(x,t,k)\begin{cases}
T, z\in i\mathbb{R}\cup\mathbb{R},  \\
U^{1}_j,  k\in A^{(1)}_j, \\
U^{2}_j,  k\in A^{(2)}_j, \\
U^{3}_j,  k\in A^{(3)}_j, \\
U^{4}_j,  k\in A^{(4)}_j,
\end{cases}$
  \item $N(x,t,k)\rightarrow I$, $k\rightarrow\infty$,
\end{itemize}
where $U^{1}_j=\left(
                   \begin{array}{cc}
                     1 & \frac{C_j{\rm e}^{-{2i}\theta{(\kappa_j)}}}{k-\kappa_j} \\
                     0 & 1 \\
                   \end{array}
                 \right), $
                 $U^{2}_j=\left(
\begin{array}{cc}
1 & 0 \\
\frac{\bar{C}_j{\rm e}^{{2i}\theta{(\bar{\kappa}_j)}}}{k-\bar{\kappa}_j} & 1 \\
\end{array}
\right)$,
$U^{3}_j=\left(
\begin{array}{cc}
1 & 0 \\
\frac{\bar{C}_j{\rm e}^{{2i}\theta{(-\bar{\kappa}_j)}}}{k+\bar{\kappa}_j} & 1 \\
\end{array}
\right)$ and
$U^{4}_j=\left(
                   \begin{array}{cc}
                     1 & \frac{C_j{\rm e}^{-{2i}\theta{(-\kappa_j)}}}{k+\kappa_j} \\
                     0 & 1 \\
                   \end{array}
                 \right).$

The jump matrix $T$ is oscillatory in ${i}\mathbb{R}\cup\mathbb{R}$, so we need to deform the jump contours of \textbf{Riemann-Hilbert problem N0}.
In soliton cases, the method for deforming the jump contours is the same as the non-soliton cases.
We need to ensure that the circles $A^{(m)}_j (m=1,2,3,4)$ do not intersect with any jump contours.
In soliton cases, we also divide the $(x,t)$-plane into three regions:
\begin{itemize}
  \item Region 1: $-c_1\sqrt{t}\leq x\leq c_1\sqrt{t}$.
  \item Region 2: $x<-c_1\sqrt{t}$.
  \item Region 3: $x>c_1\sqrt{t}$.
\end{itemize}

\subsubsection{Deformation in Region 1}
In Region 1, we directly deform the jump contours at $k=0$.
We define $N^{(1)}(x,t,k)$ by
\begin{equation}
N^{(1)}(x,t,k)=N(x,t,k)
\begin{cases}
I, & {\rm if}\ k\in \mathrm{E}_3\cup \mathrm{E}_6\cup  \mathrm{E}_9\cup \mathrm{E}_{12}, \\
C^{-1}, &  {\rm if}\ k\in \mathrm{E}_1\cup \mathrm{E}_7,\\
A, & {\rm if}\ k\in \mathrm{E}_2\cup \mathrm{E}_8,\\
L, & {\rm if}\ k\in \mathrm{E}_4\cup \mathrm{E}_{10},\\
U^{-1}, & {\rm if}\ k\in \mathrm{E}_5\cup \mathrm{E}_{11},
\end{cases}
\end{equation}
where $\mathrm{E}_{j}\ (j=1,2, ... , 12)$ is shown in Figure \ref{area_2}.

We get \textbf{Riemann-Hilbert problem N1}:
\begin{itemize}
\item $N^{(1)}(x,t,k)$ is meromorphic in $\mathbb{C}\setminus(\Sigma_{1}\cup A^{(0)}_j)$,
\item $N^{(1)}_{+}(x,t,k)=N^{(1)}_{-}(x,t,k)\begin{cases}
L, \quad &{\rm if}\ k\in \Sigma^{(3)}_1\cup\Sigma^{(7)}_1,\\
A, \quad &{\rm if}\ k\in \Sigma^{(4)}_1\cup\Sigma^{(8)}_1,\\
C, \quad &{\rm if}\ k\in \Sigma^{(1)}_1\cup\Sigma^{(5)}_1,\\
U, \quad &{\rm if}\ k\in \Sigma^{(2)}_1\cup\Sigma^{(6)}_1,\\
B, \quad &{\rm if}\ k\in \mathbb{R},\\
U^m_j,\quad &{\rm if}\ k\in A^{(m)}_j, (m=1,2,3,4),
\end{cases}$
  \item $N^{(1)}(x,t,k)\rightarrow I$, $k\rightarrow\infty$,
\end{itemize}
where $A^{(0)}_j=A^{(1)}_j\cup A^{(2)}_j\cup A^{(3)}_j\cup A^{(4)}_j$, $\Sigma^{(j)}_1(j=1,2,...,8)$ is shown in  Figure \ref{jump_t1}.
\subsubsection{Deformation in Region 2}
In Region 2, we need to deform the jump contours at the saddle points 0, $k_1$ and $k_2$.

We define $N^{(2)}(x,t,k)$ by
\begin{equation}
N^{(2)}(x,t,k)=N(x,t,k)
\begin{cases}
I, & {\rm if}\ k\in \mathrm{E}^2_3\cup \mathrm{E}^2_6\cup  \mathrm{E}^2_9\cup \mathrm{E}^2_{12}, \\
C^{-1}, &  {\rm if}\ k\in \mathrm{E}^2_1\cup\mathrm{E}^2_7\cup\mathrm{E}^2_{14}\cup\mathrm{E}^2_{16},\\
A, & {\rm if}\ k\in \mathrm{E}^2_2\cup \mathrm{E}^2_8\cup\mathrm{E}^2_{13}\cup\mathrm{E}^2_{15},\\
L, & {\rm if}\ k\in \mathrm{E}^2_4\cup \mathrm{E}^2_{10},\\
U^{-1}, & {\rm if}\ k\in \mathrm{E}^2_5\cup \mathrm{E}^2_{11},
\end{cases}
\end{equation}
where the area $\mathrm{E}^{2}_{j}\ (j=1,2, ... , 16)$ is shown in Figure \ref{area_t1}.

We get \textbf{Riemann-Hilbert problem N2}:\begin{itemize}
  \item $N^{(2)}(x,t,k)$ is meromorphic in $\mathbb{C}\setminus(\Sigma_2\cup A^{(0)}_j)$,
\item $N^{(2)}_{+}(x,t,k)=N^{(2)}_{-}(x,t,k)\begin{cases}
A, \quad &{\rm if}\ k\in \Sigma^{(1)}_2\cup\Sigma^{(5)}_2,\\
C, \quad &{\rm if}\ k\in \Sigma^{(4)}_2\cup\Sigma^{(8)}_2,\\
U, \quad &{\rm if}\ k\in \Sigma^{(2)}_2\cup\Sigma^{(6)}_2
\cup\Sigma^{(13)}_1\cup\Sigma^{(14)}_2\cup\Sigma^{(15)}_2\cup\Sigma^{(16)}_2,\\
L, \quad &{\rm if}\ k\in\Sigma^{(3)}_2\cup\Sigma^{(7)}_2
\cup\Sigma^{(9)}_1\cup\Sigma^{(10)}_2\cup\Sigma^{(11)}_2\cup\Sigma^{(12)}_2,\\
B, \quad &{\rm if}\ k\in\cup(-\infty,k_2)\cup(k_1,+\infty),\\
U^m_j,\quad &{\rm if}\ k\in A^{(m)}_j, (m=1,2,3,4),
\end{cases}$
  \item $N^{(2)}(x,t,k)\rightarrow I$,\qquad $k\rightarrow\infty$,
\end{itemize}
where the jump contour $\Sigma^{(j)}_2 (j=1,2, ... ,16)$ is shown in Figure \ref{jump_t1}.
\subsubsection{Deformation in Region 3}
In Region 3, we need to deform the jump contours at the saddle points $k=0$, $k=k_1$ and $k=k_2$.

We define $N^{(3)}(x,t,k)$ by
\begin{equation}
N^{(3)}(x,t,k)=N(x,t,k)
\begin{cases}
I,      &  {\rm if}\ k\in \mathrm{E}^3_3\cup \mathrm{E}^3_6\cup  \mathrm{E}^3_9\cup \mathrm{E}^3_{12}, \\
C^{-1}, &  {\rm if}\ k\in \mathrm{E}^3_1\cup\mathrm{E}^3_7\cup\mathrm{E}^3_{14}\cup\mathrm{E}^3_{16},\\
A,      &  {\rm if}\ k\in \mathrm{E}^3_2\cup \mathrm{E}^3_8\cup\mathrm{E}^3_{13}\cup\mathrm{E}^3_{15},\\
L,      &  {\rm if}\ k\in \mathrm{E}^3_4\cup \mathrm{E}^3_{10},\\
U^{-1}, &  {\rm if}\ k\in \mathrm{E}^3_5\cup \mathrm{E}^3_{11}, \\
\end{cases}
\end{equation}
where the area $\mathrm{E}^{3}_{j}\ (j=1,2, ... , 16)$ is shown in Figure \ref{area_t2}.

We get \textbf{Riemann-Hilbert problem N3}:
\begin{itemize}
  \vspace{-0.0cm}
  \item $N^{(3)}(x,t,k)$ is meromorphic in $\mathbb{C}\setminus(\Sigma_3\cup A^{(0)}_j)$, $\Sigma_3=\Sigma^{(1)}_3\cup \cdots \cup\Sigma^{(16)}_3\cup(-\infty{i}, k_2)\cup(\infty{i},k_1)$,
  \vspace{-0.0cm}
  \item $N^{(3)}_{+}(x,t,k)=N^{(3)}_{-}(x,t,k)\begin{cases}
L
, &  {\rm if}\ k\in\Sigma^{(3)}_3\cup\Sigma^{(7)}_3\cup\Sigma^{(9)}_3
\cup\Sigma^{(10)}_3\cup\Sigma^{(11)}_3\cup\Sigma^{(12)}_3,\\
U,&{\rm if}\ k\in\Sigma^{(4)}_3\cup\Sigma^{(8)}_3\cup\Sigma^{(13)}_3
\cup\Sigma^{(14)}_3\cup\Sigma^{(15)}_3\cup\Sigma^{(16)}_3,\\
C, & {\rm if}\ k\in\Sigma^{(1)}_3\cup\Sigma^{(5)}_3,\\
A, &{\rm if}\ k\in\Sigma^{(2)}_3\cup\Sigma^{(6)}_3,\\
B,&  {\rm if}\ k\in(-\infty{i}, k_2)\cup(\infty{i}, k_1), \\
U^m_j,\quad &{\rm if}\ k\in A^{(m)}_j, (m=1,2,3,4),
\end{cases}$
  \item $N^{(3)}(x,t,k)\rightarrow I$, $k\rightarrow\infty$,
\end{itemize}
where the jump contour $\Sigma^{(j)}_3 (j=1,2, ... ,16)$ is shown in Figure \ref{jump_t2}.
\subsection{Extra deformations for the RHP (soliton cases)}\label{deform_2_2}
When the initial value for the DNLS equation can evolve into the soliton, some extra deformations are needed.
\begin{remark}
For some $x$ and $t$, $|{\rm e}^{2i\theta{(\kappa_j)}}|$ and $|{\rm e}^{-2i\theta{(\bar{\kappa}_j)}}|$ may become large, an extra deformation is used to deal with $U^m_j$.
\end{remark}

We define $N^{(4)}(x,t,k)$ by
\begin{equation}
N^{(4)}(x,t,k)=N^{( m_1 )}(x,t,k)
\begin{cases}
\left(
                   \begin{array}{cc}
                     1 & \frac{k-\kappa_j}{C_j{\rm e}^{-{2i}\theta{(\kappa_j)}}} \\
                     -\frac{C_j{\rm e}^{-{2i}\theta{(\kappa_j)}}}{k-\kappa_j} & 0 \\
                   \end{array}
                 \right)K,&  {|k-\kappa_j|<\varepsilon}, \\
\left(
                   \begin{array}{cc}
                     1 & \frac{k+\kappa_j}{C_j{\rm e}^{2i\theta{(-\kappa_j)}}} \\
                   -\frac{C_j{\rm e}^{2i\theta{(-\kappa_j)}}}{k+\kappa_j} & 0 \\
                   \end{array}
                 \right)K,&  {|k+\kappa_j|<\varepsilon}, \\
\left(
\begin{array}{cc}
0 & \frac{\bar{C}_j{\rm e}^{2i\theta{(\bar{\kappa}_j)}}}{k-\bar{\kappa}_j} \\
-\frac{k-\bar{\kappa}_j} {\bar{C}_j{\rm e}^{2i\theta{(\bar{\kappa}_j)}}}& 1 \\
\end{array}
\right)K,&  {|k-\bar{\kappa}_j|<\varepsilon},  \\
\left(
\begin{array}{cc}
0 & \frac{\bar{C}_j{\rm e}^{-2i\theta{(-\bar{\kappa}_j)}}}{k+\bar{\kappa}_j} \\
-\frac{k+\bar{\kappa}_j} {\bar{C}_j{\rm e}^{-2i\theta{(-\bar{\kappa}_j)}}}& 1 \\
\end{array}
\right)K,&  {|k+\bar{\kappa}_j|<\varepsilon},  \\
K,& k\in {\rm others},
\end{cases}
\end{equation}
where $K=\left(
           \begin{array}{cc}
             \frac{(z-\kappa_j)(z+\kappa_j)}{(z-\bar{\kappa}_j)(z+\bar{\kappa}_j)} & 0 \\
             0 & \frac{(z-\bar{\kappa}_j)(z+\bar{\kappa}_j)}{(z-\kappa_j)(z+\kappa_j)} \\
           \end{array}
         \right)$ and $m_1$ is equal to 1, 2, or 3.

In Region 1, we get \textbf{Riemann-Hilbert problem N4}:
\begin{itemize}
\item $N^{(4)}(x,t,k)$ is meromorphic in $\mathbb{C}\setminus(\Sigma_{1}\cup A^{(0)}_j)$,
\item $N^{(4)}_{+}(x,t,k)=N^{(4)}_{-}(x,t,k)\begin{cases}
K^{-1}LK, \quad &{\rm if}\ k\in \Sigma^{(3)}_1\cup\Sigma^{(7)}_1,\\
K^{-1}AK, \quad &{\rm if}\ k\in \Sigma^{(4)}_1\cup\Sigma^{(8)}_1,\\
K^{-1}CK, \quad &{\rm if}\ k\in \Sigma^{(1)}_1\cup\Sigma^{(5)}_1,\\
K^{-1}UK, \quad &{\rm if}\ k\in \Sigma^{(2)}_1\cup\Sigma^{(6)}_1,\\
K^{-1}BK, \quad &{\rm if}\ k\in i\mathbb{R},\\
K^{-1}V^m_jK,\quad &{\rm if}\ k\in A^{(m)}_j, (m=1,2,3,4),
\end{cases}$
  \item $N^{(4)}(x,t,k)\rightarrow I$, $k\rightarrow\infty$,
\end{itemize}
where $V^{1}_j=\left(
                   \begin{array}{cc}
                     1 & 0 \\
                     \frac{k-\kappa_j}{C_j{\rm e}^{2i\theta{(\kappa_j)}}} & 1 \\
                   \end{array}
                 \right), $
$V^{2}_j=\left(
\begin{array}{cc}
1 & \frac{k-\bar{\kappa}_j}{\bar{C}_j{\rm e}^{-2i\theta{(\bar{\kappa}_j)}}} \\
0 & 1 \\
\end{array}
\right)$,
$V^{3}_j=\left(
\begin{array}{cc}
1 & \frac{k+\bar{\kappa}_j}{\bar{C}_j{\rm e}^{-2i\theta{(-\bar{\kappa}_j)}}} \\
0 & 1 \\
\end{array}
\right)$ and
$V^{4}_j=\left(
                   \begin{array}{cc}
                     1 & 0 \\
                     \frac{k+\kappa_j}{C_j{\rm e}^{2i\theta{(-\kappa_j)}}} & 1 \\
                   \end{array}
                 \right)$.

The Riemann-Hilbert problem of Region 2 or Region 3 is similar to Region 1, and we do not give its specific expressions.
\begin{remark}
When $t$ is large, we also need to remove the jump matrix $B$, which is the same as methods of section \ref{deform_1_2}.
\end{remark}
\begin{remark}
When choosing the jump contours, we need to consider the analytical property of the scattering data, which have been widely applied in previous NIST\cite{nist_Trogdon1,Cui}.
\end{remark}
\subsection{Numerical solution for the RHP}\label{sec_cheby}
After deforming the Riemann-Hilbert problem, it is necessary to solve the deformed problem using a numerical method. The Chebyshev collocation method, developed by Olver, is commonly employed for solving the Riemann-Hilbert problem numerically\cite{Olver1,Olver2,Olver3}.

Give $x$ and $t$, and we get the following RHP,
\begin{equation}\label{rhp_fai}
\begin{array}{lr}
\Phi^{+}(k)=\Phi^{-}(k)G(k), k\in\Gamma, \\
\Phi(\infty)=I,
\end{array}
\end{equation}
where $\Gamma$=$\Gamma_1\cup\Gamma_2\cup\cdots\cup\Gamma_{n_1}$, and $\Gamma_1,\ldots, \Gamma_{n_1}$ are non-self-intersecting arcs.

$\Phi(k)$ can be written into the following form
\begin{equation}\label{phi_1}
\Phi(k)=I+\mathcal{C}_{\Gamma}V(k),
\end{equation}
where $\mathcal{C}_{\Gamma}$ represents the Cauchy transformation
$$\mathcal{C}_{\Gamma}V(z)=\frac{1}{2i\pi}\int_{\Gamma}\frac{V(\xi)}{\xi-z}{\rm d}\xi.$$
According to the Plemelj formula, we get
\begin{numcases}{}
\Phi^{+}(k)-\Phi^{-}(k)=V(k), \\
\Phi^{+}(k)+\Phi^{-}(k)=-i\mathcal{H}V(k),
\end{numcases}
where $\mathcal{H}$ is the Hilbert transform.

We need to search the suitable $\rm M\ddot{o}bius$ transform $M_j (j=1, 2, \ldots, n_1)$
which mapping $\Gamma_{j}$ into unit interval $\mathbb{I}$.

Eq. (\ref{phi_1}) is written into
\begin{equation}
\Phi(k)=I+\mathcal{C}_{\Gamma_1}V_1(k)+\mathcal{C}_{\Gamma_2}V_2(k)+\ldots+\mathcal{C}_{\Gamma_{n_1}}V_{n_1}(k).
\end{equation}

For $V_{j}(k)=\left(
         \begin{array}{cc}
           V^{(11)}_{j}(k) & V^{(12)}_{j}(k) \\
           V^{(21)}_{j}(k) & V^{(22)}_{j}(k) \\
         \end{array}
       \right),$ we appropriate it by the Chebyshev polynomials
\begin{equation}\label{V_series}
\begin{cases}
V^{(1,1)}_{j}(k)=T(M_{j}^{-1}(k))\mathcal{F}v^{(11)}_{j}, \\
V^{(1,2)}_{j}(k)=T(M_{j}^{-1}(k))\mathcal{F}v^{(12)}_{j}, \\
V^{(2,1)}_{j}(k)=T(M_{j}^{-1}(k))\mathcal{F}v^{(21)}_{j}, \\
V^{(2,2)}_{j}(k)=T(M_{j}^{-1}(k))\mathcal{F}v^{(22)}_{j},
\end{cases}
\end{equation}
where $v^{(\ )}_{j}$ is the unknown coefficient.

The unknown coefficient $v^{(\ )}_{j}$ can be calculated by Eq. (\ref{rhp_fai}).
According to Eq. (\ref{phi_1}),
we get
\begin{equation}
\begin{aligned}
&\lim_{k\rightarrow\infty}k\Phi(k)_{(1,2)}
=\lim_{k\rightarrow\infty}\frac{1}{2i\pi}\int_{\Gamma}\frac{k}{\xi-k}V^{(12)}(\xi){\rm d}\xi \\
&=-\frac{1}{2i\pi}\Big[\int_{\Gamma_1}V_1^{(12)}(\xi){\rm d}\xi+\int_{\Gamma_2}V_2^{(12)}(\xi){\rm d}\xi+\cdots+\int_{\Gamma_{k_1}}V_{k_1}^{(12)}(\xi){\rm d}\xi\Big].
\end{aligned}
\end{equation}

We only summarize the key steps of the Chebyshev collocation method and do not give specific details of the method.
The specific details of the Chebyshev collocation method can be found in Olver's paper\cite{Olver1,Olver2}.
We give the following remarks to help readers to implement this method better.
\begin{remark}
The keys for the Chebyshev collocation method are calculating the Cauchy transformation of the Chebyshev polynomials and constructing Cauchy matrices, which allows us to turn (\ref{rhp_fai}) into a linear system with unknown coefficient $v^{(\ )}_{j}$.
\end{remark}
\begin{remark}
There are some relationships between certain Cauchy matrices, perhaps they are completely equal, or their imaginary parts are equal and their real parts are opposite. These relationships are determined by the value of $k$ and the jumping contour $\Gamma_j(j=1,2,...,k_1)$. Using these relationships can reduce computational costs.

\end{remark}

\section{Numerical Results and Comparison}\label{sec_result1}
In this section, the NIST is used to solve the DNLS equation, and we compare the NIST with the traditional numerical method.
It is necessary to introduce traditional numerical methods before presenting numerical results.
One of the most effective methods for solving the DNLS equation is the Fourier spectral method (FSM).
In the FSM, the $x$-space is truncated into $[-L, L]$ and the periodic boundary $q(L,t)=q(-L,t)$ is applied.
Doing the Fourier transformation for the the DNLSE $\rm \uppercase\expandafter{\romannumeral3}$, we get
\begin{equation}\label{fft_1}
\mathbb{F}[q]_t=-ik^2\mathbb{F}[q]+\mathbb{F}\Big[q^2\mathbb{F}^{-1}\big[ ik\mathbb{F}[\bar{q}]\big]\Big]+\frac{ i}{2}\mathbb{F}[|q|^4q],
\end{equation}
where $\mathbb{F}$ represents the Fourier transformation and $\mathbb{F}^{-1}$ represents the inverse Fourier transformation.
We rewrite Eq. (\ref{fft_1}) into
\begin{equation}\label{fft_2}
\frac{{\rm d}{\rm e}^{ik^2t}\mathbb{F}[q]}{{\rm d}t}={\rm e}^{ik^2t}\mathbb{F}\Big[q^2\mathbb{F}^{-1}\big[ ik\mathbb{F}[\bar{q}]\big]\Big]+{\rm e}^{ik^2t}\frac{i}{2}\mathbb{F}[|q|^4q].
\end{equation}
Eq. (\ref{fft_2}) can be solved by the finite difference method (FDM) or the Runge-Kutta method (RKM). The solution for the DNLS equation is given by the inverse Fourier transformation  of $\mathbb{F}[q]$.

\begin{figure}[h]
  \centering
  \subfigure[$t=0$]{\includegraphics[width=1.4in,height=1.1in]{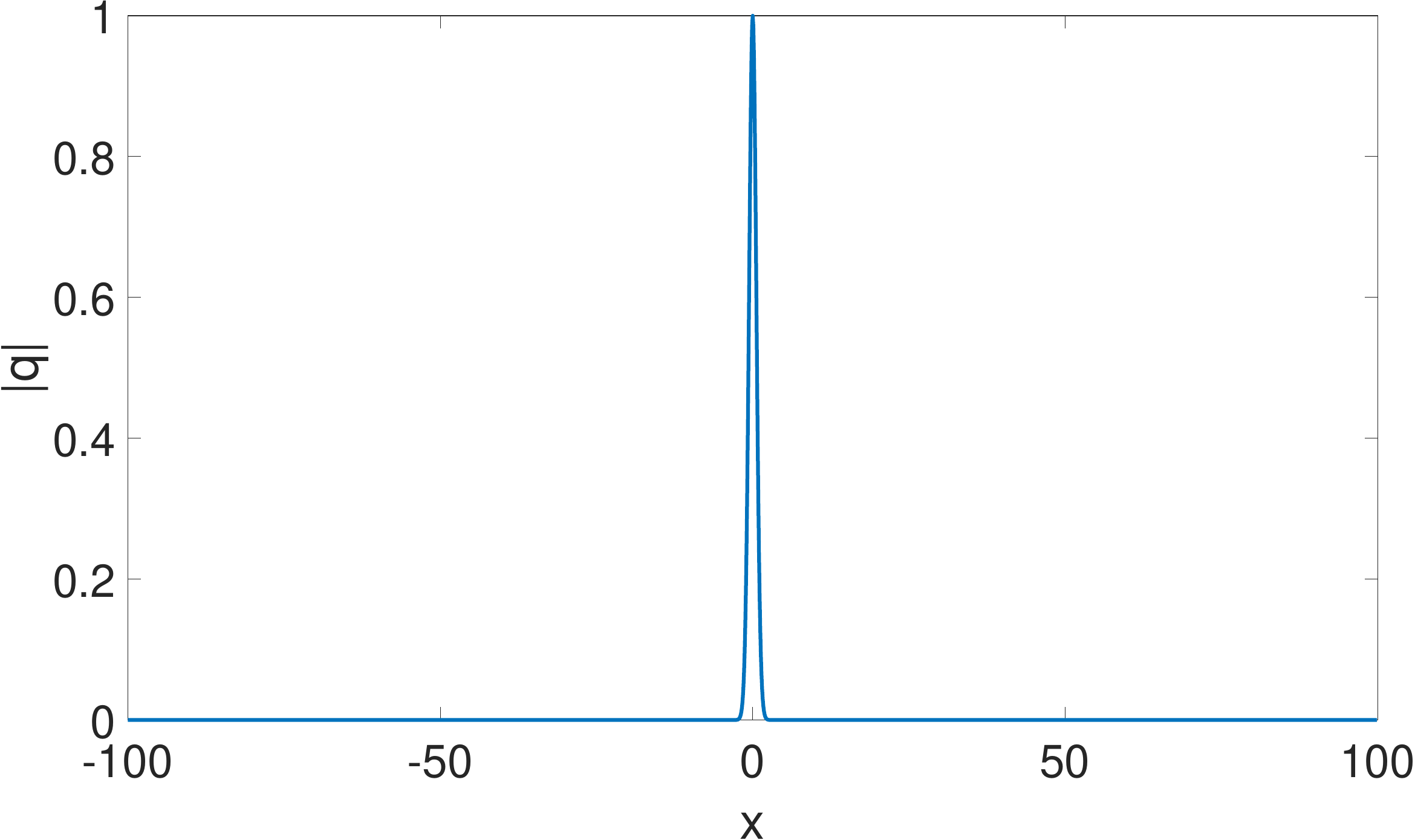}}
  \subfigure[$t=5$]{\includegraphics[width=1.4in,height=1.1in]{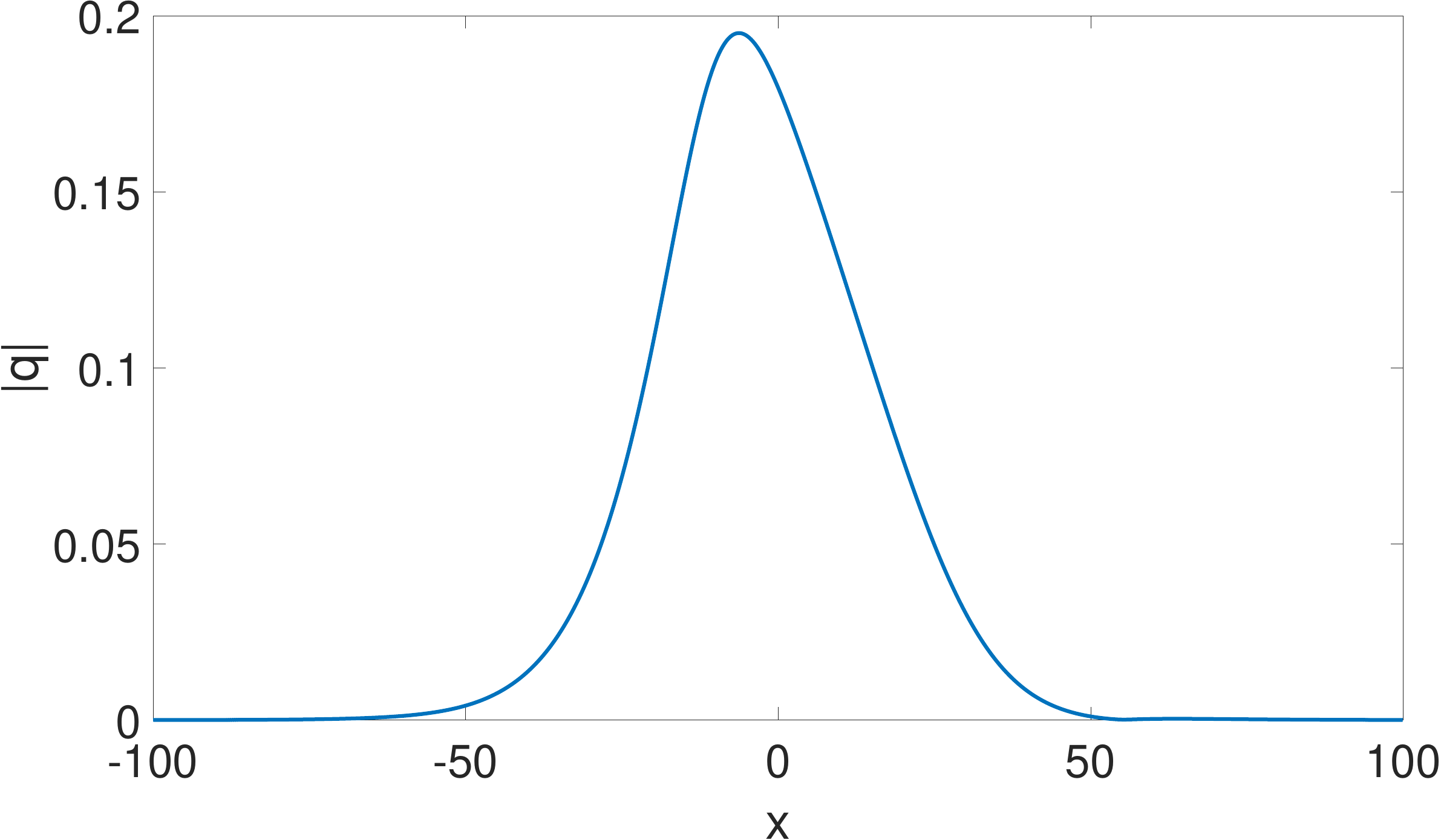}}
  \subfigure[$t=10$]{\includegraphics[width=1.4in,height=1.1in]{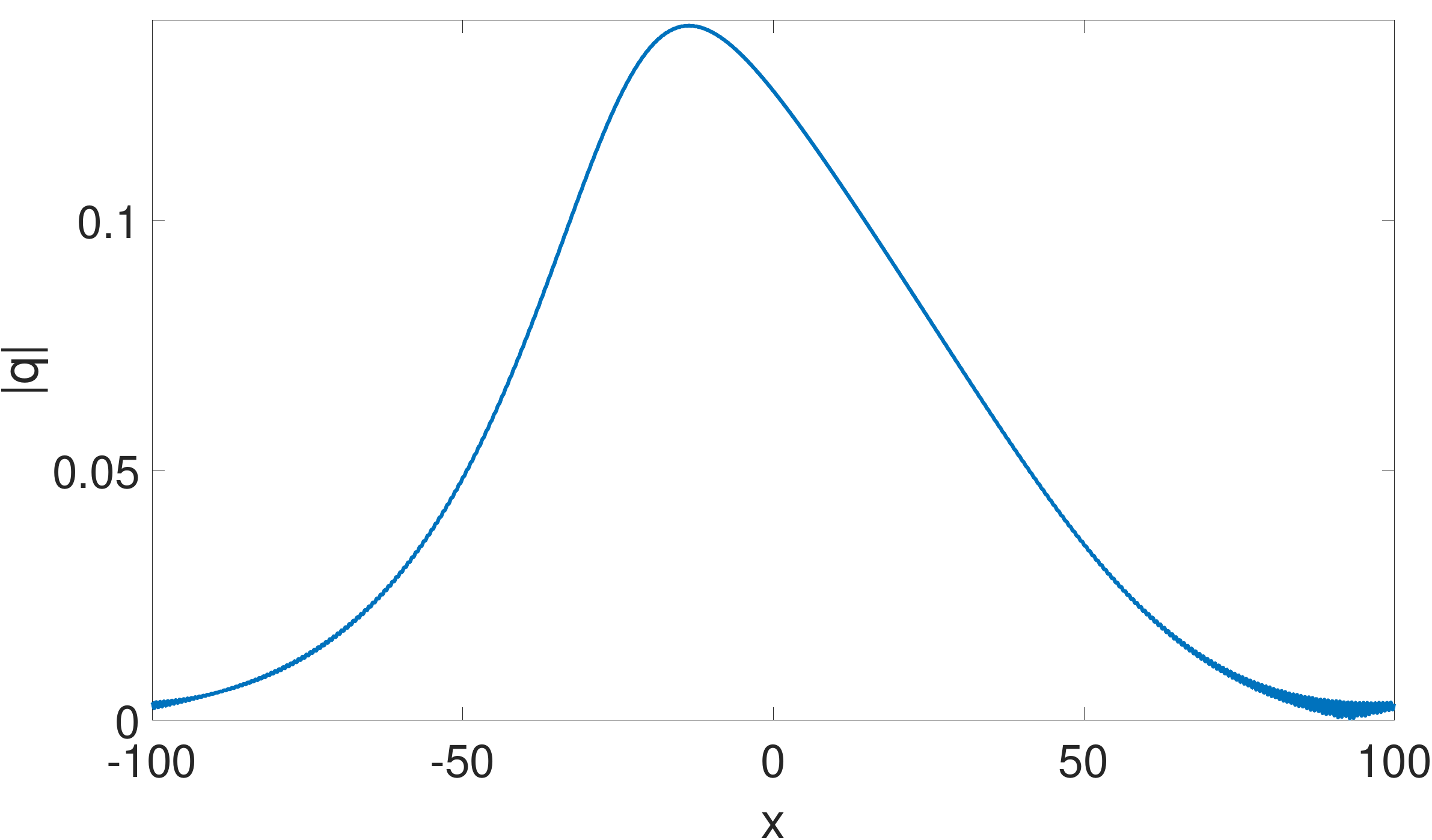}}
  \subfigure[$t=15$]{\includegraphics[width=1.4in,height=1.1in]{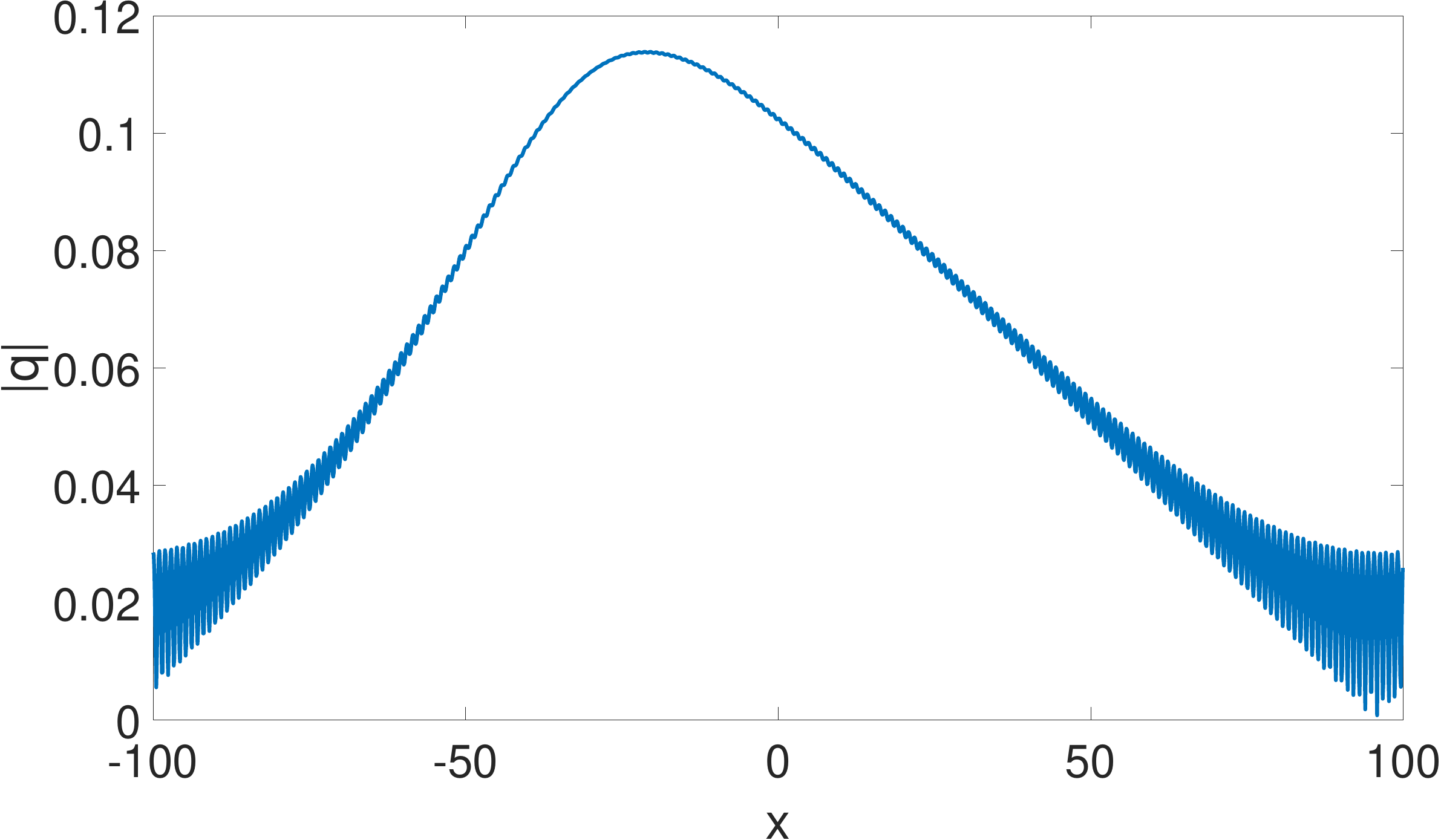}}
  \subfigure[$t=20$]{\includegraphics[width=1.4in,height=1.1in]{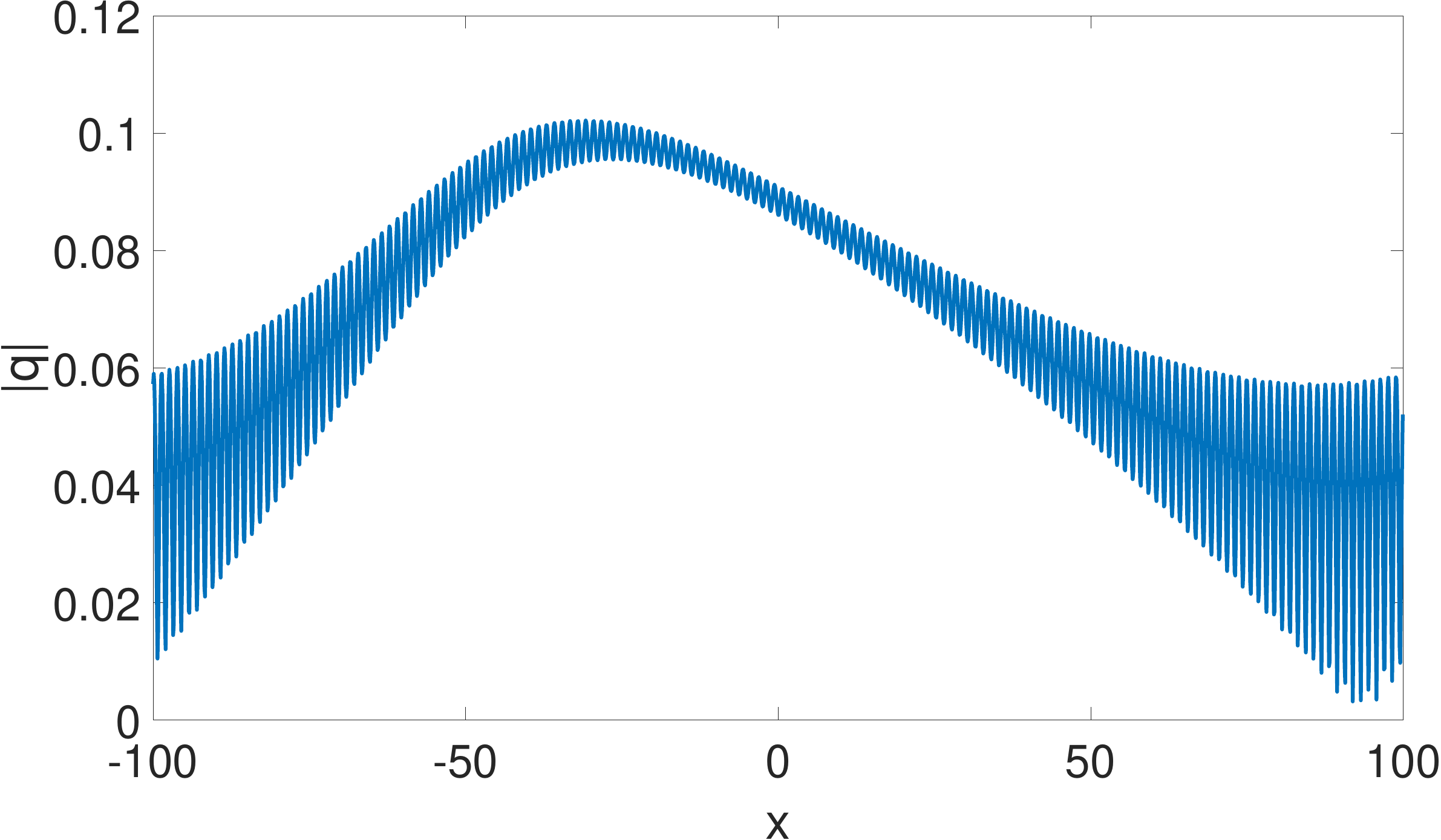}}
  \subfigure[$t=25$]{\includegraphics[width=1.4in,height=1.1in]{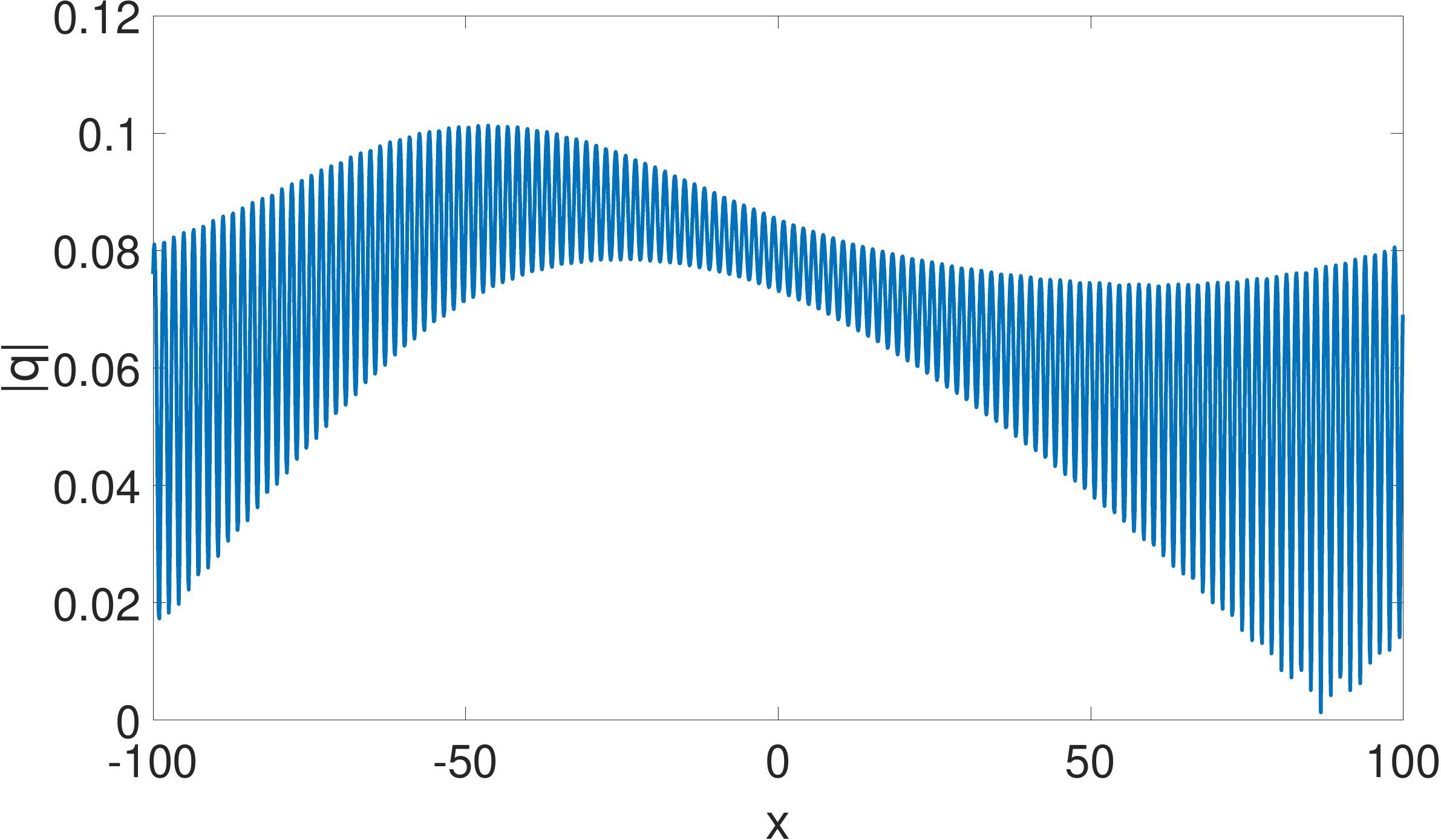}}
  \subfigure[$t=30$]{\includegraphics[width=1.4in,height=1.1in]{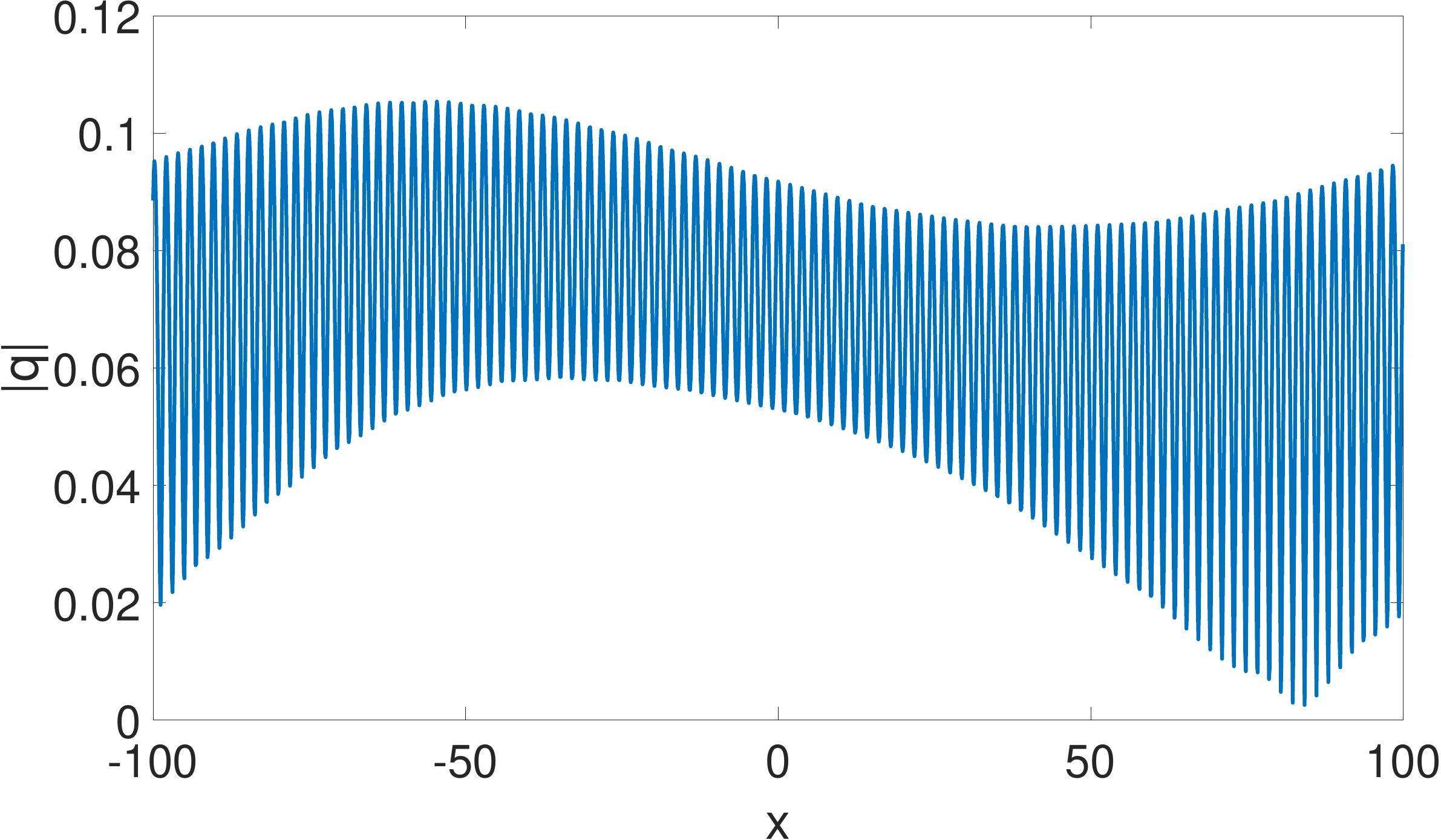}}
  \subfigure[$t=35$]{\includegraphics[width=1.4in,height=1.1in]{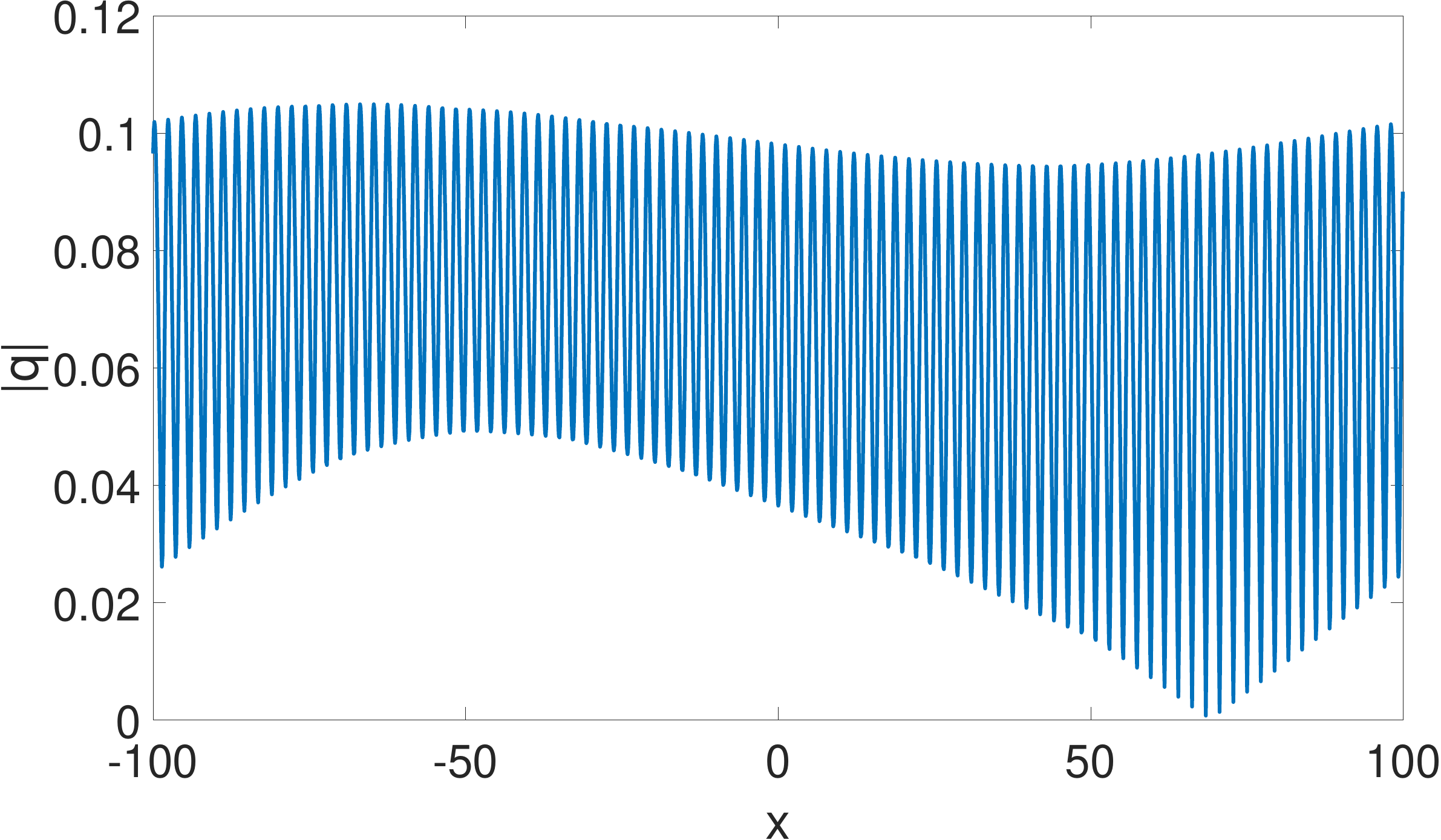}}
  \caption{The numerical solution for the $\exp(-x^2){\rm sech}(x)$ initial profile calculated by the FSM.}\label{fsm_ans1}
\end{figure}

There is no doubt that the FSM is efficient in the short time, but the FSM will lose accuracy once the wave arrives the calculated boundary. The error generated by the FSM will increase over time which will lead to numerical convergence issues. The FSM is used to calculate the DNLS equation with $\exp(-x^2){\rm sech}(x)$ initial profile, and the calculated results are shown in Figure \ref{fsm_ans1}.
In Figure \ref{fsm_ans1}(a)-(c), the calculated results are normal.
When the time reaches 15, there is a significant oscillation in the numerical results, and the oscillation becomes more severe as time increases.
When the time reaches 30, the FSM cannot effectively calculate the solution of the DNLS equation.
These phenomena can be clearly observed from Figure \ref{fsm_ans1} (d)-(h).

There is no discrete eigenvalue for the $\exp(-x^2){\rm sech}(x)$ initial profile, so we don't need to calculate the discrete eigenvalues and the normalized constant.
The NIST is used to calculate the evolution for $\exp(-x^2){\rm sech}(x)$ initial profile in the short time, and the calculated results are marked as $q_{\rm nist1}(x,t)$.

\begin{figure}[H]
  \centering
  \subfigure[Density]{\includegraphics[width=2.in,height=1.7in]{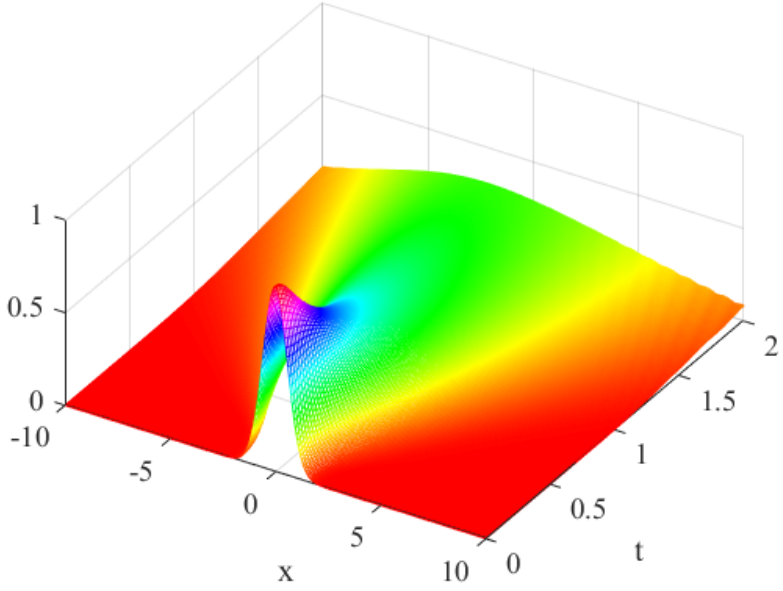}}
  \subfigure[Real part]{\includegraphics[width=2.in,height=1.7in]{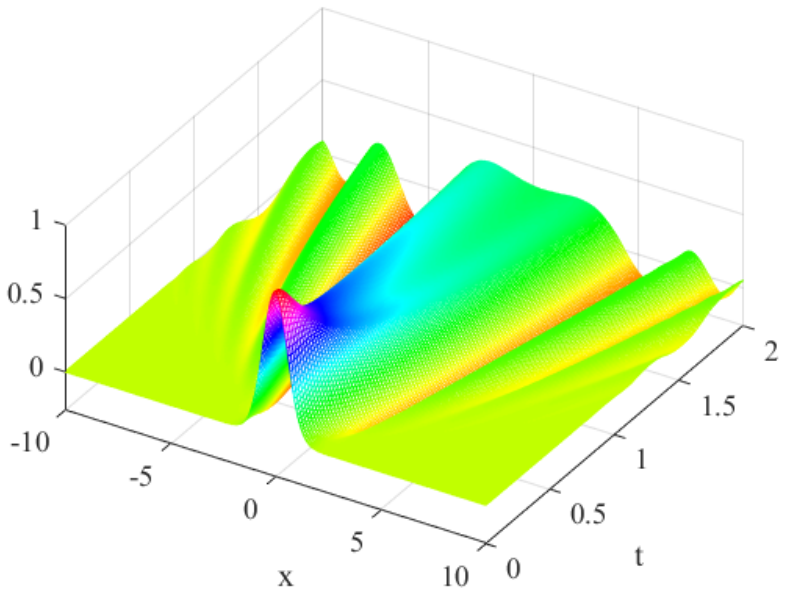}}
  \subfigure[Imaginary part]{\includegraphics[width=2.in,height=1.7in]{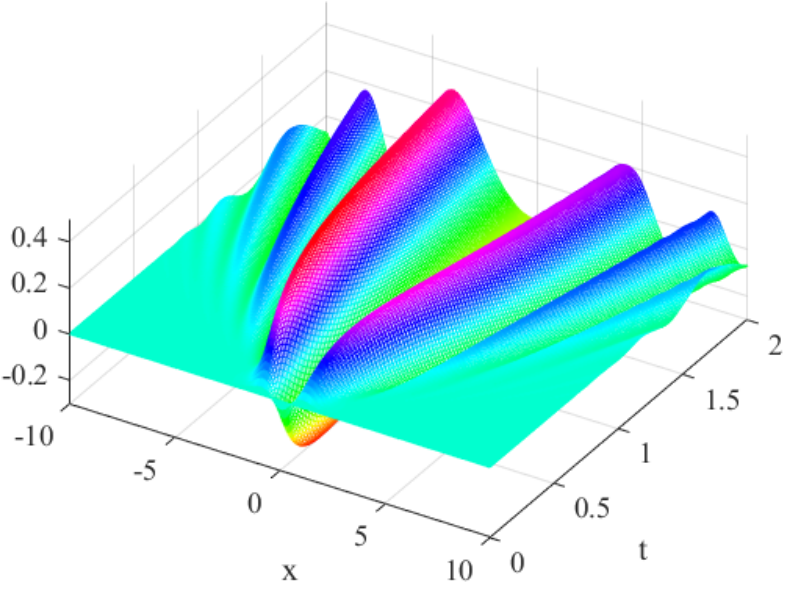}}
  \caption{The spatiotemporal evolution of the $\exp(-x^2){\rm sech}(x)$ initial profile calculated by the NIST in the short time.}\label{small_time_1}
\end{figure}

In the short time, we don't need to deform the Riemann-Hilbert problem.
The evolution of $\exp(-x^2){\rm sech}(x)$ initial profile calculated by the NIST is shown in Figure \ref{small_time_1}.
Figure \ref{small_time_1}(a) shows the density for the calculated result;
Figure \ref{small_time_1}(b) shows the real part for the calculated result;
Figure \ref{small_time_1}(c) shows the imaginary part for the calculated result.
From Figure \ref{small_time_1}, we can clearly see the spatiotemporal evolution of the initial value in the short time.

\begin{figure}[h]
  \centering
  \subfigure[$t=0$]{\includegraphics[width=1.6in,height=1.35in]{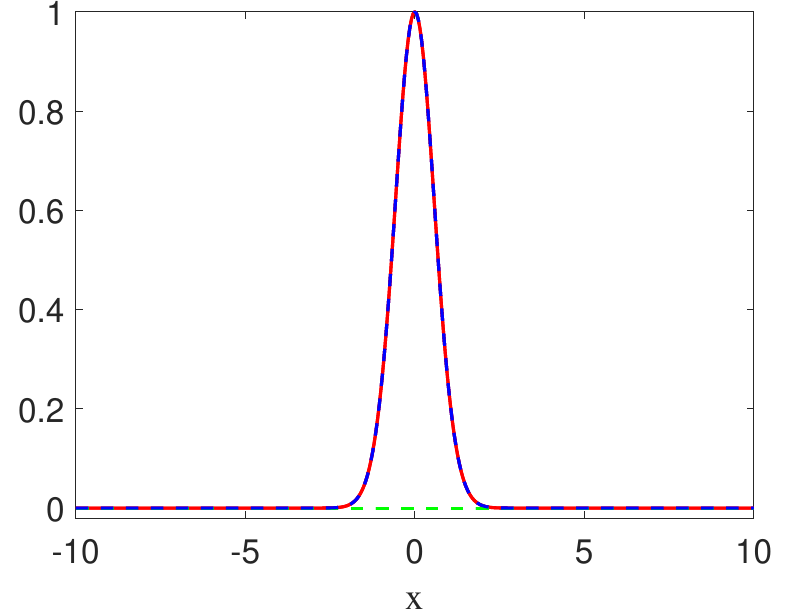}}
  \subfigure[$t=1$]{\includegraphics[width=1.6in,height=1.35in]{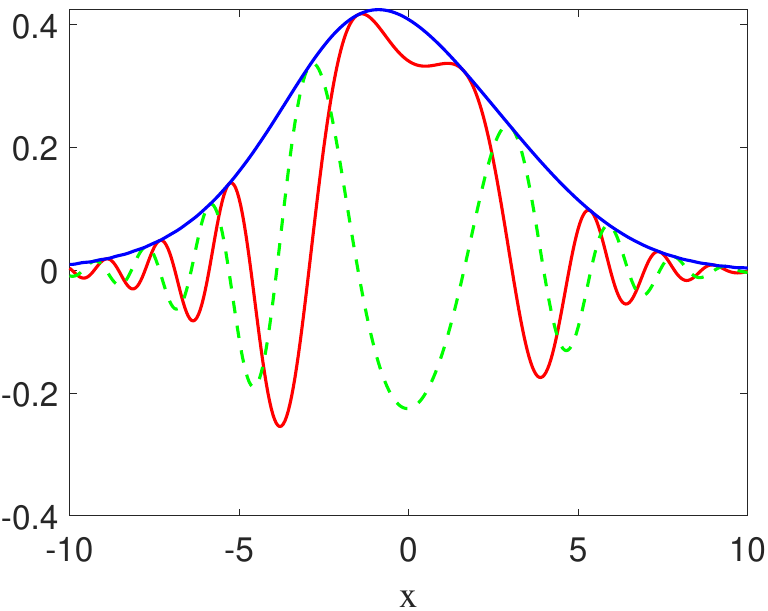}}
  \subfigure[$t=2$]{\includegraphics[width=1.6in,height=1.35in]{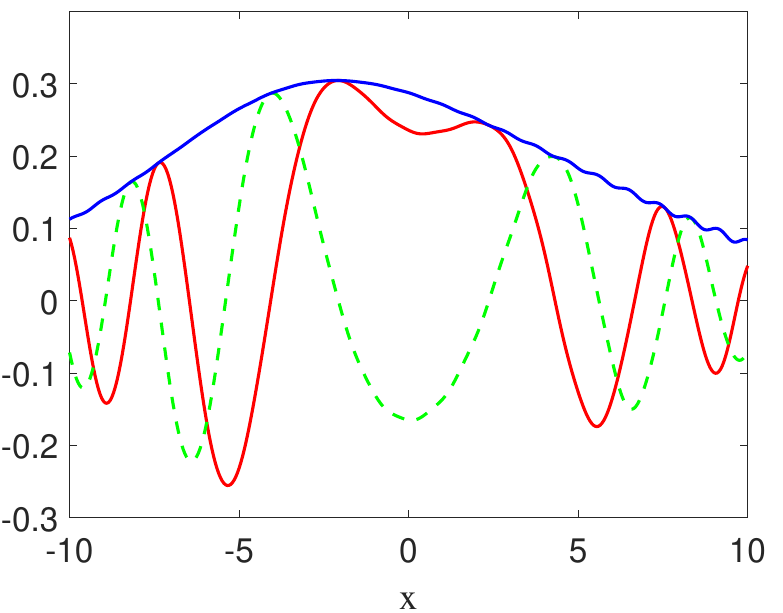}}
  \caption{The numerical solution $q_{\rm nist1}(x,t)$ for $\exp(-x^2){\rm sech}(x)$ initial profile calculated by the NIST. Blue line: the density for $q_{\rm nist1}(x,t)$. Red line: the real part for $q_{\rm nist1}(x,t)$. Green line: the imaginary part for $q_{\rm nist1}(x,t)$.}\label{small_time_2}
\end{figure}

\begin{figure}[H]
  \centering
  \subfigure[$t=0.5$]{\includegraphics[width=1.6in,height=1.35in]{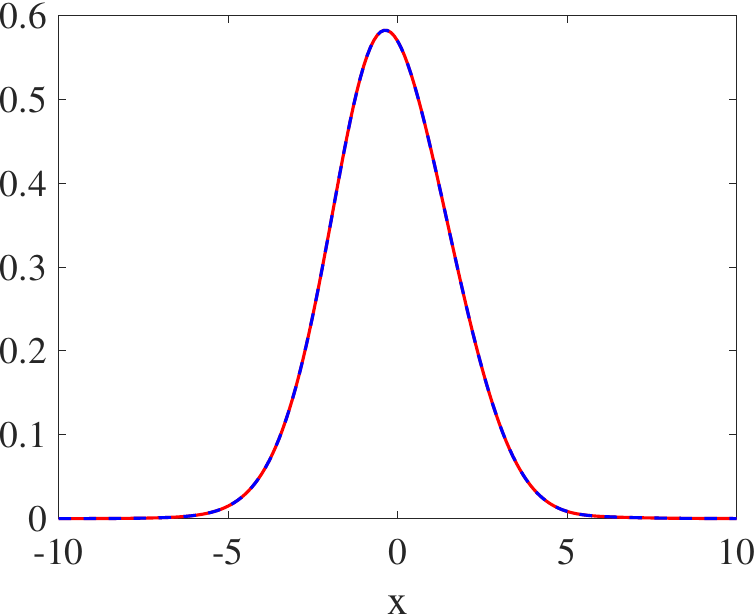}}
  \subfigure[$t=1$]{\includegraphics[width=1.6in,height=1.35in]{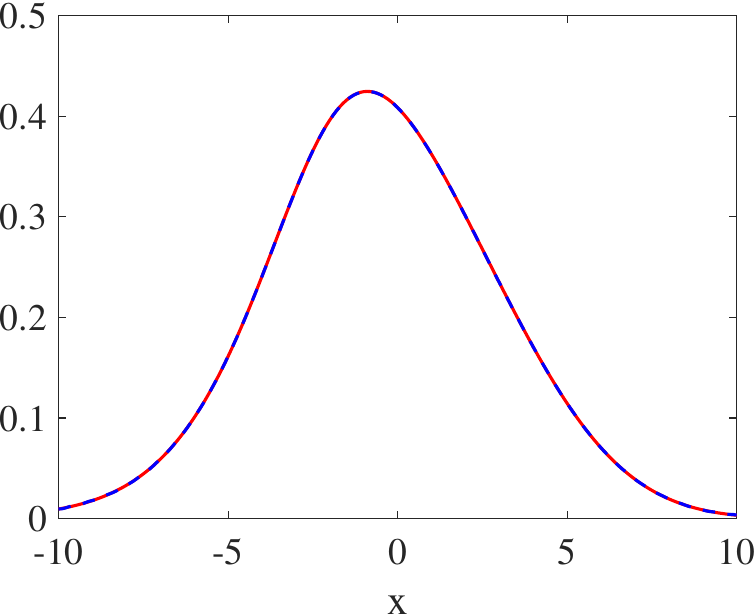}}
  \subfigure[$t=1.5$]{\includegraphics[width=1.6in,height=1.35in]{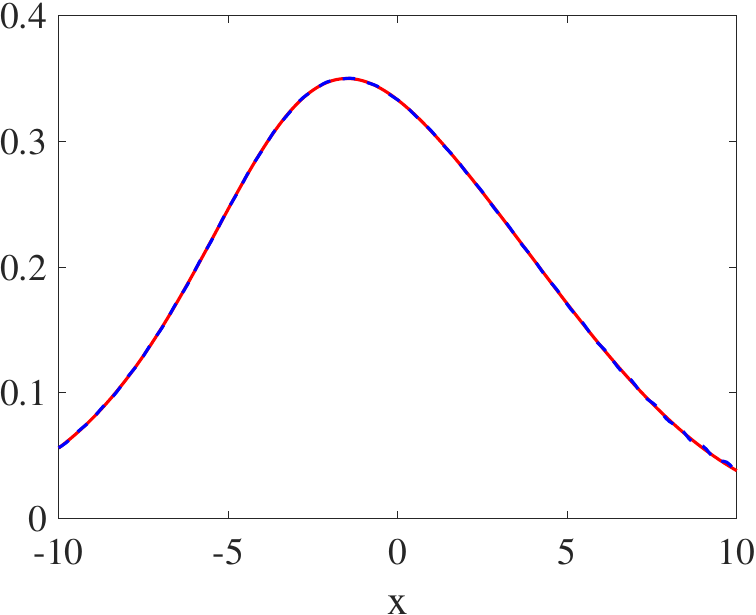}}
  \caption{The density of the numerical solution for $\exp(-x^2){\rm sech}(x)$ initial profile. Red line: the numerical result calculated by the FSM. Green line: the numerical result calculated by the NIST.}\label{small_time_2_1}
\end{figure}

\begin{figure}[H]
  \centering
  \subfigure[$t=0.5$]{\includegraphics[width=1.6in,height=1.35in]{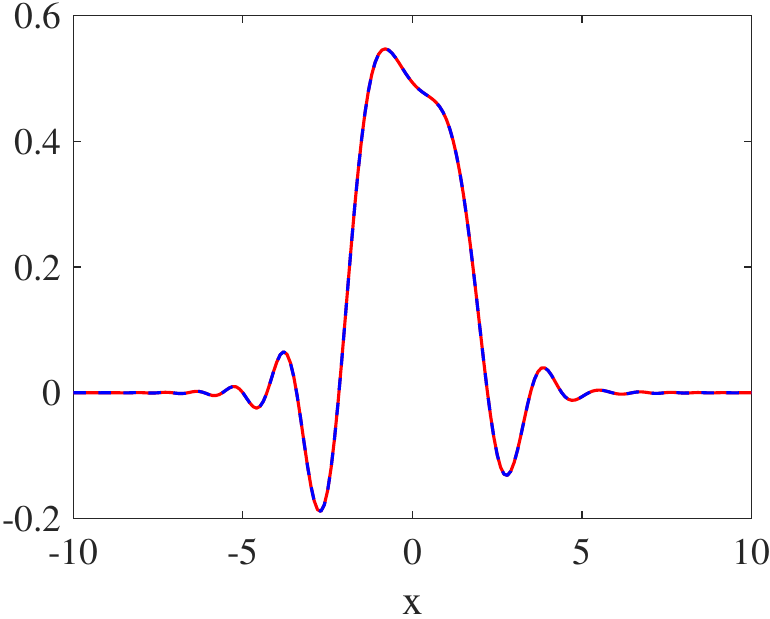}}
  \subfigure[$t=1$]{\includegraphics[width=1.6in,height=1.35in]{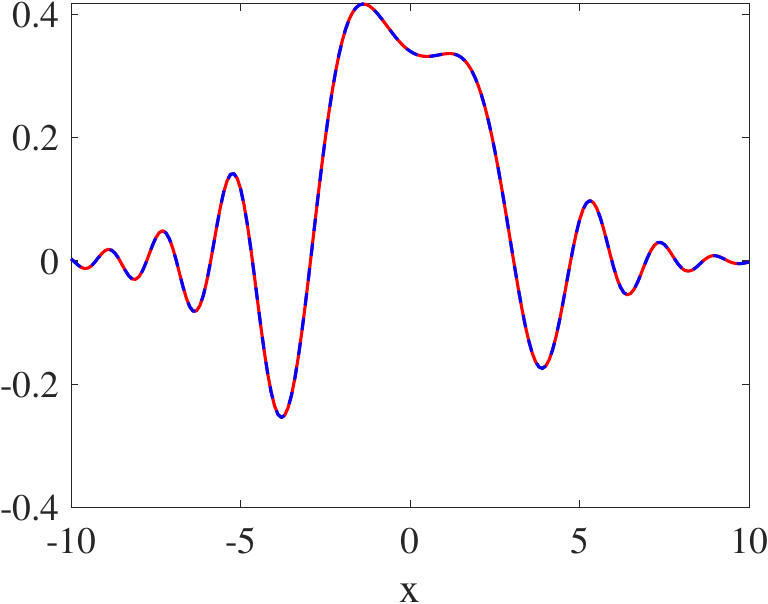}}
  \subfigure[$t=1.5$]{\includegraphics[width=1.6in,height=1.35in]{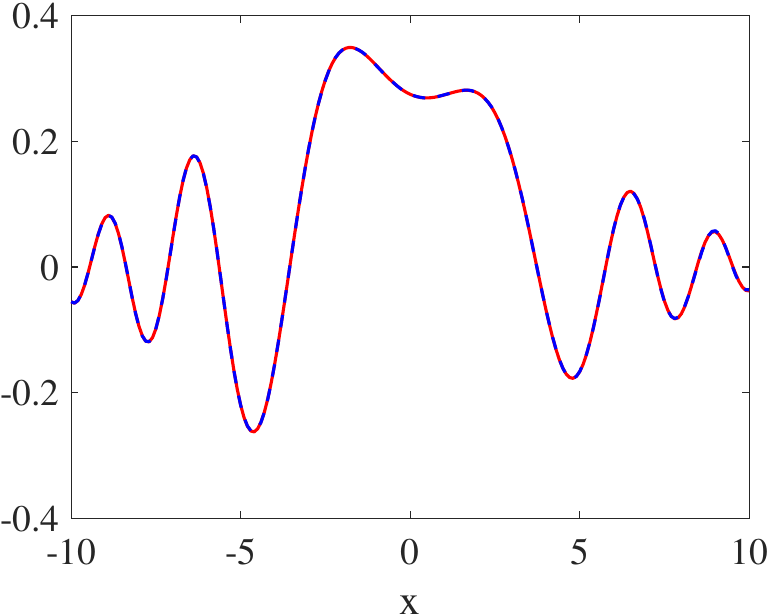}}
  \caption{The real part of the numerical solution for $\exp(-x^2){\rm sech}(x)$ initial profile. Red line: the numerical result calculated by the FSM. Green line: the numerical result calculated by the NIST.}\label{small_time_2_2}
\end{figure}

\begin{figure}[H]
  \centering
  \subfigure[$t=0.5$]{\includegraphics[width=1.6in,height=1.35in]{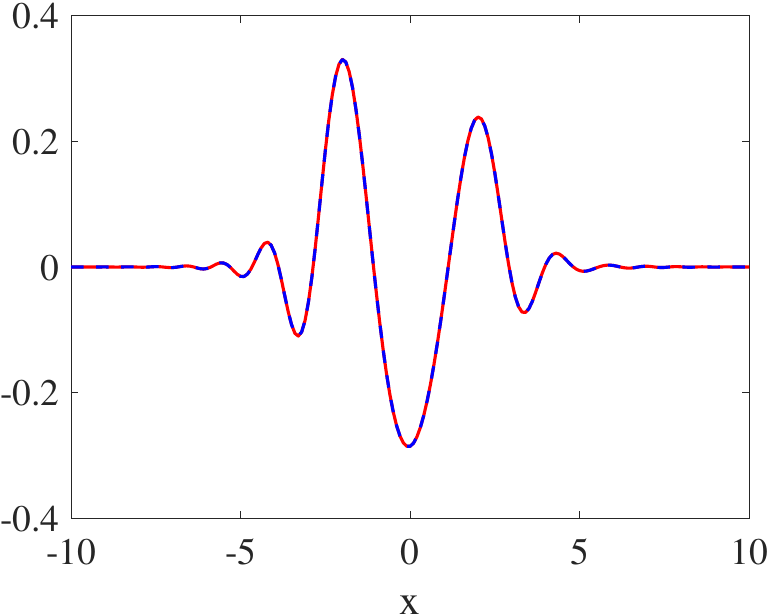}}
  \subfigure[$t=1$]{\includegraphics[width=1.6in,height=1.35in]{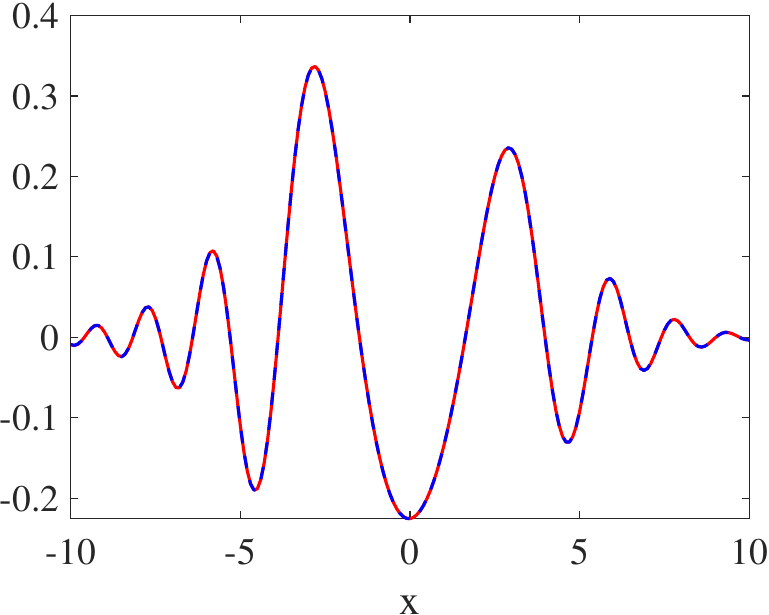}}
  \subfigure[$t=1.5$]{\includegraphics[width=1.6in,height=1.35in]{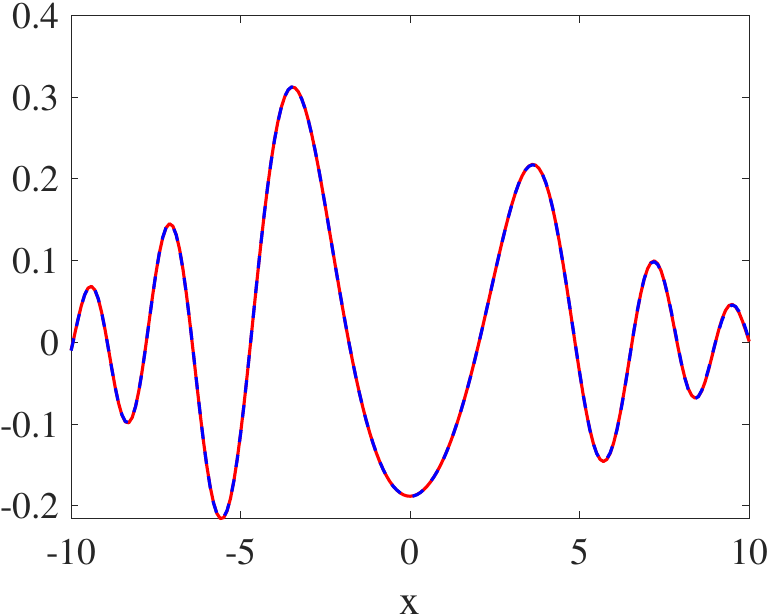}}
  \caption{The imaginary part of the numerical solution for $\exp(-x^2){\rm sech}(x)$ initial profile. Red line: the numerical result calculated by the FSM. Green line: the numerical result calculated by the NIST.}\label{small_time_2_3}
\end{figure}

\begin{figure}[H]
  \centering
  \subfigure[$t=0.5$]{\includegraphics[width=1.6in,height=1.35in]{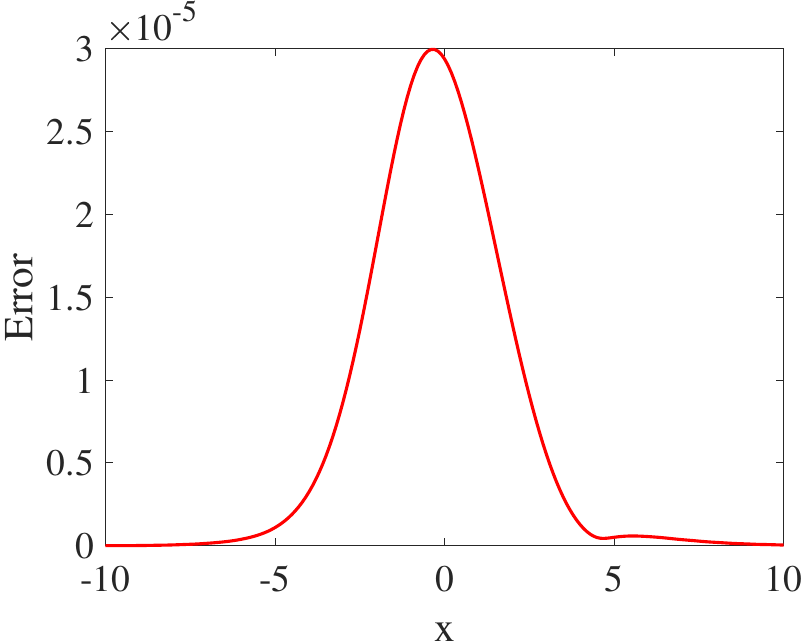}}
  \subfigure[$t=1$]{\includegraphics[width=1.6in,height=1.35in]{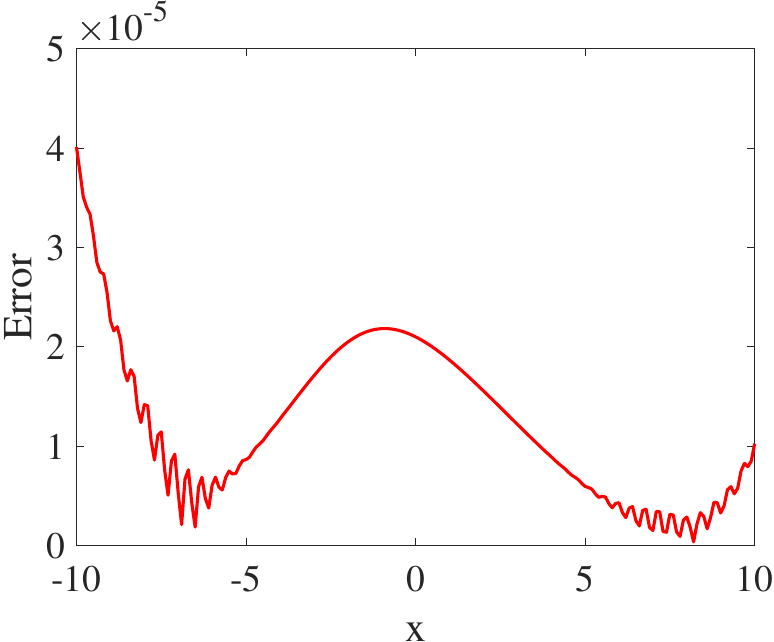}}
  \subfigure[$t=1.5$]{\includegraphics[width=1.6in,height=1.35in]{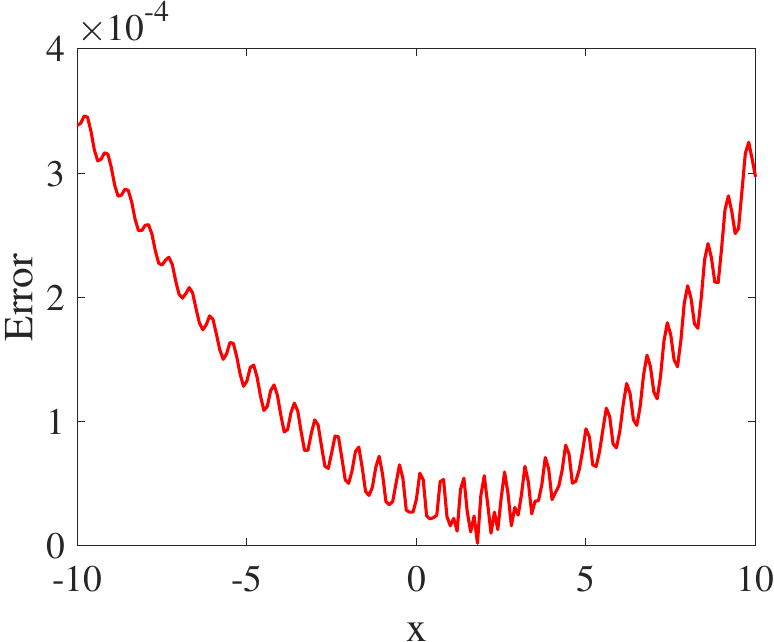}}
  \caption{The absolute errors between the NIST and the FSM at different time.}\label{small_time_2_4}
\end{figure}

Figure \ref{small_time_2} shows the numerical result $q_{\rm nist1}(x,t)$ for the $\exp(-x^2){\rm sech}(x)$ initial profile at $t=0$, $t=1$ and $t=2$. The blue line is the density for $q_{\rm nist1}(x,t)$, the red line is the real part for $q_{{\rm nist1}}(x,t)$ and the green line is the imaginary part for $q_{\rm nist1}(x,t)$.
When time is small, the calculation results of these two methods match theoretically.
To verify this conclusion, we used two methods to solve the numerical solution for $\exp(-x^2){\rm sech}(x)$ initial profile in the short time, and the density for the numerical solutions are shown in Figure  \ref{small_time_2_1}, Figure \ref{small_time_2_2} and Figure \ref{small_time_2_3}.
Figure \ref{small_time_2_1} shows the density of the numerical solutions for $\exp(-x^2){\rm sech}(x)$ initial profile.
Figure \ref{small_time_2_2} shows the real part of the numerical solutions for $\exp(-x^2){\rm sech}(x)$ initial profile.
Figure \ref{small_time_2_3} shows the imaginary part of the numerical solutions for $\exp(-x^2){\rm sech}(x)$ initial profile. The absolute errors calculated by the two methods are shown in Figure \ref{small_time_2_4}, and `Error' is defined by $${\rm Error}=\big|q_{\rm fsm1}(x,t)-q_{\rm nist1}(x,t)\big|,$$
where $q_{\rm fsm1}(x,t)$ is the numerical solution for $\exp(-x^2){\rm sech}(x)$ initial profile calculated by the FSM.
It is obvious that the two method are highly consistent, so both of them are effective.

\begin{figure}[H]
  \centering
  \includegraphics[width=1.6in,height=1.35in]{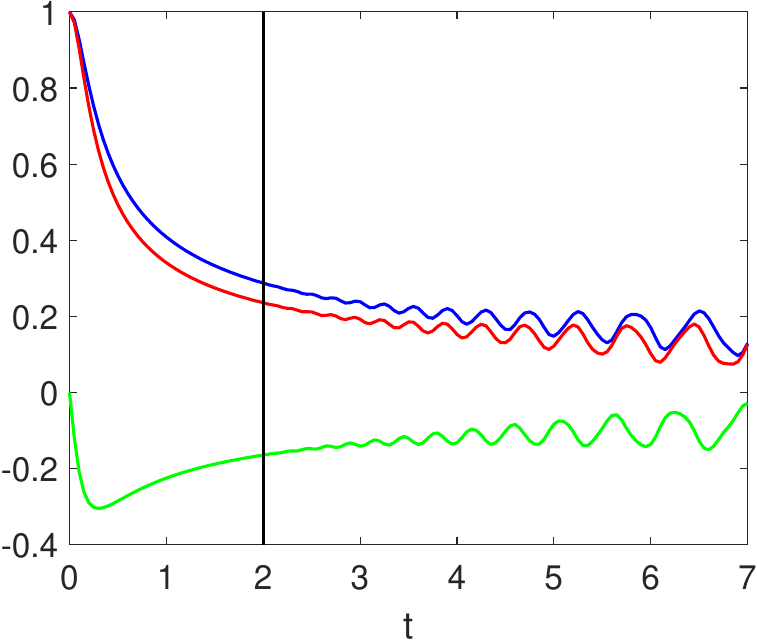}
  \includegraphics[width=1.6in,height=1.35in]{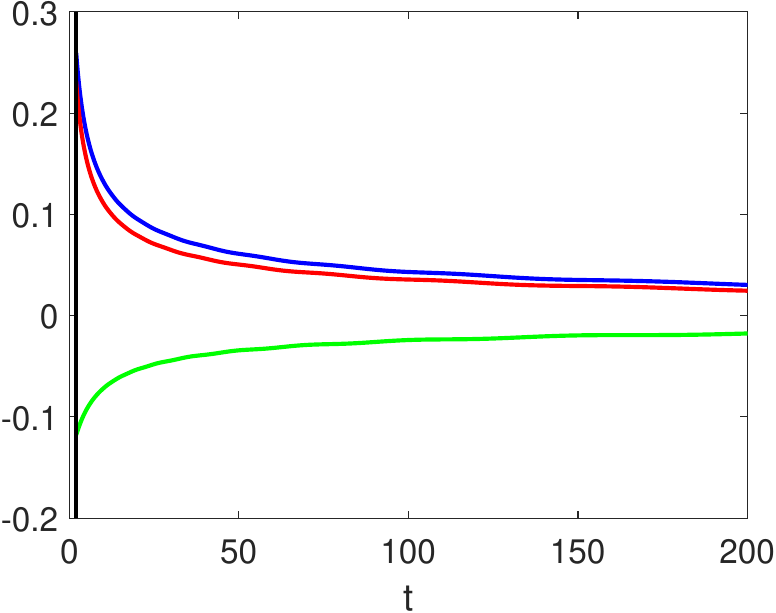}
  \caption{Subfigure(a): the numerical results of $q_{\rm nist1}(0,t)$ initial profile without deformation. Subfigure(b): the numerical results for $q_{\rm nist1}(0,t)$ initial profile after deformation. Blue line: the density for $q_{\rm nist1}(0,t)$. Red line: the real part for $q_{\rm nist1}(0,t)$. Green line: the imaginary part for $q_{\rm nist1}(0,t)$. Black line: $t=2$.}\label{small_time_3}
\end{figure}

The NIST is used to calculate $q_{\rm nist1}(0,t)$ in $t\in[0,7]$, the calculated result is shown in Figure \ref{small_time_3}(a).
When $t>2$, $q_{\rm nist1}(0,t)$ changes from flat to oscillating.
As time increases, the oscillation for $q_{\rm nist1}(0,t)$ becomes more severe.
We have mentioned the reasons for the above phenomenon in section \ref{sec_inverse}.
As time increases, the oscillation of the jump matrix becomes more severe. For highly oscillatory functions, it is difficult to achieve their uniform approximation, and we can not solve the original Riemann-Hilbert problem with high accuracy.
So when $t\geq2$, we will deform the original Riemann-Hilbert problem to ensure numerical accuracy.
After deformation, the NIST is used to calculate $q_{\rm nist1}(0,t)$, the calculated result is shown in Figure \ref{small_time_3}(b). After deformation, the oscillation disappears and the NIST becomes effective.

\begin{figure}[H]
  \centering
  \subfigure[$t=100$]{\includegraphics[width=2in,height=1.5in]{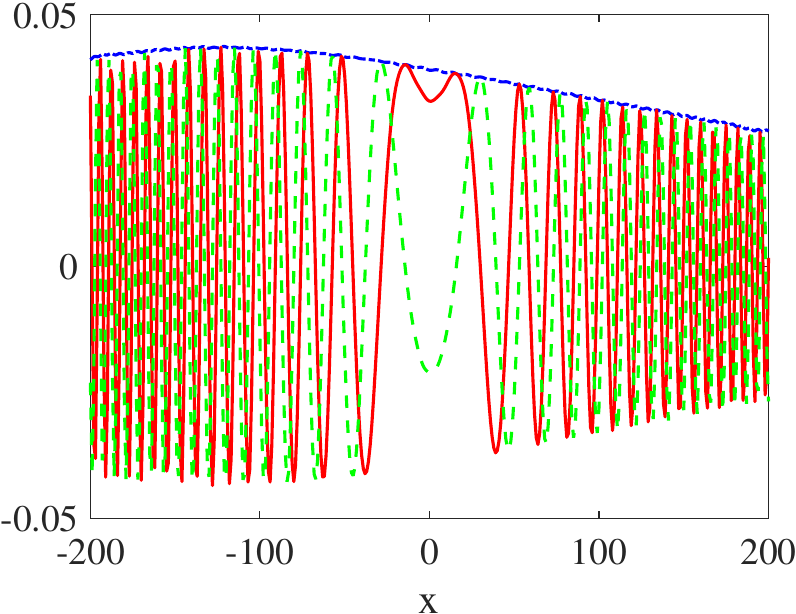}}
  \subfigure[$t=200$]{\includegraphics[width=2in,height=1.5in]{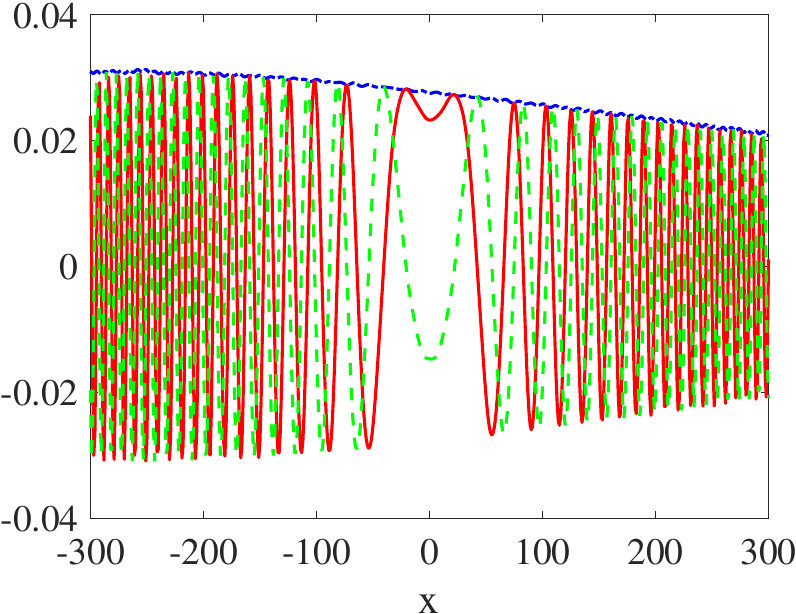}}
  \subfigure[$t=300$]{\includegraphics[width=2in,height=1.5in]{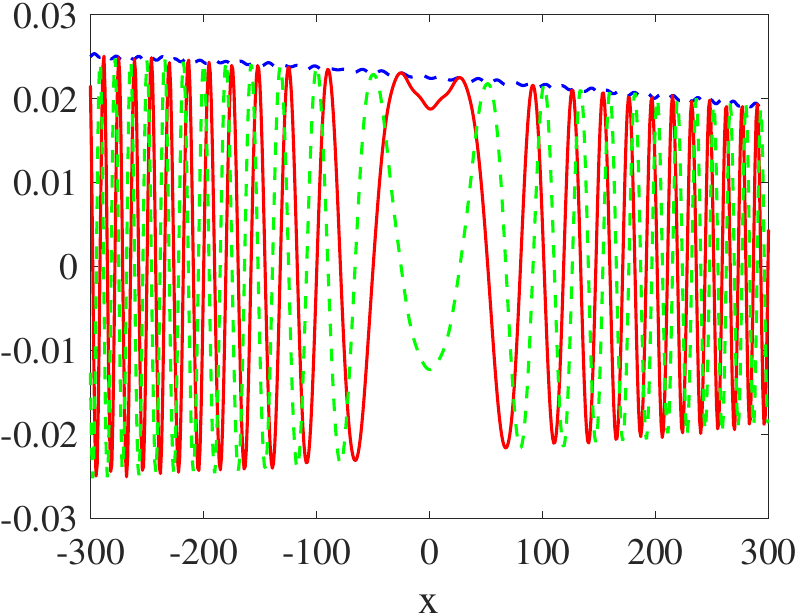}}
  \caption{The long-time evolution $q_{\rm nist1}(0,t)$ for $\exp(-x^2){\rm sech}(x)$ initial profile calculated by the NIST. Blue line: the density for $q_{\rm nist1}(0,t)$. Red line: the real part for $q_{\rm nist1}(0,t)$. Green line: the imaginary part for $q_{\rm nist1}(0,t)$.}\label{long_time_1}
\end{figure}

The NIST (Nonlinear Inverse Scattering Transform) method is utilized to calculate the long-term evolution of the DNLS equation. The calculated results are visualized in Figure \ref{long_time_1}. The reliability of these calculations highlights the effectiveness of the NIST approach. In comparison, traditional numerical methods like the finite difference method and the FSM struggle to achieve high-precision numerical calculations over extended periods.


\section{Conclusions and Discussions}\label{sec_con}



In summary, we have successfully demonstrated the effectiveness of the numerical inverse scattering transform (NIST) for the DNLS equation with zero-boundary condition. Through numerical examples, we have validated the method and compared it with the FSM.
The DNLS equation exhibits a continuous spectrum of $i\mathbb{R}\cup\mathbb{R}$, with three saddle points excluding $x=0$. These characteristics posed challenges in applying the numerical inverse scattering transform, but we have successfully addressed them. By dividing the $(x,t)$-plane into three regions, we devised specific approaches to deform the jump contours accordingly.
To ensure the accuracy of our calculations, we have mathematically established that the deformed jump matrices remain uniformly bounded. This provides a guarantee for the reliability of our results. Moreover, we have established a feasible framework for computing the scattering data of the DNLS system.

During our calculations, we observed that the NIST method imposes strict limitations on the space of initial values. Specifically, the initial value must belong to the Schwartz space $\mathcal{S}(\mathbb{R})$ or even spaces with faster decay. This restriction narrows down the range of valid initial values for the NIST.
A significant focus of future research will involve expanding the space of initial values to encompass a more general class within the NIST framework. By reducing the NIST's dependence on the decay rate of the initial value, we can enhance its applicability and accommodate a wider range of initial conditions.
Addressing and overcoming these limitations will be a key area of focus in our future work. Such advancements will help advance the capabilities of the numerical inverse scattering transform for the DNLS equation and broaden its range of applications.
\section*{Acknowledgment}
\addcontentsline{toc}{section}{Acknowledgment}
This work is supported by National Science Foundation of China (52171251, U2106225, 52231011) and Dalian Science and Technology Innovation Fund (2022JJ12GX036).
\section*{Declaration of competing interest}
The authors declare that they have no known competing financial interests or personal relationships that could have appeared to
influence the work reported in this paper.
\section*{Data availability}
Data will be made available on request.

\end{document}